\documentclass{amsart}
\usepackage{amsmath}%
\usepackage{amsthm}
\usepackage{amsfonts}%
\usepackage{amssymb}%
\usepackage{graphicx}
\usepackage{tikz-cd}
\usetikzlibrary{cd}
\usepackage{bm}
\usepackage{comment}
\usepackage{hyperref}
\usepackage{todonotes}

\usepackage{mathrsfs}

\DeclareFontFamily{OT1}{pzc}{}
\DeclareFontShape{OT1}{pzc}{m}{it}{<-> s * [1.10] pzcmi7t}{}
\DeclareMathAlphabet{\mathpzc}{OT1}{pzc}{m}{it}

\newtheorem{theorem}{Theorem}[section]

\newtheorem{definition}[theorem]{Definition}

\newtheorem{lemma}[theorem]{Lemma}

\newtheorem{proposition}[theorem]{Proposition}
\newtheorem{remark}[theorem]{Remark}

\numberwithin{equation}{section}

\def\XXint#1#2#3{{\setbox0=\hbox{$#1{#2#3}{\int}$}
\vcenter{\hbox{$#2#3$}}\kern-.5\wd0}}

\newcommand{\R}{\mathbb{R}}
\renewcommand{\P}{\mathbb{P}}
\newcommand{\Q}{\mathbb{Q}}
\newcommand{\N}{\mathbb{N}}

\newcommand{\M}{\mathcal{M}}
\newcommand{\F}{\mathcal{F}}
\newcommand{\Ha}{\mathcal{H}}
\newcommand{\leb}{\mathcal{L}}

\newcommand{\essinf}{\operatorname{essinf}}
\newcommand{\esssup}{\operatorname{esssup}}

\newcommand{\spt}{\operatorname{spt}}

\newcommand{\Mod}{\operatorname{Mod}}
\newcommand{\dist}{\operatorname{dist}}
\newcommand{\diam}{\operatorname{diam}}

\newcommand{\biLip}{\operatorname{biLip}}
\newcommand{\Lip}{\operatorname{Lip}}
\newcommand{\LIP}{\operatorname{LIP}}
\newcommand{\im}{\operatorname{im}}
\newcommand{\gap}{\operatorname{gap}}
\newcommand{\spanv}{\operatorname{span}}
\newcommand{\dom}{\operatorname{dom}}
\newcommand{\Fr}{\operatorname{Fr}}

\newcommand{\conv}{\operatorname{conv}}
\newcommand{\sh}{\operatorname{sh}}

\newcommand{\ud}{\mathrm {d}}

\newcommand{\inv}{^{-1}}

\newcommand{\defl}{\mathrel{\mathop:}=}
\newcommand{\defr}{=\mathrel{\mathop:}}

\newcommand{\loc}{\mathrm{loc}}

\makeatletter
\newcommand*{\cone}{%
	{%
		\mathpalette\@coneOf{\times}%
	}%
}
\newcommand*{\@coneOf}[2]{%
	\sbox0{$\m@th#1\mathsf{#2}$}%
	\mathsf{#2}%
	\kern-\wd0 %
	\mkern2.00mu\relax
	\nonscript\mkern0.25mu\relax
	\mathsf{#2}%
}
\makeatother

\title{Fragment-wise differentiable structures}

\author{David Bate}
\address{Zeeman Building, University of Warwick, Coventry CV4 7AL, UK}
\email{david.bate@warwick.ac.uk}

\author{Sylvester Eriksson-Bique}
\address{Department of Mathematics and Statistics, University of Jyv\"askyl\"a,
Seminaarinkatu 15, PO Box 35,  FI-40014 University of Jyv\"askyl\"a, Finland}
\email{sylvester.d.eriksson-bique@jyu.fi}

\author{Elefterios Soultanis}
\address{Department of Mathematics and Statistics, University of Jyv\"askyl\"a,
Seminaarinkatu 15, PO Box 35,  FI-40014 University of Jyv\"askyl\"a, Finland}
\email{elefterios.e.soultanis@jyu.fi}

\thanks{The first and third authors were supported by the European Union's Horizon 2020 research and innovation programme grant number 948021. The second author was supported by the Research Council of Finland grants 354241 and 356861. The third author was also supported by the Research Council of Finland grant 355122. The research was started during a visit of the second named author to University of Warwick. He thanks the institution for its hospitality. %
For the purpose of open access, the authors have applied a Creative Commons Attribution (CC BY) licence to any Author Accepted Manuscript version arising from this submission.}
\subjclass[2020]{30L99, 51F30, 28A75}

\date{\today}

\begin{document}

\begin{abstract}
    We give a definition of a fragment-wise chart and a fragment-wise differential of a Lipschitz function in general complete metric measure spaces, and show that they exist whenever the space has finite Hausdorff dimension - without any further topological or analytic assumptions. The differential has a natural definition which extends the definition of a differential in Lipschitz differentiability spaces. We also show a cornerstone in the theory of tangent modules: density of directions of curve fragments in the associated tangent space. 

    The key to our approach is the combination of modulus and upper gradient-based methods that draws analogies between the Sobolev and Lipschitz theories. These concepts bridge approaches based on abstract modules, PI-spaces,  decomposability bundles and differentiability spaces. In this unified perspective, we give a sharp approximation result,  develop duality between $\infty$-modulus and plans and show the existence of a minimal $*$-upper gradient.
\end{abstract}

\maketitle

\tableofcontents

\section{Introduction}

\subsection{Overview} 

Analysis on metric spaces, pionereed e.g. in \cite{che99,hei98} at the turn of the millenium, has become ubiquitous in the last decades. It is based on the deep observation of the crucial role {curves} in the ambient space play in analysis. \emph{Upper gradients} -- introduced and developed through the works of Heinonen--Koskela, Shanmugalingam, Cheeger and others \cite{hei98, che99, sha00, hei01} -- capture this idea by controlling the oscillation of functions along curves, leading to a fruitful theory of Sobolev spaces over metric measure spaces. Given a metric measure space $(X,d,\mu)$, a Borel function $\rho:X\to [0,\infty]$ is called an upper gradient if of $f:X\to \R$, if 
\begin{align}\label{eq:ug}
    \operatorname{osc}_\gamma f\le \int_\gamma \rho\ud s
\end{align}
for all curves $\gamma$ in $X$ (see also Definition \ref{def:*-ug} below).

In \cite{che99}, Cheeger showed that, in spaces satisfying a \emph{Poincar\'e inequality}, the theory of upper gradients leads to a weak differentiable structure, where every Lipschitz function has a differential \cite{che99,kei02,keith04}. Poincar\'e inequalities are expressed in terms of upper gradients and mean value integrals, and understanding their relationship to differentiability led to the development of \emph{Lipschitz differentiability spaces} (LDS) \cite{bate15}, a strictly weaker condition, which is closely related to (independent) \emph{Alberti representations}; these are decompositions of the measure into \emph{curve fragments} providing (independent) \emph{directions} in the space.

\medskip\noindent In this work, we take as the starting point analogues of upper gradients suited for curve fragments -- so-called $\ast$-upper gradients introduced in \cite{bate18} -- and show that, together with \emph{fragment-wise charts} (see Definition \ref{def:frag-chart} below), they lead to the existence of a differential of Lipschitz functions, along with a pointwise norm recovering almost everywhere the minimal $\ast$-upper gradient of a function, whenever the underlying space has finite Hausdorff dimension. The differential satisfies the usual calculus rules and, moreover, extends the differential defined in Lipschitz differentiability spaces. 

Our approach gives a unifying perspective on the Sobolev and Lipschitz theories. Among its applications is a new sharp approximation result (Theorem \ref{thm:approximation}), the density of directions of curve fragments at generic points (Theorem \ref{thm:density-of-directions-intro}), and a description of the module of Weaver derivations as $L^\infty$-vector fields with respect to the differential structure arising from Theorem \ref{thm:diff-and-ptwise-norm} (Theorem \ref{thm:weaver-der}).

\subsection{$\ast$-Upper gradients and differential of Lipschitz functions}
Throughout the paper, we consider a metric measure space $X=(X,d,\mu)$, i.e. a complete separable metric space $(X,d)$ equipped with a Borel regular measure $\mu$ that is finite on bounded sets. A curve fragment is a bi-Lipschitz map $K\to X$ from a compact set $K\subset \R$, and we denote their collection by $\Fr(X)$.

\begin{definition}\label{def:*-ug}
A Borel function $\rho:X\to [0,\infty]$ is a called a $\ast$-upper gradient of a Lipschitz function $f:X\to\R$ if 
\begin{align}\label{eq:*-ug}
\operatorname{osc}_\gamma f\le \int_\gamma \rho \ud \Ha^1+\LIP(f)\gap(\gamma)
\end{align}
for every curve fragment $\gamma\in \Fr(X)$.
\end{definition}
Here $\displaystyle \operatorname{osc}_\gamma f=\sup_{s,t\in \dom(\gamma)}|f(\gamma_t)-f(\gamma_s)|$, and $\gap(\gamma)$ is defined in Section \ref{sec:frag+star-ug}. Observe that the $\ast$-upper gradient is defined only for Lipschitz functions. A further difference compared to Sobolev spaces is the following fact: 
\begin{align}\label{eq:min-*-ug}
    \textrm{Given $f\in \LIP(X)$, there exists a \emph{minimal} $\ast$-upper gradient } |Df|_\ast\textrm{ of }f,
\end{align} 
i.e. $|Df|_*\le \rho$ $\mu$-a.e. for any other $\ast$-upper gradient $\rho$ of $f$ (Proposition \ref{prop:min-star-ug-exists}).

Our first theorem establishes the existence of fragment-wise differentials and pointwise norm arising from -- and recovering -- the minimal $\ast$-upper gradient. Given a Borel set $U\subset X$ with $\mu(U)>0$ and $\varphi\in \LIP(X;\R^n)$, a Borel map $\bm\xi:U\to (\R^n)^*$ is a fragment-wise differential of $f$ with respect to $(U,\varphi)$, if there is a $\mu$-null set $N\subset X$ so that 
\begin{align}\label{eq:diff}
(f\circ\gamma)_t'=\bm\xi((\varphi\circ\gamma)'_t)\quad a.e.\ t\in \gamma\inv(U)
\end{align}
for all fragments $\gamma\in \Fr(X)$ for which $|\gamma\inv(N)|=0$. If $\bm \xi=\bm \zeta$ $\mu$-a.e. on $U$ whenever both satisfy \eqref{eq:diff} (for some null sets $N_\xi, N_\zeta\subset X$), we say $f$ has a unique fragment-wise chart with respect to $(U,\varphi)$. If $n=0$, we interpret $\R^0 =\{0\}$ and only consider the map $\varphi \equiv 0$, in which case $\bm\xi\equiv 0$ is a fragment-wise differential if there is a $\mu$-null set $N\subset X$ so that 
\begin{align}\label{eq:diff}
(f\circ\gamma)_t'=0\quad a.e.\ t\in \gamma\inv(U)
\end{align}
for all fragments $\gamma\in \Fr(X)$ for which $|\gamma\inv(N)|=0$.

\begin{theorem}\label{thm:diff-and-ptwise-norm}
Suppose $(U,\varphi)$ is an $n$-dimensional fragment-wise chart of $(X,d,\mu)$. Then every $f\in \LIP(X)$ admits a unique fragment-wise differential $\ud f:U\to(\R^n)^*$. Moreover for $\mu$-a.e. $x\in U$ there exists a norm $|\cdot|_x$ on $(\R^n)^*$ such that $x\mapsto |\xi|_x$ is Borel for every $\xi\in (\R^n)^*$ and
\begin{align}\label{eq:ptwise-norm}
|Df|_\ast(x)=|\ud_xf|_x\quad\mu-a.e.\ x\in U
\end{align}
for every $f\in \LIP(X)$.
\end{theorem}
We refer to Section \ref{sec:basic-prop+exist} for the usual calculus rules the differential satisfies, and the arising fragment-wise differential structure. Here, we define fragment-wise charts in terms of independent Alberti representations (see Definition \ref{def:alberti-rep}).
\begin{definition}\label{def:frag-chart}
    An $n$-dimensional fragment-wise chart $(U,\varphi)$ of $(X,d,\mu)$ consists of a Borel set $U\subset X$ with $\mu(U)>0$ and a Lipschitz function $\varphi:X\to \R^n$ satisfying the following.
    \begin{itemize}
        \item[(i)] {\bf Independence}: $\mu|_U$ admits $n$ $\varphi$-independent Alberti representations;
        \item[(ii)] {\bf Maximality}: If $V\subset U$ and $\mu|_V$ admits $k$ independent Alberti representations, then $k\le n$.
    \end{itemize}
\end{definition}
A $0$-dimensional fragment-wise chart is a pair $(U,\varphi:X\to \{0\})$, where $\varphi\equiv 0$ and which is maximal in the sense that for no positive measure subset $V\subset U$ $\mu|_V$ admits $k$ independent Alberti representations with $k\geq 1$. We remark that, if $\dim_HX<\infty$, then there exist countably many fragment-wise charts (of possibly varying dimension) covering $X$ (Theorem \ref{thm:cover-by-frag-charts}).

\subsubsection*{Relationship with Sobolev spaces and further discussion} While \eqref{eq:ug} naturally leads to Newton--Sobolev spaces $N^{1,p}(X)$ for each exponent $p\in [1,\infty)$, this theory trivializes (in the sense that $N^{1,p}(X)=L^p(X)$) e.g. when $X$ is totally disconnected, while Lipschitz functions can still have non-trivial differentiability properties. (Consider, for example, a fat cantor set $K\subset \R^n$.) In \cite{bate15} the first named author identified curve fragments as the natural replacement for curves and \eqref{eq:*-ug}, with the extra term necessitated by the presence of gaps in curve fragments, was considered in \cite{bate18,bkt19}. In light of Sobolev theory, the existence of a minimal genuine $\ast$-upper gradient \eqref{eq:min-*-ug} is surprising: Sobolev functions $f\in N^{1,p}(X)$ need not have minimal upper gradients - only a minimal $p$-\emph{weak} upper gradient. This difference is small, but arises in many proofs, as a relaxed upper gradient always has to be modified to be a true upper gradient. This difference ultimately stems from the difference between null-families for $\Mod_\infty$ and for $\Mod_p$ for finite $p$, cf. the discussion in Section \ref{sec:duality}.

We also note that in \cite{teriseb} a similar construction of $p$-weak charts and differentials for Sobolev functions along curves appears. The charts and differentials in \cite{teriseb} need not agree with fragment-wise charts and differentials (unless the space supports a Poincar\'e inequality). Nevertheless, these are both instances of a more general underlying principle that an upper gradient-like object (together with dimension bounds on the space) always gives rise to a differential, cf. \cite{eidconj,teriseb,bkt19}. A similar idea has appeared much earlier without a dimension bound in the work of Gigli and Kell that shows how to associate an upper gradient-like object a differential which is an element of a $L^p$-module \cite{gig15, kell} see also \cite{AIDP1, AIDP2} for surveys discussing this point further. A difference in the finite dimensional setting is that one can also associate a chart, similar to the work of Cheeger in \cite{che99}. In the fractal setting, a similar notion in the guise of martingale dimension and index of a Dirichlet form has also been studied, see e.g. \cite{hino} and the references therein.

Compared to the Cheeger differential of $f$, fragment-wise differentials \eqref{eq:diff} give control only along curve fragments and not along arbitrarily changing directions. That the Cheeger differential agrees with the fragment-wise differential on LDS's follows from the equality $|Df|_\ast=\Lip f$ ($\mu$-a.e.) for Lipschitz functions, which in fact characterizes a space being LDS (cf. Theorem \ref{thm:LDS-char}). 

Theorem 1.2 is also related to the decomposability bundle of Alberti-Marchese \cite{almar16} for measures on Euclidean space. A difference in the two theories is that we do not declare the non-differentiability of a Lipschitz function in the directions "outside" of the chart (respectively, bundle), since this concept has no meaning without further ambient structure. However, it is possible to do this given a ``section" of a fragment-wise chart.

Finally, we mention that fragment-wise charts have been implicit in the literature since the appearance of Alberti representations in \cite{bate15}, and there are many equivalent ways to define them, see e.g. \cite{sch16b,bate15} and the references therein. Fragment-wise charts are related to \cite[Theorem 2.1]{batepure} and \cite[Theorem 1.6 and 1.7]{schulli}. The key difference to many definitions in the literature is that we are not in a Lipschitz differentiability space, or a subset of Euclidean space. 

\subsection{Duality and representation of the minimal $\ast$-upper gradient}\label{sec:duality} A novel aspect of our work is connecting Alberti representations and differentiability of Lipschitz functions to $\ast$-upper gradients. This connection relies on a duality between $\Mod_\infty$ and 1-plans, and yields a representation formula (see Theorem \ref{thm:min_star_ug}) for the minimal $\ast$-upper gradient of a Lipschitz function $f\in \LIP(X)$: there are probability measures $\P_x$ on $\Fr(X)$ for $\mu$-a.e. $x\in \{|Df|_\ast>0\}$ such that
\begin{align}\label{eq:min-star-ug-rep-formula}
|Df|_{\ast}(x)=\left\|\frac{(f\circ\gamma)'(0)}{|\gamma_0'|}\right\|_{L^\infty(\P_x)}\quad \mu-\textrm{a.e. }x\in \{|Df|_\ast>0\}.
\end{align}
The challenge in proving \eqref{eq:min-star-ug-rep-formula} is that -- unlike the duality of $\Mod_p$ and $p'$-plans for finite $p$ \cite{amb13}-- the duality of $\Mod_\infty$ and 1-plans is not known, neither for curve families nor for families of curve fragments. We instead invoke a weaker form of duality (Lemma \ref{lem:rainwater}), a Rainwater-type result \cite{almar16,bate15,sch16b}, which yields plans on sub-fragments, and have to use further arguments to pass from plans to Alberti representations.

\subsection{Density of directions} 
A core result in Lipschitz differentiability spaces is the density of directions of curve fragments,  \cite[Theorem 5.3]{che-kle-sch16}, cf. \cite[Theorem 9.5]{bate15}. This result implies structure for the Gromov-Hausdorff blow-ups of the space, and allows one to obtain surjectivity of blown-up Lipschitz maps, see \cite{che-kle-sch16} for more discussion. We establish a stronger version of this result here without any assumption on the metric measure space.

\begin{theorem}\label{thm:density-of-directions-intro}
Let $(U,\varphi)$ be a fragment-wise chart of dimension $n\in\N$ in a metric measure space $(X,d,\mu)$. For $\mu$-a.e. $x\in U$, the set of directions
\begin{align*}
S_x(\varphi)=\bigg\{ \frac{(\varphi\circ\gamma)'(t)}{|(\varphi\circ\gamma)'(t)|}:\ \gamma_t=x\bigg\}\subset S^{n-1}
\end{align*}
is dense in $S^{n-1}$.
\end{theorem}
This theorem is proven in a stronger form as Theorem \ref{thm:density_of_directions}. The key idea for its proof comes from approximations.

\subsection{Approximation} 
Both the theory of Sobolev spaces and that of Weaver differentials rests heavily on finding good approximations of functions. Indeed, for us, the proof of the density of directions, Theorem \ref{thm:density_of_directions}, is based on the following approximation result.
\begin{theorem}\label{thm:approximation}
Let $f\in \LIP_{bs}(X)$. Then there exists a sequence $(f_j)\subset \LIP_{bs}(X)$ with $\LIP(f_j)\le \LIP(f)$ and $|f_j|\le |f|$ for each $j\in\N$, such that
\begin{itemize}
    \item[(1)] $f_j\to f$ pointwise everywhere (and consequently uniformly on compact sets);
    \item[(2)] $\Lip_af_j\to |Df|_\ast$, $\Lip f_j\to |Df|_\ast$, and $|Df_j|_\ast\to |Df|_\ast$ pointwise $\mu$-a.e.
\end{itemize}
\end{theorem}
Theorem \ref{thm:approximation} is the Lipschitz analogue of the result in \cite{amb13}, although one should note that the proof method here is very different, and is of independent interest. We adapt to the Lipschitz setting the argument from \cite{seb2020}. It gives slightly sharper conclusions than its many predecessors in different contexts, see e.g. \cite[Lemma 6.3]{bate15}, \cite[Theorem 1.21]{sch16b}, \cite[Sections 4 and 5]{che99}  and \cite{almar16}. The proof method in \cite{seb2020} comes from adapting the methods in these references. A similar approximation was also used in \cite{eidconj} to resolve a conjecture of Bouleau and Hirsch. 

\subsection{Lipschitz differentiability spaces and Weaver derivations}

A countable covering of $X$ by fragment-wise charts gives rise to a differentiable structure as well as \emph{tangent} and \emph{cotangent} bundles $TX$ and $T^*X$, respectively, see Section \ref{sec:diff-str}. The differential of $f\in \LIP(X)$ given by Theorem \ref{thm:diff-and-ptwise-norm}  naturally defines an element $\ud f\in L^\infty(T^*X)$ in the space of essentially bounded sections of $T^*X$.

An important property of the linear map $f\mapsto \ud f:\LIP(X)\to L^\infty(T^*X)$ is its weak$^*$-continuity (cf. Theorem \ref{thm:diff-w-star-cont} and Remark \ref{rmk:normed module}), which hints at a connection with \emph{Weaver derivations} (see Section \ref{sec:weaver-der}). Introduced originally by Weaver \cite{wea99}, Weaver derivations have inspired the theory of tangent modules in \cite{gig15}, and turn out to be closely connected with Alberti representations and metric 1-currents \cite{sch16,sch16b}. The space $\mathcal X(\mu)$ of Weaver derivations equipped with the norm $\|\bm b\|_{\mathcal X(\mu)}=\||\bm b|_{\mathcal X,\loc}\|_{L^\infty(\mu)}$ is a Banach space. See Section \ref{sec:weaver-der}, and \eqref{eq:loc-norm} for the definition of the pointwise norm $|\bm b|_{\mathcal X,\loc}$ of a Weaver derivation $\bm b\in\mathcal X(\mu)$.  

In Theorem \ref{thm:weaver-der} below, we identify Weaver derivations with \emph{vector fields} on $X$, i.e. sections of the tangent bundle $TX$. This is used in Theorem \ref{thm:LDS-char} to characterize Lipschitz differentiability spaces in the spirit of \cite{sch16b} by using the minimal $\ast$-upper gradients, respectively. A fragment-wise differential structure is a countable system of fragment-wise charts $(U_i,\phi_i)$ s.t. $\mu(X\setminus \bigcup_{i\in \N} U_i)=0$; see Subsection \ref{sec:diff-str} for further discussion and the definition of the tangent module $L^\infty(TX)$ and Theorem \ref{thm:cover-by-frag-charts} for a condition for the existence of such a structure.

\begin{theorem}\label{thm:weaver-der}
	Suppose $X=(X,d,\mu)$ admits a fragment-wise differentiable structure. Then the map $\iota:L^\infty(TX)\to \mathcal X(\mu)$ defined by
	\begin{align*}
		\iota(\bm X)(f)\defl \ud f(\bm X),\quad f\in\LIP(X)
	\end{align*}
	is an isomorphism of normed $L^\infty(\mu)$-modules satisfying
	\begin{align}\label{eq:isom}
		|\iota(\bm X)|_{\mathcal X,\loc}=|\bm X| \quad\mu-a.e.
	\end{align}
	for every $\bm X\in L^\infty(TX)$. 
\end{theorem}
Many similar identifications of vector fields and derivations, and different tangent modules have appeared. In the Sobolev setting, with a different definition of tangent module, one was given by Gigli in \cite[Theorem 2.3.3]{gig18}, see also \cite[Theorem 4.2.3]{gigenr}. See also \cite{AIDP2} for further discussion on this identification and different notions of tangent modules. We note that here the $p=\infty$ setting is a key difference between the existing work and ours.

In the claim below, a modulus of continuity is a non-decreasing continuous function $\omega:[0,\infty)\to [0,\infty)$ with $\omega(0)=0$.
\begin{theorem}\label{thm:LDS-char}
	Let $(X,d,\mu)$ be a metric measure space. The following are equivalent.
	\begin{itemize}
		\item[(i)] $X$ is a Lipschitz differentiability space;
		\item[(ii)] $\Lip f=|Df|_{\ast}$ $\mu$-a.e. for every $f\in \LIP(X)$;
		\item[(iii)] there is a collection $\{\omega_x\}_{x\in X}$ of moduli of continuities such that \[\Lip f\le \omega(|Df|_\ast)\quad \mu-a.e.\] for every $f\in \LIP(X)$.
	\end{itemize}
\end{theorem}

We remark that the equivalence (ii)$\iff$(iii) is a self-improving property of the inequality $\Lip f\le \omega(|Df|_\ast)$, cf. the self-improvement of the Lip--lip condition \cite{sch16b,che-kle-sch16}. 

\begin{remark}
    Tangent modules, differential structures and charts have many further applications we do not pursue here. Among these we mention the very active study of metric currents and their structure \cite{sch16, almarflat,wil12}, connections to rectifiability \cite{marchesemerlo,schulli}, closability of operators \cite{bateclosability, pasqualetto20}, and the recent resolution of the energy image density conjecture of Bouleau and Hirsch \cite{eidconj}. 
\end{remark}

\section{Preliminaries}\label{sec:preli}

\subsection{Notation and conventions}

In this paper, we will only consider complete and separable metric measure spaces $X$ that are equipped with  a Radon measure $\mu$ that is finite on bounded sets.

We say that a function $f:X\to Y$ is Lipschitz, if $\LIP(f)=\sup_{x\neq y} \frac{d(f(x),f(y))}{d(x,y)}<\infty$. We denote by $\LIP(X)$ the collection of Lipschitz functions $f:X\to \R$, and treat it simply as a set without a topology. Further, $\LIP(X,Y)$ denotes the space of Lipschitz functions $f:X\to Y$.  We denote by $\LIP_\infty(X)$ the collection of bounded Lipschitz functions $f:X\to \R$, equipped with the norm $\|f\|_{\LIP_\infty}=\max(\|f\|_{\infty}, \LIP(f))$. We let $\LIP_{bs}(X)$ be the space of Lipschitz functions with bounded support, and we write 
\[
\Lip_a f(x)=\lim_{r\to 0} \LIP(f|_{B(x,r)}) \text{ and }  \Lip f(x)=\lim_{r\to 0}\sup_{y\in B(x,r)\setminus \{x\}} \frac{|f(x)-f(y)|}{d(x,y)}.\]
For isolated points $x\in X$, we set $\Lip_a f(x)=\Lip f(x)=0$. 

If $\mu$ is a Radon measure on $X$, then the pushforward of $\mu$ under a measurable map $f:X\to Y$ is denoted $f_*\mu$ and is given by $f_*\mu(A)=\mu(f^{-1}(A))$.
In a few instances, we will wish to use $0$-dimensional charts, and we consider $\R^0=\{0\}$ to be a trivial $0$-dimensional vector space equipped with the trivial norm.

\subsection{Lipschitz free space}\label{sec:lip-free}

Let $\LIP_0(X)$ be the space of Lipschitz functions vanishing at a given base point $x_0\in X$, equipped with the Lipschitz norm $\LIP(\cdot)$. For every $x\in X$ the associated evaluation map $\iota_x(f):=f(x)$ belongs to $\LIP(X)^*$. Note that $\iota_x$ is moreover continuous with respect to pointwise convergence. The Lipschitz free space $\F(X)$ is defined as
\[
\F(X)=\overline\spanv \{\iota_x:\ x\in X\}\subset \LIP_0(X)^*
\]
equipped with the dual norm. The canonical embedding $\iota:X\hookrightarrow \mathcal F(X)$, given by $x\mapsto \iota_x$, is an isometry. The Lipschitz free space $\F(X)$ over $X$ is a predual of $\LIP_0(X)$ -- i.e. $\F(X)^*=\LIP_0(X)$; moreover on bounded sets, the $w^*$-topology on $\LIP_0(X)$ induced by 
$\F(X)$ coincides with the topology of pointwise convergence, cf. \cite[Theorem 2.2.2]{wea99}.  We refer to \cite{wea99} for the more detailed definition of and thorough introduction to $\F(X)$ and its properties. 

We next define the \emph{semihull} $\sh(X)$ of $X$. Set
\begin{align*}
\sh(X)=\bigcup\{[x,y]:\ x,y\in X  \},
\end{align*}
where $[x,y]=\{\iota_x+t(\iota_y-\iota_x):\ t\in [0,1]\}$ is the line segment in $\mathcal F(X)$ connecting the isometric images of $x$ and $y$. Every element $v\in \sh(X)\setminus X$ has two representations $v=\iota_x+t(\iota_y-\iota_x)=\iota_{y}+(1-t)(\iota_{x}-\iota_{y})$ that is, the representation is unique up to switching $x$ and $y$. We equip $\sh(X)$ with the metric from $\mathcal F(X)$. If $v,w\in \sh(X)$ this metric can be written as $d(v,w)=\sup_{\stackrel{f\in \LIP_0(X)}{\LIP[f]\leq 1}} |v(f)-w(f)|$.  It directly follows that $X\subset \sh(X)$ isometrically.

The semihull provides a convenient technical tool which will be used in Section \ref{sec:approximation} to obtain lower semicontinuity of suitable integrals over fragments of curves. Its key property is a canonical extension of functions $X\to \R$ to functions $\sh(X)\to \R$ which preserves lower semicontinuity.

\begin{remark}\label{AE-better}

The standard Kuratowski embedding $X\hookrightarrow \ell^\infty$ is not appropriate for our canonical extension scheme. Indeed, the natural ``convex'' extension we use (see Section \ref{sec:line-integral_lsc}) 
is not well-defined for the Kuratowski embedding, and we do not know whether a canonical semicontinuity preserving extension exists for it. 

\end{remark}

The following lemma gives a concrete description of convergence in $\sh(X)$.
\begin{lemma}\label{lem:sh-conv}
    Let $v=\iota_{x}+t(\iota_{y}-\iota_{x})$ and suppose $(v_j)\subset \sh(X)$ converges to $v$. Then there exist sequences $(x_j),(y_j)\subset X$ and $(t_j)\subset [0,1]$ with $v_j=\iota_{x_j}+t_j(\iota_{y_j}-\iota_{x_j})$ such that the following holds.
    \begin{itemize}
        \item[(1)] If $v\in \sh(X)\setminus\iota(X)$, then $t_j\to t$ and $(x_j,y_j)\to (x,y)$;
        \item[(2)] If $v\in\iota(X)$, then up to a subsequence, we have one of the following options: (a) $(x_j,y_j)\to (x,y)$, (b) $x_j\to x $ and $t_j\to 0$, or (c) $y_j\to y$ and $t_j\to 1$.
    \end{itemize}
\end{lemma}
\begin{proof}
Assume first that $v=\iota_x+t(\iota_y-\iota_x)\notin \iota(X)$, i.e. $x\ne y$ and $t\in (0,1)$. Choose a representation $v_j=\iota_{x_j}+t_j(\iota_{y_j}-\iota_{x_j})$ such that $d(x,x_j)+d(y,y_j)$ is minimal. The functions $f=\dist(\{x,y\},\cdot)$, $\tilde f=f-f(x_0)$ are 1-Lipschitz, $\tilde f$ vanishes at $x_0$, and
$v_j(\tilde f)-v(\tilde f)=t_jf(y_j)+(1-t_j)f(x_j)$. 

We first argue $\liminf_{i\to \infty} \min\{t_j,1-t_j\} >0$. If this were not the case, we could pass to a subsequence such that $t_j\to 0$ or $t_j\to 1$. Assume the first case. Since 
\[
\lim_{j\to \infty} v_j(\tilde f)-v(\tilde f)=0,
\]
and $t_j\to 0$, we get $\lim_{j\to \infty} d(\{x,y\},x_j)=0$, i.e. a subsequence of $x_j$ converges to either $z=x$ or $z=y$. By passing to a subsequence we may assume convergence and consider $g=\dist(z,\cdot)$ and $\tilde g=g-g(x_0)$. We have
$v_j(\tilde g)-v(\tilde g)=t_jf(y_j)+(1-t_j)f(x_j)-((1-t)d(z,x)+td(z,y))$, and
\[
\lim_{j\to \infty} v_j(\tilde g)-v(\tilde g) = -((1-t)d(z,x)+td(z,y))\neq 0,
\]
which is a contradiction. The case $t_j\to 1$ is similar, with $y_j$ converging to either $z=x$ or $z=y$.

Thus, we can assume $\liminf_{i\to \infty} \min\{t_j,1-t_j\} >0$.  We have
    \begin{align*}
    0\le &\min\{t_j,1-t_j\}[\dist(\{x,y\},x_j)+\dist(\{x,y\},y_j)]\\ 
    \le &t_jf(y_j)+(1-t_j)f(x_j) = |v_j(\tilde f)-v(\tilde f)|\le \|v-v_j\|\to 0,
    \end{align*}
    i.e. the sequences $(x_j)$ and $(y_j)$ converge either to $x$ or $y$. By the minimality of $d(x,x_j)+d(y,y_j)$ it follows that $x_j\to x$ and $y_j\to y$. Similarly, if $g=d(x,\cdot)$ and $\tilde g=g-g(x_0)$, then $v_j(\tilde g)-v(\tilde g)=d(x_j,x)+t_j(d(y_j,x)-d(x_j,x))-td(y,x)\to 0$, implying $t_j\to t$.

    Suppose next that $v\in \iota(X)$. Without loss of generality we may assume that $v=\iota_x$ and $t\le 1/2$.  Considering the functions $f=\dist(\{x,y\},\cdot)$ as above we obtain that $x_j\to x$ or $y_j\to x$. Up to relabeling $x_j$ and $y_j$ we may assume that $x_j\to x$. If a subsequence of  $(y_j)$ converges to $x$ then (a) is satisfied. If not, there exists $\delta>0$ such that $d(y_j,x)\ge 2\delta $ for all $j$. Consequently $d(x_j,y_j)\ge \delta$ for all large $j$. It follows that the functions $g_j=\frac{d(x_j,\cdot)}{d(x_j,y_j)}$, $\tilde g_j=g_j-g_j(x_0)$ are uniformly Lipschitz and $v_j(\tilde g_j)-v(\tilde g_j)=t_j-\frac{d(x,x_j)}{d(y_j,x_j)}\to 0$. Since $\frac{d(x,x_j)}{d(y_j,x_j)}\to 0$, it follows that $t_j\to 0$, completing the proof.
\end{proof}

\begin{remark}\label{rmk:edge-length}
Define the \emph{edge-length function} $l_{\sh}:\sh(X)\to [0,\infty)$ by
\begin{align*}
    l_{\sh}(\iota_x+t(\iota_y-\iota_x))\defl d(x,y).
\end{align*}
In particular, $l_{\sh}(x)=0$ for $x\in X\subset \sh(X)$. 
This is well-defined and $\{l_{\sh}>0\}=\sh(X)\setminus X$. It follows from Lemma \ref{lem:sh-conv} that $l_{\sh}$ is lower semicontinuous.
\end{remark}

\subsection{Curve fragments and $\ast$-upper gradients}\label{sec:frag+star-ug}  A curve fragment in $X$ is a bi-Lipschitz map $\gamma:\dom(\gamma)\to X$ from a compact set $\dom(\gamma)\subset \R$, called the domain of $\gamma$. The biLipschitz constant of $\gamma$ is denoted $\biLip(\gamma)$ and it is the smallest constant $L$ so that $L^{-1}|a-b| \leq d(\gamma(a),\gamma(b)) \leq L |a-b|$ holds for all $a,b\in \dom(\gamma)$.

Let $I_\gamma\defl \conv(\dom(\gamma))=[a,b]$ be the smallest interval containing $\operatorname{dom}(\gamma)$ and, writing $\displaystyle [a,b]\setminus \dom(\gamma)=\bigcup_i(a_i,b_i)$, the gap of the curve fragment $\gamma$ is defined as  
\[
\operatorname{gap}(\gamma)= \sum_i d(\gamma_{b_i},\gamma_{a_i}).
\]
We also call the intervals $(a_i,b_i)$ appearing in this context the gaps of the curve fragment $\gamma$.  The gaps are indexed by a finite or countable set. 
We define an extension $\bar\gamma:I_\gamma\to \sh(X)$ of $\gamma$:
\begin{align*}
\bar\gamma_t\defl
\left\{
\begin{array}{ll}
   \gamma_t\ ,  & t\in \dom(\gamma) \\
    \frac{t-a_i}{b_i-a_i}\gamma_{b_i}+\frac{b_i-t}{b_i-a_i}\gamma_{a_i} , & t\in (a_i,b_i).
\end{array}
\right.
\end{align*}
This extension is directly seen to be Lipschitz, and thus absolutely continuous.
The \emph{metric speed} $|\bar\gamma_t'|$ exists for $\leb^1$-a.e. $t\in I_\gamma$, see \cite[Theorem 4.1.6]{ambrosiotilli}.
In particular, curve fragments also admit a metric speed: for $\leb^1$-a.e. $t\in \dom(\gamma)$ the limit
\[
|\gamma_t'|:= \lim_{\dom(\gamma)\ni s\to t}\frac{d(\gamma_t,\gamma_s)}{|t-s|}
\]
exists.
At any non isolated point of $\dom(\gamma)$, this value is independent of the chosen extension.

If $f:X\to \R$ is Lipschitz, then it can be extended to a Lipschitz function $\bar f:\sh(X)\to \R$. Similarly, the derivative $(\bar f\circ \bar \gamma)'(t)$ exists for $\leb^1$-a.e. $t\in I_\gamma$. Consequently, for $\leb^1$-a.e. $t\in \dom(\gamma)$ the limit
\[
(f \circ \gamma)'(t):= \lim_{\dom(\gamma)\ni s\to t}\frac{f(\gamma_s)-f(\gamma_t)}{s-t}
\]
exists. We will often abbreviate this for visual purposes as $(f \circ \gamma)'_t$. Again, the limit is independent of the choice of extension at non-isolated points of the domain $\dom(\gamma)$. 

The line integral of a Borel function $\rho:X\to [0,\infty]$ over $\gamma$ satisfies
\[
\int_\gamma\rho\ud s:=\int_{\dom(\gamma)}\rho(\gamma_t)|\gamma_t'|\ud t=\int_{\im(\gamma)}\rho\ud\Ha^1.
\]

\begin{remark}\label{rmk:reparam}
 Observe that
\begin{align*}
\ell^\ast(\gamma)\defl \int_{\dom(\gamma)}|\gamma_t'| \ud t+\gap(\gamma)=\ell(\bar\gamma).
\end{align*}
The curve $\bar\gamma$ can be parametrized by constant speed, see e.g. \cite[Theorem 4.2.1]{ambrosiotilli}.
Consequently also the curve fragment $\gamma$ has a bi-Lipschitz reparametrization $\tilde\gamma:\dom(\tilde\gamma)\to X$ with $[0,1]=I_{\tilde\gamma}$ such that $|\tilde\gamma_t'|=\ell^*(\gamma)$ for every $t\in \dom(\tilde\gamma)$. 
\end{remark}

The set of curve fragments on $X$ is denoted by $\Fr(X)$ and is given a topology as follows. Using the injective map $\gamma\mapsto {\rm graph}(\gamma)=\{(t,\gamma_t):\ t\in \dom(\gamma)\}$,
we identify $\Fr(X)$ with a subset of $\mathcal H(\R\times X)$ (the space of closed bounded subsets of $\R\times X$) and equip it with the Hausdorff metric. The following lemma follows by a straightforward application of the definitions and Lemma \ref{lem:sh-conv} and we omit the proof.
\begin{lemma}\label{lem:sh-ext-curve}
We have that $\bar\gamma|_{\dom(\gamma)}=\gamma$ and $\ell(\bar\gamma)=\Ha^1(\im(\gamma))+\gap(\gamma)$ for $\gamma\in\Fr(X)$. Moreover, $\gamma_j\to\gamma$ in $\Fr(X)$ if and only if $\bar\gamma_j\to\bar\gamma$ in $\Fr(\sh(X))$.
\end{lemma}

Next we define curve fragment analogues of upper gradients, called $\ast$-upper gradients.

\begin{definition}\label{def:star_ug}
Let $f\in \LIP(X)$. A Borel function $g:X\to[0,\infty]$ is a $\ast$-upper gradient of $f$ if
\begin{equation}\label{eq:star_ug_ineq-defn}
{\rm osc}_\gamma f\le \int_\gamma g\ud s+\LIP(f){\rm gap}(\gamma)
\end{equation}
for every $\gamma\in \Fr(X)$. 
\end{definition}
Here, we use the notation $\displaystyle {\rm osc}_\gamma f=\sup_{x,y\in \im(\gamma)}|f(x)-f(y)|$. Note that $\ast$-upper gradients are only defined for Lipschitz functions, and that $\LIP(f)$ is always a $\ast$-upper gradient of $f$.

\subsection{$\infty$-Modulus and minimal $\ast$-upper gradients}

We relax the notion of $\ast$-upper gradients by introducing negligible families of curve fragments in the sense of the \emph{$\infty$-modulus} $\Mod_\infty$. The resulting class of \emph{weak $\ast$-upper gradients} behaves better with respect to limits, and automatically contains minimal elements, making it a suitable analogue of minimal $p$-weak upper gradients in the context of curve fragments; compare \cite{HKST07,bjo11}. However, we will also see that compared to the theory for finite $p$, the exponent $p=\infty$ is special: a minimal $\ast$-upper gradient also exists.

The modulus $\Mod_\infty(\Gamma)$ of a family of curve fragments $\Gamma\subset \Fr(X)$ is defined by
\begin{align}
\Mod_\infty(\Gamma)=\inf\left\{\|\rho\|_{L^\infty(\mu)}:\ \int_\gamma\rho\ud s\ge 1\textrm{ for all }\gamma\in\Gamma \right\}.
\end{align}
The infimum is taken over all Borel functions $\rho:X\to [0,\infty]$. The $\infty$-modulus is an outer measure on $\Fr(X)$, and all the basic properties of $\Mod_p$ for finite $p$ remain valid for $\Mod_\infty$, see \cite{car-jar10}. Families $\Gamma_0\subset \Fr(X)$ of zero $\infty$-modulus have a particularly simple description, given below; compare \cite[Lemma 5.7]{car-jar10}.
\begin{lemma}\label{lem:Nexistence}
A family $\Gamma_0\subset \Fr(X)$ satisfies $\Mod_\infty(\Gamma_0)=0$ if and only if there exists a Borel set $N\subset X$ with $\mu(N)=0$ such that $\Gamma_0\subset \Gamma_N^+$, where 
\[
\Gamma_N^+:=\{\gamma\in\Fr(X):\ \Ha^1(N\cap \im(\gamma))>0\}.
\]
\end{lemma}
\begin{proof}
    The proof follows closely that of \cite[Lemma 5.7]{car-jar10}. If $\mu(N)=0$, let $\rho=\infty \cdot \chi_N$. Then $\int_\gamma \rho \ud s= \infty$ for all $\gamma\in \Gamma_N^+$. Further, $\|\rho\|_{L^\infty(\mu)}=0$, and thus $\Mod_\infty(\Gamma_N^+)=0$. This gives the converse direction of the lemma.

    If on the other hand $\Mod_\infty(\Gamma_0)=0$, then for every $\epsilon>0$, there exists a $\rho_\epsilon$ with $\|\rho_\epsilon\|_{L^\infty(\mu)}\leq \epsilon$ and $\int_\gamma \rho_\epsilon ds\geq 1$ for all $\gamma \in \Gamma_0$. Let $N = \bigcup_{n\in \N} \rho_{n^{-1}}^{-1}((n^{-1},\infty))$. From the definition of the $L^\infty(\mu)$-norm it is direct to see $\mu(N)=0$. Now, if $\gamma \in \Gamma_0$, we have some $n\in \N$ such that $\Ha^1(\im(\gamma))< n$. Then, $\Ha^1(\rho_{n^{-1}}^{-1}((n^{-1},\infty))\cap \im(\gamma))>0$ since otherwise $\int_\gamma \rho_{n^{-1}} ds \leq \Ha^1(\im(\gamma)) n^{-1} < 1$. Therefore, $\Ha^1(N\cap \im(\gamma))>0$ and $\gamma\in \Gamma_N^+$. Since this holds for all $\gamma\in \Gamma_0$, the claim follows.
\end{proof}

\begin{definition}
A Borel function $\rho\colon X\to [0,\infty]$ is a weak $\ast$-upper gradient of a function $f\in \LIP(X)$ if \eqref{eq:star_ug_ineq-defn} holds for every curve fragment outside a family $\Gamma_0\subset \Fr(X)$ with $\Mod_\infty(\Gamma_0)=0$.
\end{definition}

\begin{remark}\label{rmk:genuine-star-ug}
	In contrast with finite exponents, a weak $\ast$-upper gradient $\rho$ of $f\in\LIP(X)$ has a $\mu$-representative which is a genuine upper gradient. Indeed, if $\Gamma_0\subset \Fr(X)$ is the set where \eqref{eq:star_ug_ineq-defn} fails for $\rho$ and $N\subset X$ is the Borel $\mu$-null set with $\Gamma_0\subset \Gamma_N^+$, then
	\[
	\tilde\rho=\LIP(f)\chi_N+\rho\chi_{X\setminus N}
	\]
	is the required representative.  Indeed, if $\gamma$ is any curve fragment, then we can divide the domain as $\dom(\gamma)=\bigcup_{i=1}^\infty K_i \cup B$, where $K_i$ are compact sets such that  $\gamma(K_i)\subset N$ or $\gamma(K_i)\cap N=\emptyset$, and $B$ is a null set. Then $\tilde{\rho}$ is an $\ast$-upper gradient for each $\gamma|_{K_i}$. From this one can conclude that $\tilde{\rho}$ is an $\ast$-upper gradient for $\gamma$.
\end{remark}

\begin{proposition}\label{prop:min-star-ug-exists}
    Let $f\in \LIP(X)$. Then there exists a minimal $\ast$-upper gradient $|Df|_\ast$ of $f$ in the following sense: $|Df|_\ast$ is a $\ast$-upper gradient of $f$ and, if $\rho\in L^\infty(\mu)$ is a weak $\ast$-upper gradient of $f$, then $|Df|_\ast\le \rho$ $\mu$-a.e. Moreover, any two minimal $\ast$-upper gradients of $f$ agree $\mu$-a.e.
\end{proposition}
\begin{remark}\label{rmk:uppergradchar}
In the proof below we use the observation that a Borel function $\rho\in L^\infty(\mu)$ is a weak $\ast$-upper gradient of $f\in \LIP(X)$ if and only if 
\begin{align*}
|(f\circ\gamma)'(t)|\le \rho(\gamma_t)|\gamma'_t|\quad a.e.\ t\in \operatorname{dom}(\gamma)
\end{align*}
for $\Mod_\infty$-a.e. $\gamma\in \Fr(X)$. The proof also uses the concept of a measure theoretic infimum, which is also sometimes referred to as a lattice infimum.
\end{remark} 

\begin{definition}\label{def:measinfsup}
If $\mathcal{F}$ is any collection of measurable functions, then the measure theoretic infimum ${\rm essinf} \mathcal{F}$ is a measurable function $\rho$ such that $\rho \leq f$ for all $f\in \mathcal{F}$ and such that if $\rho'$ is any other such function, then $\rho' \leq \rho$ (a.e.).

Similarly, the measure theoretic supremum ${\rm esssup} \mathcal{F}$ is a measurable function $\rho$ such that $\rho \geq f$ for all $f\in \mathcal{F}$ and such that if $\rho'$ is any other such function, then $\rho' \geq \rho$ (a.e.).
\end{definition}

It can be shown that the measure theoretic infimum exists for $\sigma$-finite spaces, and that there exists always a sequence of $\rho_n\in \mathcal{F}$ such that $\inf_{n\in \N} \rho_n = {\rm essinf} \mathcal{F}$. The proof is quite standard, see e.g. \cite[Appendix A.2, Essential supremum]{teriseb}. By replacing functions by their negations, we can obtain the same results for measure theoretic suprema.

\begin{proof}[Proof of Proposition \ref{prop:min-star-ug-exists}]
Two minimal weak $\ast$-upper gradient must necessarily agree $\mu$-almost everywhere: if $\rho_1,\rho_2$ are both pointwise minimal, then $\rho_1\le \rho_2$ and $\rho_2\le \rho_1$ $\mu$-a.e. Thus it suffices to show its existence. 

We show the following lattice property of weak $\ast$-upper gradients: if $\rho_i\in L^\infty(\mu)$, $i=1,2,\ldots$ are weak $\ast$-upper gradients of $f$, then $\rho_\infty=\inf_i\rho_i$ is a weak $\ast$-upper gradient of $f$. Indeed, if $\Gamma_i\subset \Fr(X)$ is the $\Mod_\infty$-null set where \eqref{eq:star_ug_ineq-defn} fails for $\rho_i$ and $\Gamma_0=\Gamma_1\cup\Gamma_2\cup\cdots$, then for $\gamma\notin \Gamma_0$ there exist null sets $N_i\subset \dom(\gamma)$ such that $|(f\circ\gamma)'(t)|\le \rho_i(\gamma_t)|\gamma_t'|$ for $t\in \dom(\gamma)\setminus N_i$. It follows that 
\[|(f\circ\gamma)'(t)|\le \rho_\infty(\gamma_t)|\gamma_t'|
\]
for $t\in \dom(\gamma)\setminus(N_1\cup N_2\cup\cdots)$. Thus $\rho_\infty$ is a weak $\ast$-upper gradient of $f$ as claimed.

Now let $\rho$ be the measure theoretic infimum over all weak $\ast$-upper gradients of $f$, and let $\rho_i\in L^\infty(\mu)$ be a sequence from this family such that $\rho=\inf_i\rho_i$ $\mu$-a.e. By the preceding argument $\rho$ is a weak $\ast$-upper gradient of $f$, and by Remark \ref{rmk:genuine-star-ug} we may choose a representative of $\rho$ that is a genuine $\ast$-upper gradient. The definition of measure theoretic infimum implies that $\rho\le g$ $\mu$-a.e. for any weak $\ast$-upper gradient $g$ of $f$.  
\end{proof}

The lattice property established in the proof above yields the basic locality properties of the minimal $\ast$-upper gradient by standard arguments, see \cite[Chapter 6.3]{HKST07} or \cite[Chapter 2.2]{bjo11}. We record below the basic properties needed later on without proof, and refer the interested reader to the references above.

\begin{lemma}\label{lem:basic-prop-min-star-ug}
Let $f,g\in \LIP(X)$. Then
\begin{itemize}
	\item[(1)] $|D(f+g)|_\ast\le |Df|_\ast+|Dg|_\ast$ $\mu$-a.e.;
	\item[(2)] $|Df|_\ast=|Dg|_\ast$ $\mu$-a.e. on $\{f=g\}$; and 
	\item[(3)] $|Dh|_\ast=\chi_{\{f<g\}}|Df|_\ast+\chi_{\{f\ge g\}}|Dg|_\ast$ $\mu$-a.e. for $h=\min\{f,g\}$.
\end{itemize}
\end{lemma}

\begin{lemma}\label{lem:lim-star-ug}
Suppose $(f_j)\subset \LIP(X)$ converges pointwise to $f\in \LIP(X)$ and $\sup_j\LIP(f_j)<\infty$. If $\rho_j\in L^\infty(\mu)$ is a weak $\ast$-upper gradient of $f_j$ for each $j$ and $(\rho_j)$ converges pointwise $\mu$-a.e. to $\rho\in L^\infty$ with $\sup_j\|\rho_j\|_{L^\infty(\mu)}<\infty$, then $\rho$ is a weak $\ast$-upper gradient of $f$.
\end{lemma}
\begin{proof} 
Let $N=\bigcup_{i=0}^\infty N_i \subset X$ be the Borel set, where
\[
N_0:=\{x\in X: \lim_{j\to \infty} \rho_j(x)\neq \rho(x) \text{ or } \rho(x)> \sup_j\|\rho_j\|_{L^\infty(\mu)}\}
\]
and $N_j$ is a Borel set so that for $\gamma\not\in \Gamma_{N_j}^+$ the function $\rho_j$ is a $\ast$-upper gradient of $f_j$. The latter is guaranteed by Lemma \ref{lem:Nexistence}.

Then, whenever $\gamma\notin \Gamma_N^+$, we have that (a) $\rho_j$ is a $\ast$-upper gradient of $f_j$ along $\gamma$, and (b) $\rho_j\circ\gamma\to \rho\circ\gamma$ pointwise a.e. with $\sup_j\|\rho_j\circ\gamma\|_{L^\infty(\dom(\gamma))}<\infty$. By the dominated convergence theorem we have
\begin{align*}
|f(\gamma_t)-f(\gamma_s)|= & \lim_{j\to \infty}|f_j(\gamma_t)-f_j(\gamma_s)|\le \liminf_{j\to\infty}\bigg(\int_{\gamma}\rho_j\ud s+\LIP(f_j)\gap(\gamma)\bigg)\\
\leq & \int_\gamma\rho\ud s+C\gap(\gamma)
\end{align*}
for all $t,s\in\dom(\gamma)$, where $C=\sup_j\LIP(f_j)$. If $t$ is a Lebesgue density point of $\dom(\gamma)$ then $\leb^1([t-\varepsilon,t+\varepsilon]\setminus \dom(\gamma))/\varepsilon\to 0$ as $\varepsilon\to 0^+$. This and Lebesgue's differentiation theorem imply that $|(f\circ\gamma)_t'|\le \rho(\gamma_t)|\gamma_t'|$ a.e. $t\in \dom(\gamma)$ for $\gamma\notin\Gamma_N^+$. Thus $\rho$ is a weak $\ast$-upper gradient of $f$, by Remark \ref{rmk:uppergradchar}.
\end{proof}

The same conclusion remains true when the $\mu$-a.e. convergence $\rho_j\to \rho$ is weakened to $w^\ast$-convergence in $L^\infty(\mu)=L^1(\mu)^\ast$. The argument uses a slightly non-standard version of Mazur's lemma.
\begin{lemma}\label{lem:lim-star-ug-weak}
Suppose $(f_j)\subset \LIP(X)$ converges pointwise to $f\in \LIP(X)$ and $\sup_j\LIP(f_j)<\infty$. If $\rho_j\in L^\infty(\mu)$ is a weak $\ast$-upper gradient of $f_j$ for each $j$ and $(\rho_j)$ $w^\ast$-converges in $L^\infty(\mu)$ to $\rho$, then $\rho$ is a weak $\ast$-upper gradient of $f$.
\end{lemma}
\begin{proof}
First, by Banach-Steinhaus, we get $\sup_j\|\rho_j\|_{L^\infty(\mu)}<\infty$.
We will show that there exists non-negative numbers $\alpha_k(n), \dots \alpha_{N_n}(n)$, $N_n\in \N$ (with $n\leq N_n)$ and with $\sum_{j=n}^{N_n} \alpha_j(n) = 1$, so that
\[
\tilde{\rho}_n=\sum_{j=n}^{N_n} \alpha_j(n) \rho_j \stackrel{n\to\infty}{\longrightarrow} \rho
\]
converges pointwise almost everywhere and with $\sup_j\|\tilde{\rho}_j\|_{L^\infty(\mu)}<\infty$. If this is true, then it is direct to verify that $\tilde{\rho}_k$ is a weak $\ast$-upper gradient for $\tilde{f}_n=\sum_{j=n}^{N_n} \alpha_j(n) f_j$. Indeed, for $\Mod_\infty$-a.e. curve fragment $\gamma$, and any $s,t\in \dom(\gamma)$, we have
\begin{align*}
|\tilde{f}_n(\gamma_t)-\tilde{f}_n(\gamma_s)|&\le \sum_{j=n}^{N_n} \alpha_j(n) |f_j(\gamma_t)-f_j(\gamma_s)|
\\
& \leq \sum_{j=n}^{N_n} \alpha_j(n) \bigg(\int_{\gamma}\rho_j\ud s+\LIP(f_j)\gap(\gamma)\bigg)\\
\leq & \int_\gamma\tilde{\rho}_n \ud s+C\gap(\gamma),
\end{align*}
where $C=\sup_j \LIP(f_j)$. From this and Remark \ref{rmk:uppergradchar}, similar to the previous lemma, we get that $\tilde{\rho}_n$ is a weak $\ast$-upper gradient for $\tilde{f}_n$. Now, the claim follows from Lemma \ref{lem:lim-star-ug}. We are left to prove the existence of $\alpha_j(n)$ and $N_n$ such that the previous claims hold.

Fix $n\in \N$ and fix $x_0\in X$. Since $L^\infty(B(x_0,n),\mu) \subset L^2(B(x_0,n),\mu)\newline \subset L^1(B(x_0,n),\mu)$, we have that the sequence $\rho_j \chi_R \to \rho \chi_R$ weakly in $L^2(B(x_0,n),\mu)$. Thus, by Mazur's lemma, together with passing to a subsequence, we can find non-negative values $\beta^n_k(k), \dots \beta^n_{M^n_k}(k)$, $k\leq M^n_k\in \N$, with $\sum_{j=k}^{M^n_k} \beta_j(k) = 1$, so that
\[
\tilde{\rho}^n_k=\sum_{j=k}^{M^n_k} \beta_j(k) \rho_j \stackrel{k\to\infty}{\longrightarrow} \rho
\]
converges almost everywhere in $B(x_0,n)$. So, for some index $m_n\in \N$ we have $m_n\geq n$ and $\mu(\{x\in B(x_0,n): |\tilde{\rho}^n_{m_n}(x)-\rho(x)|\geq\frac{1}{n}\}) \leq 2^{-n}$. Set $N_n=M_{m_n}$, and $\alpha_j(n)=\beta_j(m)$ for $j\in [m,M_m]$, and otherwise $\alpha_j(n)=0$ for $j\in [n,m]$.

By the following Borel-Cantelli-type argument we have that almost everywhere
\[
\tilde{\rho}_n=\sum_{j=n}^{N_n} \alpha_j(n) \rho_j \stackrel{n\to\infty}{\longrightarrow} \rho.
\]
Indeed, let $A_n=\{x\in B(x_0,n): |\tilde{\rho}^n_{m_n}(x)-\rho(x)|\geq\frac{1}{n}\}=\{\{x\in B(x_0,n): |\tilde{\rho}_{n}(x)-\rho(x)|\geq\frac{1}{n}\}.$ Then for the limit superior set $\displaystyle \limsup_{n\to \infty} A_n = \bigcap_{m=1}^\infty \bigcup_{n\geq m} A_n$  we have 
\[
\mu\left(\bigcap_{m=1}^\infty \bigcup_{n\geq m} \limsup_{n\to \infty} A_n\right) = 0,
\]and for every $x\not\in \limsup_{n\to \infty} A_n$, we get pointwise convergence $\lim_{n\to \infty}\tilde{\rho}_n(x)=\rho(x)$. 
\end{proof}

\subsection{Plans and Alberti representations}\label{sec:1-plan+alb-rep}
Let $\P$ be a finite measure on $\Fr(X)$. The \emph{barycenter} $\P^\#$ of $\P$ is the Borel regular measure on $X$ given by
\begin{align*}
\P^\#(E):=\int\int_\gamma\chi_E\ud s\ud\P(\gamma),\quad E\subset X\textrm{ Borel.}
\end{align*}
\begin{definition}\label{def:1-plan}
A finite measure $\P$ on $\Fr(X)$ is called a 1-plan, if $\P^\#=\rho_{\P^\#}\mu$ for some $\rho\in L^1(\mu)$.
\end{definition}

\begin{definition}\label{def:alberti-rep}
Let $\mu$ be a Radon measure. An Alberti representation $\mathcal A=\{\mu_\gamma,\P\}$ of $\mu$ consists of a finite measure $\P$ on $\Fr(X)$ and, for $\P$-a.e. $\gamma$, a finite measure $\mu_\gamma\ll\Ha^1|_{\im(\gamma)}$ such that the following hold. For every Borel set $E\subset X$ we have that 
\begin{itemize}
    \item[(a)] $\displaystyle \gamma\mapsto \mu_\gamma(E)$ is Borel, and
    \item[(b)]$ \displaystyle \mu(E)=\int_{\Fr(X)}\mu_\gamma(E)\ud\P(\gamma)$.
\end{itemize}
\end{definition}
We record the close relationship between 1-plans and Alberti representations in the next remark.

\begin{remark}\label{rmk:alb-rep-vs-1-plans}
If $\P$ is a 1-plan on $X$ and $\rho$ is a Borel representative of $\frac{\ud \P^\#}{\ud\mu}$, then $\mu|_{\{\rho>0\}}$ has an Alberti representation $\mathcal A=\{\mu_\gamma,\P\}$, where $\mu_\gamma=\frac{\chi_{\{\rho>0\}}}{\rho} \Ha^1|_{\im(\gamma)}$.

Conversely, if $\mathcal A=\{\mu_\gamma,\P\} $ is an Alberti representation of $\mu$ with 
\[
\int\Ha^1(\spt(\mu_\gamma))\ud\P(\gamma)<\infty,
\]
there exists a compact $K\subset X$ such that, if $R:\Fr(X)\to \Fr(X)$ is the map $\gamma\mapsto \gamma|_{\gamma\inv(K)}$,
then $\P':=R_\ast\P$ is a 1-plan.
Such a $K$ can be obtained from the Lebesgue decomposition of $\mu$ with respect to $\P^\#$.
Note that $\P'$ is concentrated on subcurves of $\spt\P$.
\end{remark}

Depending on the application, there are advantages to both Alberti representations and plans. We highlight some crucial features to pay attention to. For example, Alberti representations are easier to restrict to subsets (we can just replace $\mu_\gamma$ with $\mu_\gamma|_A$ when $A\subset X$, while $\gamma|_{\gamma\inv(A)}$ may not be a curve fragment), whereas it is easier to take superpositions of plans (we can add $\P^1+\P^2$, but it is more difficult to combine $\mu_\gamma^1$ and $\mu_\gamma^2$ in order to preserve b) in Definition \ref{def:alberti-rep}).  Plans are a bit more flexible, since they do not impose (b) in Definition \ref{def:alberti-rep}, however, a single plan may not ``see'' all of the measure $\mu$ in the sense that $\mathbb{P}^\sharp(A)=0$ for some $\mu$-positive measure Borel set $A\subset X$. Moreover, for Alberti representations, a crucial technical issue is that $\mu_\gamma$ only ``sees'' a part of the curve fragment $\gamma$. These formal details are reflected in the details of the arguments employing the different notions, and have some subtle consequences for the formulations of different results.  A reader interested in these subtleties can see the arguments in \cite{sch16b,bate15}.

Denote by $C(v,\vartheta)\defl\{w\in \R^n:\ v\cdot w>\cos\vartheta |w| \}$ the (open) cone of width $\vartheta\in (0,\infty)$ in the direction $v\in S^{n-1}$. We say that a $k$-tuple $C_1,\ldots, C_k\subset \R^n$ of cones is \emph{independent}, if $\{w_1,\ldots, w_k\}$ are linearly independent whenever $w_j\in C_j$, $j=1,\ldots,k$ and \emph{$\varepsilon$-separated} ($\varepsilon>0$), if
\begin{align*}
	\bigg\|\sum_j^k\lambda_jw_j\bigg\|>\varepsilon \max_{1\le j\le k}\|\lambda_jw_j\|
\end{align*}
whenever $w_j\in C_j$, $\lambda_j\in\R$, $j=1,\ldots, k$. 

\begin{remark}\label{rmk:indep-vs-sep+conefields}
If $\bm v:X\to S^{n-1}$, $\bm\vartheta,\bm\varepsilon :X\to (0,\infty)$ are Borel functions, we may similarly define linearity and $\bm\varepsilon$-separatedness pointwise for \emph{cone fields}. Here, following the terminology of \cite{sch16b}, a cone field $\bm C=C(\bm v,\bm \vartheta)$ is the map $x\mapsto \bm C_x\defl C(\bm v_x,\bm\vartheta_x)$.

Observe that $\bm\varepsilon$-separated cone fields are always independent. Conversely, independent (closed) cones are $\bm\varepsilon$-separated for some Borel $\bm\varepsilon:X\to (0,\infty)$.
\end{remark}

\begin{definition}\label{def:indep-alb-rep}
Suppose $\varphi\in \LIP(X,\R^n)$. 

\begin{itemize}
	\item[(a)] An Alberti representation $\mathcal A=\{\mu_\gamma,\P\}$ of $\mu$ is in the $\varphi$-direction of a cone field $\bm C$, if $(\varphi\circ\gamma)'(t)\in \bm C_{\gamma_t}$ a.e. $t\in \dom(\gamma)$ for $\P$-a.e. $\gamma\in\Fr(X)$.
	\item[(b)] A collection $\mathcal A_1,\ldots, \mathcal A_k$ of Alberti representations is $\varphi$-independent (resp. $(\varphi,\bm\varepsilon)$-separated) if 
	they are in the $\varphi$-direction of independent (resp. $\bm\varepsilon$-separated) cone fields.
\end{itemize}
\end{definition}
 
\subsection{Weaver derivations}\label{sec:weaver-der}
A Weaver derivation is a bounded linear map 
\[
\bm b:\LIP_{\infty}(X)\to L^\infty(\mu)
\]
that (a) satisfies the Leibniz rule $\bm b(fg)=f\bm b(g)+g\bm b(f)$ $\mu$-a.e. for every $f,g\in \LIP_{\infty}(X)$,
and (b) is weak$^*$-continuous in the following sense: if a sequence $(f_j)\subset \LIP_{\infty}(X)$ 
with $\sup_j\|f_j\|_{\LIP_\infty}<\infty$ converges to $f\in\LIP_{\infty}(X)$ pointwise, then $\bm b(f_j)\stackrel{w^*}{\rightharpoonup} \bm b(f)$ in $L^\infty(\mu)$. The collection $\mathcal X(\mu)$ of Weaver derivations on $X$ is a Banach space when equipped with the norm $\|\bm b\|_{\mathcal X(\mu)}=\||\bm b|_{\mathcal X,\loc}\|_{L^\infty(\mu)}$, where
\begin{align}\label{eq:loc-norm}
	|\bm b|_{\mathcal X,\loc}\defl \esssup\{ |\bm b(f)|:\ \LIP(f)\le 1\}\in L^\infty(\mu)
\end{align}
is the pointwise norm defined as a measure theoretic supremum, see Definition \ref{def:measinfsup}. 

It is a standard fact that the estimate
\begin{align}\label{eq:deriv-ptwise-norm-bound}
	|\bm b(f)|\le |\bm b|_{\mathcal X,\loc}\Lip_af\quad \mu-a.e.
\end{align}
holds for every $f\in \LIP_\infty(X)$. Indeed, in order to introduce some of the common techniques involved in the study of Weaver derivations, we will sketch this argument. 

The Leibniz rule together with weak$^*$-continuity implies locality: if $A$ is a positive measure subset and $f|_A=g|_A$, then $\bm b(f)|_A=\bm b(g)|_A$ almost everywhere, see e.g. \cite[Lemma 13.4]{heinonen2007nonsmooth}. Now, for every ball $B(x,r)$, we can take $f_{x,r}$ to be the McShane extension of $f|_{B(x,r)}$, which is $\LIP(f|_{B(x,r)})$-Lipschitz. Then, locality implies that $\bm b(f_{x,r})|_{B(x,r)}=\bm b(f)|_{B(x,r)}$ almost everywhere. From \eqref{eq:loc-norm} we thus get 
\[
|\bm b(f)|_{B(x,r)}|\leq |\bm b|_{\mathcal X,\loc}\LIP(f|_{B(x,r)}).
\]
Then, taking a countable basis of the topology given by balls $\{B(x_i,r_i) : i\in \N\}$, we see that for a.e. $x\in X$ we have 
\[
|\bm b(f)|(x)\leq |\bm b|_{\mathcal X,\loc}(x)\inf\{\LIP(f|_{B(x_i,r_i)}): x\in B(x_i,r_i) \}. 
\]
Now, the final part is to observe that $\inf\{\LIP(f|_{B(x_i,r_i)}): x\in B(x_i,r_i) \} = \Lip_a(f)(x)$, whenever $\{B(x_i,r_i) : i\in \N\}$ forms a basis for the topology.

\begin{remark}\label{rmk:extend-deriv-to-LIP}
    The locality of Weaver derivation can also be used to define $\bm b(f)$ for any $f\in \LIP(X)$ by $\bm b(f):=\lim_{j\to\infty}\bm b(f_j)$, where $f_j=(1-\dist(\cdot,B(x_0,j))_+f$ for some arbitrary $x_0\in X$. The pointwise norm bound \eqref{eq:deriv-ptwise-norm-bound} remains valid for $f\in \LIP(X)$ and we will tacitly use this extension in the sequel.
\end{remark}

\section{Approximation}\label{sec:approximation}

Before turning to the representation of minimal $\ast$-upper gradients, in this section, we prove the main approximation result, Theorem \ref{thm:approximation}. This theorem will be used later in Sections \ref{sec:dense-dir} and \ref{sec:comparisons} in connection with the density of directions and the characterization of Lipschitz differentiability spaces (cf. Theorems \ref{thm:density_of_directions} and \ref{thm:LDS-char}).

\subsection{Curve fragments, line integrals, and lower semicontinuity}\label{sec:line-integral_lsc}

In the proof of Theorem \ref{thm:approximation}, but also of the representation Theorem \ref{thm:min_star_ug}, we will need semicontinuity properties of the line integral. However, in contrast to the situation with curves, the line integral of a continuous function need not be lower semicontinuous with respect to convergence of curve fragments, compare \cite[Proposition 4]{keith03}. Nevertheless here we obtain the semicontinuity of a closely related quantity, called the $\ast$-integral in \cite{bate18}, defined by
\begin{align*}
    \int^\ast_\gamma\rho\ud s:=\int_\gamma\rho\ud s+\gap(\rho,\gamma)
\end{align*}
for Borel functions $\rho:X\to [0,\infty]$ and $\gamma\in \Fr(X)$. Here the weighted gap-functional $\gap(\rho,\gamma)$ is defined by
\begin{align*}
\gap(\rho,\gamma)\defl \sum_i\frac{\rho(\gamma_{a_i})+\rho(\gamma_{b_i})}{2}d(\gamma_{a_i},\gamma_{b_i}),
\end{align*}
where $\bigcup_i(a_i,b_i)=I_\gamma\setminus\dom(\gamma)$. We also recall the edge length functional $l_{\sh}$ and the extension $\bar\gamma:I_\gamma\to \sh(X)$ of $\gamma$ (cf. Remark \ref{rmk:edge-length} and Section \ref{sec:lip-free}, respectively), and define $\gap_\delta$ for $\delta\ge 0$ by
\begin{align*}
\gap_\delta(\gamma)\defl \int_{\bar\gamma}\chi_{\{l_{\sh}>\delta\}}\ud s=\sum_{d(\gamma_{a_i}, \gamma_{b_i})>\delta} d(\gamma_{a_i}, \gamma_{b_i}),
\end{align*}
where $(a_i,b_i)$ are the gaps of the curve fragment $\gamma$; recall the terminology in Subsection \ref{sec:frag+star-ug}. 
Observe that $\gap(1,\gamma)=\gap(\gamma)=\gap_0(\gamma)$. 

\begin{proposition}\label{prop:liminf}
Suppose $(\gamma_j)\subset \Fr(X)$ converges to $\gamma$ in $\Fr(X)$, and $(\rho_j)$ is a pointwise non-decreasing sequence of lower semicontinuous functions with (pointwise) limit $\rho$. Then
\begin{align*}
\int_\gamma\rho\ud s+\gap(\rho,\gamma)\le \liminf_{j\to\infty}\left(\int_{\gamma_j}\rho_j\ud s+\gap(\rho_j,\gamma_j) \right)
\end{align*}
and 
\begin{align*}
    \gap_\delta(\gamma)\le \liminf_{j\to\infty}\gap_\delta(\gamma_j).
\end{align*}
\end{proposition}
We will prove this proposition after a  bit of preliminary work. Given any function $\rho:X\to \R\cup\{\infty\}$ we may define the extension $\bar \rho:\sh(X)\to \R\cup\{\infty\}$ by
\begin{align*}
\bar \rho(\iota_x+t(\iota_y-\iota_x))\defl t\rho(y)+(1-t)\rho(x).
\end{align*}
Using the properties of $\F(X)$ it is not difficult to see that the extension is well-defined. The extensions $\bar\gamma$ and $\bar\rho$ can be used to identify the $\ast$-integral with a line-integral over the extended curve in $\sh(X)$.

\begin{lemma}\label{lem:line-int}
    Let $\gamma:\dom(\gamma)\to X$ be a curve fragment and $\bar\gamma$ the associated extension. Let $\rho:X\to [0,\infty]$ be a Borel function and let $\bar\rho$ be it's extension. Then 
    \begin{align*}
    \int_{\gamma}\rho\ud s+\gap(\rho,\gamma)=\int_{\bar\gamma}\bar\rho\ud s.
    \end{align*}
\end{lemma}
\begin{proof}
    Given $\gamma$ with the decomposition using gaps, $I_\gamma\setminus \dom(\gamma)=\bigcup_i(a_i,b_i)$, consider the extended curve $\bar\gamma$ and note that, for each $i$, 
    \[
\int_{\bar\gamma|_{[a_i,b_i]}}\bar\rho \ud s=d(\gamma_{b_i},\gamma_{a_i})\int_{0}^{1}[t\rho(\gamma_{b_i})+(1-t)\rho(\gamma_{a_i})]\ud t=\frac{\rho(\gamma_{a_i})+\rho(\gamma_{b_i})}{2}d(\gamma_{b_i},\gamma_{a_i})
    \]
    so that 
\begin{align*}
\int_{\bar\gamma}\bar\rho \ud s=\int_{\dom(\gamma)}\rho\ud s+\int_{[a,b]\setminus\dom(\gamma)}\bar\rho \ud s=\int_{\dom(\gamma)}\rho\ud s + \gap(\rho,\gamma),
\end{align*}
as claimed.
\end{proof}

\begin{lemma}\label{lem:sh-ext-lsc}
    If $\rho:X\to [0,\infty]$ is lower semicontinuous, then $\bar\rho:\sh(X)\to [0,\infty]$ is lower semicontinuous.
\end{lemma}

\begin{proof}[Proof of Lemma \ref{lem:sh-ext-lsc}]
    Let $v=\iota_{x}+t(\iota_{y}-\iota_{x})$ with $t\le 1/2$ and $v_j\to v$ in $sh(X)$. Let $(x_j,y_j)$ and $(t_j)$ be as in Lemma \ref{lem:sh-conv}. If $x\ne y$ and $t\in (0,1)$, then
    \begin{align*}
\bar\rho(v)=&t\rho(y)+(1-t)\rho(x)\le \liminf_{j\to\infty} t_j\rho(y_j)+\liminf_{j\to\infty}(1-t_j)\rho(x_j)\\
\le & \liminf_{j\to\infty}[t_j\rho(y_j)+(1-t_j)\rho(x_j)]= \liminf_{j\to\infty}\bar\rho(v_j).
    \end{align*}

In the second case we have that $v=\iota_x$ and $t=0$, Thus $t_j\to 0$ and $x_j\to x$. In this case
\[
\bar\rho(v)=\rho(x)\le \liminf_{j\to\infty}(1-t_j)\rho(x_j)\le \liminf_{j\to\infty}[t_j\rho(x_j)+(1-t_j)\rho(y_j)]= \liminf_{j\to\infty}\bar\rho(v_j).
\]
This completes the proof.
\end{proof}

We now provide the proof of the main result of this subsection.

\begin{proof}[Proof of Proposition \ref{prop:liminf}]
The first estimate follows from \cite[Lemma 4.2]{davideb20} and Lemmas \ref{lem:sh-ext-curve} and \ref{lem:line-int}. The second similarly follows from \cite[Lemma 4.2]{davideb20} since $\chi_{\{l_{\sh}>\delta\}}$ is lower semicontinuous by Remark \ref{rmk:edge-length}.
\end{proof}

\subsection{Construction of an approximating sequence}

Suppose $f\in \LIP_{bs}(X)$. By Remark \ref{rmk:genuine-star-ug}, $|Df|_\ast$ has a Borel representative $\rho$ that is a genuine $\ast$-upper gradient of $f$. Since $f$ has bounded support, $\rho$ vanishes outside some bounded set. We fix $\rho$ for this construction.
 
If $j\in\N$, we set
\[
E_{j,l}\defl\rho\inv\bigg(\Big(\frac{l-1}{2^{j}},\frac{l}{2^{j}}\Big]\bigg),
\]
for $l=1,\ldots, N_j$, where $N_j$ is the unique integer in $(2^j\|\rho\|_{L^\infty(\mu)},2^j\|\rho\|_{L^\infty(\mu)}+1]$. 
Given $j,l$, let $K_{j,l}$ be a compact subsets of $E_{j,l}$, such that
\[
\mu(E_{j,l}\setminus K_{j,l})\le 2^{-j^2}.
\]
Set $\displaystyle K_j:=\bigcup_{l=1}^{N_j}K_{j,l}$. For any compact set $A\subset X$ and $r>0$, let $B(A,r)=\{x\in X: d(x,A)< r\}$.  Since $K_{j,l}$ and $K_{j,\tilde l}$ are at a positive distance for distinct $l,\tilde l$, we may find $k_j$ be such that
\[2^{-k_j}<\frac{\min_{l\ne\tilde l}\dist(K_{j,l},K_{j,\tilde l})}{2}.\]
Then, for each $k\geq k_j$ define the functions $\psi_{j,l}^{k}:X\to [0,1]$, 
\[
\psi_{j,l}^{k}:=\frac{\dist(x,X\setminus B(K_{j},2^{-k}))}{\dist(x,X \setminus B(K_{j},2^{-k}))+\dist(x,K_{j,l})},\]
which vanish outside $B(K_{j},2^{-k})$ and equal $1$ on $K_{j,l}$. Finally set
\begin{align*}
\rho_j^{k}:=\LIP(f)-\sum_{l=1}^{N_j}\Big(\LIP(f)-\frac{l}{2^{j}}\Big)_+\psi_{j,l}^{k}.
\end{align*}

\begin{remark}\label{rmk:ptwise-incr-limit}
We have for every $k\geq k_j$ that $\rho_j^{k}=\min\{l/2^j,\LIP(f)\}$ on $K_{j,l}$, $\rho_j^{k}=\LIP(f)$ outside $B(K_j,2^{-k})$ and $\rho_j^{k}\le \LIP(f)$ everywhere. Furthermore, $(\rho_j^{k})_k$ is a non-decreasing sequence of continuous functions converging pointwise to the function
\[
\rho_j\defl\LIP(f)\chi_{X\setminus K_j}+\sum_{l=1}^{N_j}\frac{l}{2^j}\chi_{K_{j,l}},
\]
which satisfies $\rho\le \rho_{j}$ and is thus a $\ast$-upper gradient of $f$. Further, $\rho_j \geq \rho_j^k \geq 2^{-j}$.
\end{remark}

\begin{lemma}\label{lem:ug+e}
    Given $j\in\N$, and $\varepsilon>0$ there exists $k=k(j,\varepsilon)\in \N$ such that 
    \begin{align*}
    |f(y)-f(x)|\le &\int_{\gamma}\rho_j^k\ud s +\sum_{d(\gamma_{a_i},\gamma_{b_i})\le 2^{-k}}\frac{\rho^k_j(\gamma_{a_i})+\rho^k_j(\gamma_{b_i})}{2}d(\gamma_{a_i},\gamma_{b_i})\\ & +\LIP(f)\gap_{2^{-k}}(\gamma)+\varepsilon
    \end{align*}
    for all $\gamma\in\Fr(X)$ with start and endpoints $x$ and $y$, respectively, and
    \[
    \displaystyle \bigcup_i(a_i,b_i)=\conv{\dom(\gamma)}\setminus \dom(\gamma).
    \]
\end{lemma} 

The proof scheme is based on a contradiction argument that has appeared much earlier. Indeed,  Cheeger attributes it to Ziemer in \cite[Remark 5.22]{che99}. For the proof we define the function $\widehat\rho_j^k:\sh(X)\to [0,\infty)$ by
\[
\widehat\rho_j^k:=\max\{\overline{\rho_j^k},\LIP(f)\chi_{\{l_{\sh(X)}>2^{-k}\}}\}
\]
for each $j,k\in\N$, where $\overline{\rho_j^k}$ denotes the extension of $\rho_j^k$ into $\sh(X)$, cf. Section \ref{sec:line-integral_lsc}, and note that
\begin{align*}
\int_{\bar\gamma}\widehat\rho_j^{k}\ud s
= \int_{\gamma}\rho_j^k\ud s &+\sum_{d(\gamma_{a_i},\gamma_{b_i})\le 2^{-k}}\frac{\rho^k_j(\gamma_{a_i})+\rho^k_j(\gamma_{b_i})}{2}d(\gamma_{a_i},\gamma_{b_i})+\LIP(f)\gap_{2^{-k}}(\gamma).
\end{align*}

\begin{remark}
The line-integral above has the following geometric interpretation: from the perspective of the curve fragment $\gamma$, the line integral of $\widehat \rho_j^k$ over $\bar\gamma$ penalizes gaps in $\gamma$ with length above a certain threshold, while short gaps only contribute an amount determined by $\rho_j^k$.
\end{remark}

\begin{proof}[Proof of Lemma \ref{lem:ug+e}]
Observe that, as the maximum of two lower semicontinuous functions, $\widehat\rho_j^k$ is lower semicontinuous, cf. Lemma \ref{lem:sh-ext-lsc} and Remark \ref{rmk:edge-length}. Moreover the sequence $(\widehat\rho_j^k)_k$ increases pointwise to $\widehat\rho_j:=\chi_{X}\rho_j+\LIP(f)\chi_{\sh(X)\setminus X}$ as $k\to \infty$, see Remark \ref{rmk:ptwise-incr-limit}. Note that
\begin{align*}
\int_{\bar\gamma}\widehat \rho_j\ud s=\int_\gamma\rho_j\ud s+\LIP(f)\gap(\gamma),\quad \gamma\in\Fr(X).
\end{align*}
Fix $\varepsilon>0$ and suppose by contradiction that, for every $k\ge k_j$, there exists $\gamma_k\in\Fr(X)$ with start and endpoints $x_k,y_k\in X$ respectively, such that
\begin{equation}\label{eq:badcurvefrag}
   |f( y_k)-f( x_k)|\ge \int_{\overline{\gamma_k}}\widehat\rho_{j}^{k}\ud s+\varepsilon .
\end{equation}

Choose $\gamma_k$ to be a minimal curve-fragment such that \eqref{eq:badcurvefrag} holds, where curve fragments are ordered according to the inclusions of their graphs. By Zorn's lemma, since the intersection of any chain of curve fragments satisfying the previous inequality will also satisfy that same inequality (by lower-semicontinuity of the integrand in \eqref{eq:badcurvefrag}), such a minimal curve fragment exists. Denote $C_k=\dom(\gamma_k)$ and $\conv C_k=[a_k,b_k]$.

We claim that $\im(\gamma_k)\subset \overline B(K_j,2^{1-k})$. The idea here is that $\widehat\rho_j^k=\LIP(f)$ on the parts of the fragment lying outside $\overline B(K_j,2^{1-k})$ and thus they can be removed while preserving \eqref{eq:badcurvefrag}. (Indeed, these modifications can be performed starting from any curve-fragment $\gamma_k$ in order to change it, in a finite number of steps to satisy this condition and \eqref{eq:badcurvefrag}. This would give another proof that avoids the use of minimality or Zorns's lemma.)

Suppose that $F_k:=\gamma_k\inv(X\setminus \overline B(K_j,2^{1-k}))\ne \varnothing$. If $t_0\in F_k$ is isolated in $C_k$, let $t_0^\pm\in (C_k\setminus\{t_0\})\cup\{a_k,b_k\}$, $t_0^-\le t_0\le t_0^+$, be its nearest points. Consider two cases $\sigma \in \{-,+\}.$ Either $\gamma_k(t_0^\sigma)\in \overline B(K_j,2^{-k})$, in which case $d(\gamma_k(t_0),\gamma_k(t_0^\sigma))>2^{-k}$, or $\gamma_k(t_0^\sigma)\notin \overline B(K_j,2^{-k})$, in which case $\rho_j^k(t_0^\sigma)=\LIP(f)$. In both cases we have $\widehat \rho_j^k(\overline\gamma_k(t)))=\LIP(f)$ for all $t\in (t_0^\pm,t_0]\cup [t_0,t_0^\pm)$. Consequently $\widehat\rho_j^k(\overline\gamma_k(t))=\LIP(f)$ for $t\in \{t_0\}\cup(t_0^-,t_0^+)$ and 
\begin{align*}
\int_{\overline\gamma_k|_{[t_0^-,t_0^+]}}\widehat \rho_j^k\ud s=\LIP(f)[d(\gamma_{t_0^-},\gamma_{t_0})+d(\gamma_{t_0},\gamma_{t_0^+})]\ge \LIP(f)d(\gamma_{t_0^-},\gamma_{t_0^+})
\end{align*}

Considering the sub-fragment $\gamma_k':=\gamma_k|_{C_k\setminus \{t_0\}}$ of $\gamma_k$, we have from the above estimate that
\begin{align*}
\int_{\overline{\gamma_k}}\widehat\rho_j^k\ud s&=\int_{\overline{\gamma_k|_{[a_k,t_0^-]}}}\widehat\rho_j^k\ud s+\int_{\overline\gamma_k|_{[t_0^-,t_0^+]}}\widehat \rho_j^k\ud s+\int_{\overline{\gamma_k|_{[t_0^+,b_k]}}}\widehat\rho_j^k\ud s\\
&\ge \int_{\overline{\gamma_k|_{[a_k,t_0^-]}}}\widehat\rho_j^k\ud s+\LIP(f)d(\gamma_{t_0^-},\gamma_{t_0^+})+\int_{\overline{\gamma_k|_{[t_0^+,b_k]}}}\widehat\rho_j^k\ud s\\
&=\int_{\overline{\gamma_k'}}\widehat\rho_j^k\ud s
\end{align*}

\noindent It follows that the proper sub-fragment $\gamma_k'$ of $\gamma_k$ satisfies \eqref{eq:badcurvefrag}, contradicting the minimality of $\gamma_k$.

Therefore, we have that every point in $F_k$ is an accumulation point of $C_k$ and $F_k$ is relatively open. Thus, there exist $a_k',b_k'\in F_k$ such that $\varnothing\ne [a_k',b_k']\cap C_k\subset F_k$. In particular $\widehat\rho_j^k(\overline\gamma_k(t))=\LIP(f)$ for all $t\in [a_k',b_k']$. As above we obtain that $|f(\gamma_k(b_k'))-f(\gamma_k(a_k'))|\le \int_{\overline\gamma_k|_{[a_k',b_k']}}\widehat \rho_j^k\ud s$, and that $\gamma_k':=\gamma_k|_{C_k\setminus(a_k',b_k')}$ satisfies \eqref{eq:badcurvefrag}, contradicting the minimality of $\gamma_k$.
 
Thus $\im(\gamma_k)\subset \overline B(K_j,2^{1-k})$ for all $k$ and, recalling that $\widehat \rho_j^k\ge 2^{-j}$, we moreover have 
\[
2^{-j}\ell^\ast(\gamma_k)\le \varepsilon+\int_{\bar\gamma_k}\hat\rho_{j}^{k}\ud s\le |f(y_k)-f(x_k)|\le \LIP(f)\diam(\spt(f))
\]
so that $\sup_{k\geq k_j}\ell^\ast(\gamma_k)\defr L<\infty$ (see Remark \ref{rmk:reparam} for the definition of $\ell^*$). By Remark \ref{rmk:reparam}, we can parametrize all curves $\bar \gamma_k$ to be $L-$Lipschitz with the parameter interval $[0,1]$. 
Then the sequence $\bar\gamma_k:[0,1]\to \sh(X)$ of uniformly Lipschitz curves satisfies $\bar\gamma_k([0,1])\subset \sh(\overline{B(K_j,2^{-k})})$ and thus, it is not hard to see that for each $t\in [0,1]$, the sequence $(\bar\gamma_k(t))\subset \sh(X)$ is precompact. By the Arzel\'a--Ascoli Theorem \cite[Theorem 10.28]{hei01} a subsequence of $(\bar\gamma_k)$ converges to a Lipschitz curve $\bar\gamma:[0,1]\to \sh(X)$, and thus the corresponding subsequence of $(\gamma_k)\subset \Fr(X)$ converges to the curve fragment $\gamma\defl \bar\gamma|_{\bar\gamma\inv(X)}$, cf. Lemma \ref{lem:sh-ext-curve}. By Remark \ref{rmk:ptwise-incr-limit} and Proposition \ref{prop:liminf} we have that
\begin{align*}
\varepsilon+\int_\gamma\rho_j\ud s+\LIP(f)\gap(\gamma)&=\varepsilon+\int_{\bar\gamma}\hat\rho_j\ud s\le \varepsilon+\liminf_{k\to\infty}\int_{\bar\gamma_k}\hat\rho_{j}^{k}\ud s\\
&\le \liminf_{k\to\infty}|f(y_k)-f(x_k)|=|f(y)-f(x)|,
\end{align*}
which contradicts the fact that $\rho_j$ is a $\ast$-upper gradient of $f$, cf. Remark \ref{rmk:ptwise-incr-limit}. This completes the proof. 
\end{proof}

By Lemma \ref{lem:ug+e} there exists $k=k(j)>j^2$ so that  $\hat g_j:=\widehat\rho_j^{k(j)}:\sh(X)\to (0,\infty)$ satisfies
\begin{align}\label{eq:ug+e}
|f(x)-f(y)|\le \frac{\LIP(f)}{2^{j+1}}+\int_{\overline\gamma}\hat g_j\ud s,\quad \gamma\in\Fr(X).
\end{align}
We define
\begin{align*}
f_j(x)=\inf\left\{ f(y)+\int_{\overline\gamma}\hat g_j\ud s\Big|\Big.\ \gamma\in\Fr(X),\ \gamma:y\curvearrowright x\right\},
\end{align*}
where we write $\gamma:y\curvearrowright x$ to denote that the curve fragment starts at $y$ and ends at $x$, that is $\gamma(\min(\dom(\gamma)))=y,\gamma(\max(\dom(\gamma)))=x$.

In the next lemma and proposition we collect some key properties of the functions $\hat g_j$ and $f_j$ defined above.

\begin{lemma}\label{lem:rho-jm}
For each $j\in\N$ we have that 
\begin{itemize}
	\item[(a)] $\hat g_j|_X$ is continuous and $\min\{2^{-j},\LIP(f)\}\le \hat g_j|_X\le\LIP(f)$;
	\item[(b)] $|\hat g_j-\rho|\le 2^{-j+1}$ on $K_j$; and 
	\item[(c)] $\mu(\{\rho>0\}\setminus K_j)\le 2^{-j^2}(2^j\|\rho\|_{L^\infty(\mu)}+1)$.
\end{itemize}
\end{lemma}
\begin{proof}
The claims are immediate from the construction, see also Remark \ref{rmk:ptwise-incr-limit}.
\end{proof}

\begin{proposition}\label{prop:fjm}
Suppose $f\in \LIP_{bs}(X)$ is non-negative. For each $j\in\N$,
\begin{itemize}
\item[(a)] the function $f_j$ satisfies $0\le f_j(x)\le f(x)$ for every $x\in X$;
\item[(b)] there exists $\delta=\delta(j)>0$ so that, for each $x,y\in X$ with $d(x,y)\le \delta$, we have
\begin{align*}
    |f_j(x)-f_j(y)|\le \frac{\hat g_j(x)+\hat g_j(y)}{2}d(x,y),
\end{align*}
while $|f_j(x)-f_j(y)|\le \LIP(f)d(x,y)$ if $d(x,y)>\delta$; and
\item[(c)] we have that $|f(x)-f_j(x)|\le 2^{-j}\LIP(f)$ for all $x\in X$.
\end{itemize}
\end{proposition}
\begin{proof}
The statement in (a) is immediate from the construction. To see this, it is useful to observe that $f,\hat{g}_j\geq 0$ and that a curve fragment $\gamma$ can be chosen to consist of only one point.

To see (b), let $k=k(j)\in\N$ be as in \eqref{eq:ug+e} and the definition of $f_j$. Suppose $x,y\in X$ and assume without loss of generality that $f_j(x)\le f_j(y)$. Given $\varepsilon>0$, let $\gamma:y_\gamma\curvearrowright x$ be a curve fragment such that $f_j(x)+\varepsilon> f(y_\gamma)+\int_{\overline{\gamma}} \hat g_j\ud s$. Set $t_y=\max\dom(\gamma)+d(x,y)$, $\dom(\tilde\gamma):=\dom(\gamma)\cup\{t_y\}$, and define $\tilde\gamma:\dom(\tilde\gamma)\to X$ by $\tilde\gamma|_{\dom(\gamma)}=\gamma$ and $\tilde\gamma(t_y)=y$. If $d(x,y)\le 2^{-k}\defr \delta$, then
\[
\int_{\overline{\tilde\gamma}}\hat g_j\ud s=\int_{\overline\gamma}\hat g_j\ud s+\frac{\hat g_j(x)+\hat g_j(y)}{2}d(x,y)
\]
Consequently
\begin{align*}
f_j(y)\le f(y_\gamma)+\int_{\overline{\tilde\gamma}}\hat g_j\ud s&=f(y_\gamma)+\int_{\overline\gamma}\hat g_j\ud s+\frac{\hat g_j(x)+\hat g_j(y)}{2}d(x,y)\\
&<f_j(x)+\varepsilon+\frac{\hat g_j(x)+\hat g_j(y)}{2}d(x,y)
\end{align*}
and thus $f_j(y)-f_j(x)\le \frac{\hat g_j(x)+\hat g_j(y)}{2}d(x,y)+\varepsilon$. If $d(x,y)> \delta$ we have 
\[
\int_{\overline{\tilde\gamma}}\hat g_j\ud s\le \int_{\overline\gamma}\hat g_j\ud s+\LIP(f)d(x,y)
\]
and the argument above yields $f_j(y)-f_j(x)\le \LIP(f)d(x,y)+\varepsilon$. Since $\varepsilon$ is arbitrary, (b) follows.

It remains to show (c). Let $x\in X$ and, for each $j$, let $\gamma:y_\gamma\curvearrowright x$ be such that
\[f(y_\gamma)+\int_{\overline\gamma}\hat g_j\ud s<f_j(x)+\frac{\LIP(f)}{2^{j+1}}.\]
By \eqref{eq:ug+e} we have that
\[
f(x)\le f(y_\gamma)+\int_{\overline\gamma}\hat g_j\ud s+\frac{\LIP(f)}{2^{j+1}}\le f_j(x)+\frac{\LIP(f)}{2^{j}}.
\]
Since $f_j\le f$ it follows that $|f(x)-f_j(x)|\le 2^{-j}\LIP(f)$, completing the proof.
\end{proof}

We now prove Theorem \ref{thm:approximation}.
\begin{proof}[Proof of Theorem \ref{thm:approximation}]
Let $f^+=\max\{f,0\}$ and $f^-=-\min\{f,0\}$, so that $f^\pm\in\LIP_{bs}(X)$ is non-negative and $f=f^+-f^-$. Let $(f_j^\pm)$ be the sequence constructed above for $f^\pm$, with minimal $\ast$-upper gradient $\tilde\rho_j^\pm$ and $\rho^\pm$, respectively, satisfying \eqref{eq:ug+e}.  Set $f_j\defl f_j^+-f_j^-$. Then $\LIP(f_j)\le \LIP(f)$,  $|f_j|\le |f|$, and $f_j\to f$ pointwise by Proposition \ref{prop:fjm}.

It remains to show (2). Set $\tilde\rho_j\defl \tilde\rho_j^++\tilde\rho_j^-$ and $\rho\defl \rho^++\rho^-$. We observe that
\begin{align}\label{squash_Lip_a}
	\Lip_af_j\le \tilde\rho_j \ \textrm{ on } X,\quad \tilde\rho_j\to \rho \ \mu-a.e. \textrm{ on }X
\end{align}
by Proposition \ref{prop:fjm}(b) (with the continuity of $\tilde\rho_j$) and Lemma \ref{lem:rho-jm}(b),(c), respectively. By Lemma \ref{lem:basic-prop-min-star-ug} we have that $|Df|_\ast=\rho$ $\mu$-a.e..

Let $A\subset X$ be a Borel set with $\mu(A)<\infty$. We claim that $\Lip_af_j\rightharpoonup \rho$ weakly in $L^2(A)$. Indeed, since the sequence $(\Lip_af_j)$ is bounded in $L^2(A)$, given any subsequence, a further subsequence converges weakly to some $g\in L^2(A)$. By Mazur's Lemma a sequence of convex combination's of the tail of this subsequence converges to $g$ in $L^2(A)$ and pointwise. Since $\Lip_af_j\le \rho_j\to \rho$ pointwise, we have that $g\le \rho$ on $A$. On the other hand by Lemma \ref{lem:lim-star-ug} $g$ is a weak $\ast$-upper gradient of $f$ (which equals the limit of the corresponding convex combination of $(f_j)$), implying $\rho\le g$. We have proven that any subsequence of $(\Lip_af_j)$ has a further subsequence converging weakly in $L^2(A)$ to $\rho$, establishing the claim. 

For every $\varepsilon>0$ we have the estimate
\begin{align*}
\mu(A\cap \{ |\rho-\Lip_af_j|\ge \varepsilon\})\le &\frac 1\varepsilon\int_A|\rho-\Lip_af_j|\ud\mu\\
\le & \frac 1\varepsilon\int_A|\rho-\rho_j|\ud\mu+\frac 1\varepsilon\int_A(\rho_j-\Lip_af_j)\ud\mu.
\end{align*}
The $\mu$-a.e. pointwise convergence $\rho_j\to \rho$, the weak $L^2(A)$-convergence $\rho_j-\Lip_af_j\rightharpoonup 0$, and the dominated convergence theorem imply that
\[
\limsup_{j\to\infty}\mu(A\cap \{ |\rho-\Lip_af_j|\ge \varepsilon\})=0,
\]
i.e. $\Lip_af_j\to \rho$ in measure (in $A$). Together with the weak $L^2(A)$-convergence this implies pointwise convergence $\mu$-a.e. in $A$ after passing to a subsequence. Since $A\subset X$ is arbitrary this shows that, after passing to a subsequence, $\Lip_af_j\to\rho$ $\mu$-a.e. on $X$. 

Since $|Dh|_\ast\le \Lip h\le \Lip_a h$ for $h\in \LIP(X)$, \eqref{squash_Lip_a} holds for $|Df_j|_\ast$ and $\Lip f_j$ in place of $\Lip_a f_j$. Therefore, the same argument yields the $\mu$-a.e. convergence $|Df_j|_\ast\to \rho$ and $\Lip f_j\to \rho$, finishing the proof. 

\end{proof}

\section{Representation of minimal $\ast$-upper gradients}\label{sec:repgradient}

In this section we establish the following representation of the minimal weak $\ast$-upper gradient as a maximal directional derivative along curve fragments, analogous to \cite[Corollary 1.10]{che-kle-sch16}, \cite[Theorem 1.1]{teriseb} and \cite[Theorem 4.2]{cheegerkleiner}, as a tool towards constructing the pointwise norm and differential of functions in Theorem \ref{thm:diff-and-ptwise-norm}.

\begin{theorem}\label{thm:min_star_ug}
Let $f\in \LIP(X)$. There exists a 1-plan $\P$ with $\mu|_{\{|Df|>0\}}\ll\P^\#$ and disintegration $\{\bar\P_x\}$, such that 
\begin{align*}
|Df|_{\ast}(x)=\chi_{\{|Df|>0\}}(x)\left\|\frac{(f\circ\gamma)'(t)}{|\gamma_t'|}\right\|_{L^\infty(\bar\P_x)}
\end{align*}
$\mu$-a.e. Here $|\gamma_t'|$ is the metric derivative of the curve fragment $\gamma$ at $t\in \dom(\gamma)$.
\end{theorem}

The representation is based on the disintegration of a measure with respect to a mapping \cite[Theorem 5.3.1]{AGS08}: if $\Phi:X\to Y$ is a Borel measurable map between two separable metric spaces, and $\mu$ is a finite measure on $X$, then the disintegration of $\mu$ with respect to $\Phi$ is a collection $\{\mu_y\}\subset \M(X)$ for $\Phi_\ast\mu$-a.e. $y\in Y$ such that 
\begin{enumerate}
    \item for $\Phi_\ast\mu$-a.e. $y\in Y$, $\mu_y$ is concentrated on $f^{-1}(y)$; 
    \item for every Borel set $B$ in $X$, the map $y\mapsto \mu_y(B)$ is Borel; and
    \item for every measurable set $A\subset X$ we have
    \[
    \int_X \chi_A(x) \ud\mu(x) = \int_Y \int_X \chi_A(x) d\mu_y(x) \ud\Phi_\ast\mu(y).
    \]
\end{enumerate}
Moreover, if $\{\tilde\mu_y\}$ is another collection satisfying (1) and (2) above, then $\mu_y=\tilde \mu_y$ for $\Phi_\ast\mu$-a.e. $y\in Y$. See also \cite[III-70]{dellacherie2011probabilities}.

Set $$\overline\Fr(X)\defl \{ (\gamma,t)\in \Fr(X)\times \R:\ \gamma\in\Fr(X),\ t\in \dom(\gamma) \},$$ and equip it with the subspace topology of $\Fr(X)\times \R$. Define the evaluation map $e:\overline\Fr(X)\to X$, $(\gamma,t)\mapsto \gamma_t$. Given a finite measure $\P \in\mathcal M(\Fr(X))$, the measure $\ud\bar\P\defl \ud s\ud\P$ on $\overline\Fr(X)$ is given by 
\begin{align*}
	\bar\P(B)=\int_{\Fr(X)}\int_\gamma\chi_B(\gamma,\cdot)\ud s\ud\P(\gamma),\quad B\subset \overline\Fr(X)\ \textrm{ Borel}.
\end{align*}
Observe that $e_\ast\bar\P=\P^\#$ so that, if $\P$ is a 1-plan, then $\P^\#=\rho_{\P^\#}\mu$, where $\rho_{\P^\#}\in L^1(\mu)$ is the Radon--Nikodym derivative of $\P^\#$ with respect to $\mu$. Consequently, if $\P$ is a 1-plan, then the disintegration $\{\bar\P_x\}$ of $\bar\P$ with respect to $e$ is defined for $\mu$-a.e. $x\in \{\rho_{\P^\#}>0\}$.

\subsection{Plan-modulus duality} The proof of Theorem \ref{thm:min_star_ug} is based on duality, and is modeled on the argument in \cite{teriseb}. However, compared to the setting in \cite{teriseb}, we face the central difficulty that the duality theory of modulus of curve families from \cite{amb13} is only applicable for a finite exponent $p\in (1,\infty)$. (See, however, \cite{exnerova2019plans, shanmunsemmes, davideb20} for duality results that are applicable for $p=1$.) This duality theory does not extend directly to the case of $p=\infty$.  As a substitute, for the exponent $p=\infty$, we prove a partial duality result (Proposition \ref{prop:rho-bad}) which is based on the use of Rainwater's lemma. Historically, this is quite natural, since this lemma was used in \cite{almar16} for the decomposability bundle of a measure on Euclidean space, and in \cite{bate15} and \cite{sch16b} to obtain results on the geometry of Lipschitz differentiability spaces. Thus, Rainwater's Lemma really should be thought of as the $p=\infty$ substitute for the duality result in \cite{amb13}. Indeed, a close examination of the proofs in \cite{amb13} and the following shows that these are related.

\begin{lemma}\label{lem:rainwater}
    Let $\mu$ be a finite measure on a metric space $X$ and $\Gamma$ a compact family of curve fragments. Then one of the following two conditions holds.
    \begin{itemize}
        \item[(1)] There exists a $\mu$-null Borel set $N\subset X$ such that $\Ha^1|_{\im(\gamma)}$ is concentrated on $N$ for every $\gamma\in\Gamma$; or
        \item[(2)] there exists a compact set $K\subset X$ with $\mu(K)>0$ and a 1-plan $\mathbb P\ne 0$ concentrated on the curve fragment family $\Gamma|_K$.
    \end{itemize}
\end{lemma}
In the proof we use a variant of the barycenter called the \emph{parametrized barycenter}. Given a measure $\P\in \mathcal M(\Fr(X))$, define a measure $\P_b$ on $X$ by 
\[
\P_b(E)=\int_\Gamma\int_{\dom(\gamma)}\chi_E(\gamma_t)\ud t\ud\P(\gamma),\quad E\subset X\textrm{ Borel}.
\]
The reason for this is that the integration over curve fragments with respect to $\ud t$ instead of $\ud s$ ensures suitable upper semicontinuity properties for the functional $\Phi$ defined in the proof below. Note moreover that $\P^\#\ll\P_b\ll\P^\#$.

\begin{proof}
    Suppose (2) does not hold. Then for each $\P\in\mathcal P(\Gamma)$ the barycenter $\P^\#$ of $\P$ is singular with respect to $\mu$.  
    Indeed, if there existed a set $K\subset X$ with $\mu(K)>0$ (which could be taken compact by the inner regularity of $\mu$) and $\P\in\mathcal P(\Gamma)$ such that $0\ne \P^\#|_K\in L^1(\mu)$, then $\P|_K\in\mathcal P(\Gamma|_K)$ given by
    \[
    \P|_K=(R_K)_\ast \P,\text{ where } R_K(\gamma)=\gamma|_{\gamma^{-1}(K)},
    \]  
would yield a 1-plan, because
\begin{align*}
(\P|_K)^\#(E)=\int\int_\gamma\chi_E\ud s\ud\P|_K=\int\int_{\gamma|_{\gamma^{-1}(K)}}\chi_E \ud s\ud\P=\P^\#|_K(E)
\end{align*}
for any Borel set $E\subset X$.

To complete the proof we follow the argument in \cite[Lemma 9.4.3]{rud} to show that if $\P^\#\perp\mu$ for all $\P\in\mathcal P(\Gamma)$, then (1) holds. Since $\P_b$ and $\P^\#$ are mutually absolutely continuous, it is equivalent to show that $\P_b\perp\mu$ for all $\P\in\mathcal P(\Gamma)$ implies (1).  Let $\Gamma_L=\{\gamma\in\Gamma:\biLip(\gamma)\le L\}$ for $L\in \N$. It suffices to show that there exists a $\mu$-null set $N_L\subset X$ so that $\Ha^1|_{\im(\gamma)}$ is concentrated on $N_L$ for every $\Gamma_L$, because then $N=\bigcup_{L=1}^\infty N_L$ is the set we search for. 

Let $L\in \N$. Equip $\mathcal P(\Gamma_L)$ with the $w^\ast$-topology and $G=C(X;[0,1])$ with the topology of uniform convergence, and define the functional $\Phi:\mathcal P(\Gamma_L)\times G\to \R$ by
    \begin{align*}
    \Phi(\P,g):=\int(1-g)\ud\mu+\int g\ud\P_b.
    \end{align*}
 $\Phi$ is bi-linear and continuous in $g$ for any fixed $\P\in \mathcal P(\Gamma_L)$. Since $\gamma\mapsto \int_{\dom(\gamma)}g(\gamma_t)\ud t$ is upper semicontinuous for any $g\in G$, cf. \cite[Proof of Lemma 2.28]{sch16b}, it follows that $\Phi(\cdot,g)$ is upper semicontinuous for any fixed $g\in G$. As  $\mathcal P(\Gamma_L)$ is compact, we may apply Sion's minimax theorem (see \cite{sion58} and \cite[Theorem 4.7]{davideb20}) to obtain
 \begin{equation}\label{eq:minimax}
     \sup_\P\inf_g\Phi(\P,g)=\inf_g\sup_\P\Phi(\P,g).
 \end{equation}
 For every $\P\in \mathcal P(\Gamma_L)$ the fact that $\P_b\perp\mu$ implies $\displaystyle \inf_g\Phi(\P,g)=0$ by an application of Urysohn's lemma. By \eqref{eq:minimax} there are functions $g_j\in G$, for $j\in \N$, such that $\displaystyle \sup_\P\Phi(\P,g_j)<2^{-j}$. Set $N_L=\{x\in X:\ \sum_{j\in\N} g_j(x)<\infty\}$. For each $\P\in\mathcal P(\Gamma_L)$
 \begin{equation}\label{eq:rudinsarg}
0\le \sum_{j\in\N}\Phi(\P,g_j)=\int\bigg(\sum_{j\in\N}[1-g_j]\bigg)\ud\mu+\int\bigg(\sum_{j\in\N}g_j\bigg)\ud\P_b\le 1.
 \end{equation}
 Since $\sum_{j\in\N}(1-g_j)=\infty$ on $N_L$, \eqref{eq:rudinsarg} implies that $\mu(N_L)=0$ and moreover that $\P_b$ is concentrated on $N_L$. The set $N_L$ satisfies the claim in (1) since for any $\gamma\in \Gamma_L$ we may choose $\P=\delta_{\gamma}$ to obtain that $\Ha^1|_{\im(\gamma)}\simeq_L\P_b$ is concentrated on $N_L$.
\end{proof}

Recall that a Borel function $\rho:X\to [0,\infty)$ is a weak $\ast$-upper gradient of $f\in \LIP(X)$ if and only if, for $\Mod_\infty$-a.e. $\gamma\in \Fr(X)$, we have that $(f\circ\gamma)'(t)\le \rho(\gamma_t)|\gamma_t'|$ a.e. $t\in \dom(\gamma)$, cf. the discussion before the proof of Proposition \ref{prop:min-star-ug-exists}. In the next proposition we apply Lemma \ref{lem:rainwater} to the family of curve-fragments failing this condition. Namely, we define
\[
\Gamma_\rho:=\{\gamma\in\Fr(X):\ (f\circ\gamma)'(t)> \rho(\gamma_t)|\gamma_t'|\textrm{ on a positive measure subset of }\dom(\gamma)\}.
\]
Notice that $\rho$ fails to be a weak $\ast$-upper gradient of $f$ if and only if $\Mod_\infty\Gamma_\rho>0$.

\begin{proposition}\label{prop:rho-bad}
Let $f\in \LIP(X)$, $\rho\in L^\infty(\mu)$, and suppose that $\rho$ is not a weak $\ast$-upper gradient of $f$, i.e. that we have $\mu(\{\rho<|Df|_\ast\})>0$, or that, equivalently, $\Mod_\infty\Gamma_\rho>0$. Then there exists a 1-plan concentrated on $\Gamma_\rho$.

\end{proposition}

\begin{proof}

Let $(K_j)$ be a sequence of compact sets in $X$ such that $K_j\subset K_{j+1}$, $N=X\setminus \bigcup_jK_j$ is $\mu$-null, and moreover $\rho|_{K_j}$ is continuous for each $j$. Given $j,n\in\N$, we define $\Gamma_{n,j}$ as the family of curve fragments $\gamma:\dom(\gamma)\to K_j$ for which 
\[
|f(\gamma_s)-f(\gamma_t)|\ge (1+1/n)\int_{\gamma|_{\dom(\gamma)\cap [t,s]}}^*\rho\ud s=(1+1/n)\int_{\bar\gamma|_{[t,s]}}\bar \rho\ud s
\]
for all $t,s\in \dom(\gamma)$. Note that for each $n,j\in\N$, the family $\Gamma_{n,j}$ is compact by the lower semicontinuity of the $\ast$-integral, cf. Proposition \ref{prop:liminf}. We first claim that every fragment in $\Gamma_\rho\setminus \Gamma_N^+$ contains a sub-fragment in $\Gamma_{n,j}$ for some $j,n\in\N$.
Let $\gamma\in \Gamma_\rho\setminus\Gamma_N^+$. By a standard exhaustion argument there exists $\lambda>1$ and a (compact) $A\subset \dom(\gamma)$ with $|A|>0$ such that $$(f\circ\gamma)'_t\ge \lambda\rho(\gamma_t)|\gamma_t'|$$ for all $ t\in A$. Since $0=|A\cap \gamma\inv(N)|=\lim_{j\to\infty}|A\setminus \gamma\inv(K_j)|$, we may assume that $\im(\gamma)\subset K_j$ for some suitably large $j$ by considering the restriction $\gamma|_{\gamma\inv(K_j)}$. If $\rho(\gamma_t)|\gamma_t'|=0$ for a.e. $t\in A$ then $\gamma\in \Gamma_{n,j}$ for all $n$. Otherwise we may pass to a compact subset $C\subset A$ with $|C|>0$ such that
\begin{itemize}
\item $t\mapsto \rho(\gamma_t)|\gamma_t'|$ and $t\mapsto (f\circ\gamma)_t'$ restricted to $C$ are continuous and bounded below away from zero;
\item $t\mapsto \frac 1h\int_{\bar\gamma|_{[t,t+h]}}\bar\rho\ud s$  and $t \mapsto \frac{\bar f(\bar\gamma_{t+h})-\bar f(\bar\gamma_t)}{h}$ converge uniformly in $C$ to $\rho(\gamma_t)|\gamma_t'|$ and $(f\circ\gamma)_t'$, respectively, as $h\to 0$.
\end{itemize}
From these it follows that there exists $h_0>0$ such that, for some $0<\lambda'<\lambda$, we have 
\[
\lambda'\frac 1h \int_{\bar\gamma|_{[t,t+h]}}\bar\rho\ud s \le \frac{\bar f(\bar \gamma_{t+h})-\bar f(\bar \gamma_t)}{h}
\]
whenever $t\in C$ and $h\in (0,h_0]$. Choosing $t_0\in C$ such that $|C\cap [t_0,t_0+h_0]|>0$ and $n\in \N$ such that $1+1/n<\lambda'$ we obtain that the restriction of $\gamma$ to $C\cap [t_0,t_0+h_0]$ belongs to $\Gamma_{n,j}$, proving the first claim.

Since every fragment in $\Gamma_\rho\setminus\Gamma_N^+$ contains a subfragment in $\displaystyle \bigcup_{n,j}\Gamma_{n,j}$, we have that $$\Mod_\infty(\Gamma_\rho)=\Mod_\infty(\Gamma_\rho\setminus\Gamma_N^+)\le \Mod_\infty(\bigcup_{n,j}\Gamma_{n,j})\le \sum_{n,j}\Mod_\infty(\Gamma_{n,j}).$$ Thus $\Mod_\infty(\Gamma_{n,j})>0$ for some $n,j\in \N$. In particular the measures $\Ha^1|_{\im(\gamma)}$, $\gamma\in \Gamma_{n,j}$, cannot be all concentrated on a $\mu$-null set. Since $\Gamma_{n,j}$ is compact, Rainwater's lemma (Lemma \ref{lem:rainwater}) yields the existence of a compact set $K\subset K_j$ and a 1-plan $\P$ concentrated on $\Gamma_{n,j}|_K$. This implies the claim since $\Gamma_{n,j}|_K\subset\Gamma_{n,j}\subset \Gamma_\rho$.
\end{proof}

\subsection{Proof of Theorem \ref{thm:min_star_ug}}
In order to prove Theorem \ref{thm:min_star_ug} recall that, given a 1-plan $\P$, we set $\bar\P=\ud s\ud\P=|\gamma_t'|\ud t\ud\P$. Note that $\P^\#=e_\ast\bar\P$ and $\rho_\mathbb{\P^\#}=\frac{d\P^\#}{d\mu}$ is its density. We define
\begin{align*}
    |Df|_\P(x):=\chi_{\{\rho_\mathbb{\P^\#}>0\}}(x)\left\|\frac{(f\circ\gamma)_t'}{|\gamma_t'|}\right\|_{L^\infty(\bar\P_x)},
\end{align*}
where $\{\bar\P_x\}_x$ denotes the disintegration of $\bar\P$ with respect to $e$. The measurability, and a.e. definedness of the expression $(\varphi \circ \gamma)_t'$ follows from the disintegration theorem and standard arguments, see e.g. \cite[Corollary A.3]{teriseb}. These, and similar measurability arguments, are also used later in this paper without mention. Since $\rho_{\P^\#}\in L^1(\mu)$, $\bar\P_x$ is defined for $\mu$-a.e. $x\in \{\rho_{\P^\#}>0\}$, and consequently $|Df|_\P$ is well-defined $\mu$-a.e. on $X$. We record the following properties, whose proofs are similar to \cite[Lemma 3.2]{teriseb} and are left to the interested reader.

\begin{lemma}\label{lem:plan-ug-properties}
Let $f\in \LIP(X)$ and $\P$ a 1-plan. Then $|Df|_\P$ has the following properties.
\begin{itemize}
    \item[(1)] $|(f\circ\gamma)_t'|\le |Df|_\P(\gamma_t)|\gamma_t'|$ a.e. $t\in \dom(\gamma)$ for $\P$-a.e. $\gamma$;
    \item[(2)] $|Df|_\P\le |Df|_\ast$ $\mu$-a.e. 
    \item[(3)] If $\tilde\P$ is a 1-plan with $\tilde\P\ll\P$, then $|Df|_{\tilde\P}\le |Df|_\P$.
\end{itemize}
\end{lemma}

\begin{proof}[Proof of Theorem \ref{thm:min_star_ug}]
    Define the function $\rho=\esssup_{\P}|Df|_\P$, where the measure theoretic supremum is taken over all 1-plans $\P$. We have $\rho\le |Df|_\ast$ $\mu$-a.e by Lemma \ref{lem:plan-ug-properties}(2). We claim that $\rho=|Df|_\ast$. Suppose to the contrary that the set $\{\rho<|Df|_\ast\}$ has positive $\mu$-measure. 
 
    By Proposition \ref{prop:rho-bad} there exists a 1-plan $\tilde \P$ concentrated on $\Gamma_\rho$. By the definition of $\Gamma_\rho$, (1) in Lemma \ref{lem:plan-ug-properties} and the properties of disintegration we have that $|Df|_{\tilde\P}>\rho$ on a set of positive measure, contradicting the definition of $\rho$. Arguing as in \cite[Lemma 3.2]{teriseb} we obtain the existence of a 1-plan $\P$ such that $\rho=|Df|_\P$. This completes the proof of Theorem \ref{thm:min_star_ug}.

\end{proof}

We finish this section by recording the following vector valued analogue of Theorem \ref{thm:min_star_ug}. The two propositions below are in the spirit of \cite[Propositions 4.1--4.3]{teriseb} and the proofs therein establish their claims, except for (iii), which is special for the $p=\infty$ case.  We therefore omit the details except for this  claim. 

\begin{proposition}\label{prop:canonical rep}
	Let $\varphi:X\to \R^n$ be a Lipschitz map. Then there exists a Borel set $D\subset X$ and a 1-plan $\P$ with $\mu|_D\ll \P^\#$ and disintegration $\{\bar\P_x\}$ so that the map $\Phi:(\R^n)^*\times X\to [0,\infty)$ given by 
	\[
	\Phi(\lambda,x)=\chi_D(x)\left\|\frac{\lambda((\varphi\circ\gamma)_t')}{|\gamma_t'|}\right\|_{L^\infty(\bar\P_x)}
	\]
	is Borel and has the following properties.
	\begin{itemize}
		\item[(i)] For each $\lambda\in (\R^n)^*$, $\Phi_\lambda:=\Phi(\lambda,\cdot)$ is a Borel representative of $|D(\lambda\circ\varphi)|_\ast$;
		\item[(ii)] For $\mu$-a.e. $x\in X$, $\Phi^x:=\Phi(\cdot,x)$ is a seminorm on $(\R^n)^*$;
        \item[(iii)] there exists a null set $E$ so that $\Phi(\lambda, \cdot) + \LIP(\varphi) \chi_E$ is a $\ast$-upper gradient for $\lambda \circ \varphi$ for every $\lambda \in (\R^n)^*$ and any curve fragment $\gamma$. 
	\end{itemize}
\end{proposition}
\begin{proof}
We only prove (iii) as its proof is distinct from those in \cite{teriseb}. First, let $S\subset (\R^n)^*$ be a countable dense subset. Then by (i) and Remark \ref{rmk:genuine-star-ug} there exists for every $s\in S$ a set $E_s$ so that $\Phi(s, \cdot) + \LIP(\varphi) \chi_{E_s}$ is a $\ast$-upper gradient for $s \circ \varphi$. If we set $E=\bigcup_{s\in S} E_s$, we observe that $\Phi(s, \cdot) + \LIP(\varphi) \chi_{E}$ is a $\ast$-upper gradient for $s \circ \varphi$ for all $s\in S$. Since $\lambda \mapsto \Phi(\lambda, x)$ is Lipschitz continuous, and since $\|s-s'\|\LIP(\varphi)$ is a $\ast$-upper gradient for $(s-s')\circ \varphi$,
we obtain directly that $\Phi(\lambda, \cdot) + \LIP(\varphi) \chi_{E}$ is an $\ast$-upper gradient for all $\lambda \circ \varphi$. 
\end{proof}
We note that in the previous proposition, $D$ may be a null-set itself, which occurs if $|D(\lambda \circ \varphi)|_*=0$ a.e. for every $\lambda \in (\R^n)^*$. Let $\lambda \mapsto \|\lambda\|_*$ denote the dual norm to the standard Euclidean norm on $\R^n$. In the next statement, we define
\begin{align*}
	I(\varphi,\cdot):= \underset{\|\lambda\|_*=1}\essinf \, |D(\lambda\circ\varphi)|_\ast.
\end{align*}
This is the measure theoretic essential infimum of $|D(\lambda\circ\varphi)|_\ast$.

\begin{proposition}\label{prop:useful_properties}
	$\Phi$ has the following useful properties.
	\begin{itemize}
		\item[(i)] For any Borel map $\bm\lambda:X\to (\R^n)^*$ we have that $\Phi(\bm\lambda_x,x)=0$ $\mu$-a.e. $x\in X$ if and only if $\bm\lambda_{\gamma_t}((\varphi\circ\gamma)_t')=0$ a.e. $t\in \dom(\gamma)$ for $\infty$-a.e. $\gamma$;
		\item[(ii)] $\displaystyle I(\varphi,x)=\underset{\|\lambda\|_*=1}{\inf}\, \Phi(\lambda,x)$ for $\mu$-a.e. $x\in X$.
	\end{itemize}
\end{proposition}

\section{Fragment-wise differentiable structure}

Throughout this section we fix a Borel set $U\subset X$ with $\mu(U)>0$ and $\varphi\in \LIP(X,\R^n)$, $n\ge 1$, and let $\Phi$ be as in Proposition \ref{prop:canonical rep}.

\subsection{Charts and differentials}

For what follows we define the Borel map ${\rm Dir}_\varphi:\overline\Fr(X)\to S^{n-1}$ by ${\rm Dir}_\varphi(\gamma,t)=\frac{(\varphi\circ\gamma)'(t)}{|(\varphi\circ\gamma)'(t)|}$ whenever $(\varphi\circ\gamma)'(t)$ exists and is non-zero, and ${\rm Dir}_\varphi(\gamma,t)=e_1$ otherwise. 

\begin{lemma}\label{lem:ess-span}
	Let $\P'$ be a 1-plan with $U\subset \{\rho_{(\P')^\#}>0\}$ and set $$\Phi'(\lambda,x):=\left\|\frac{\lambda((\varphi\circ\gamma)'(t))}{|\gamma_t'|}\right\|_{L^\infty(\bar\P'_x)},\quad I'(\varphi,x):=\inf_{\|\lambda\|_*=1}\Phi'(\lambda,x)$$ for $\mu$-a.e. $x\in U$. Let

\[
G:=\Big\{(\gamma,t):\ \gamma_t\in U,\ \frac{| (\varphi\circ\gamma)'(t)|}{|\gamma_t'|}\ge \frac{I'(\varphi,\gamma_t)}{2} \Big\}\subset\overline\Fr(X).
\]

Then for $\mu$-a.e. $x\in U$ we have that $I'(\varphi,x)>0$ if and only if ${\rm Dir}_{\varphi\ast}(\chi_G\bar\P'_x)\in \mathcal P(S^{n-1})$ is supported on a set spanning $\R^n$.
\end{lemma}
In particular if $I'(\varphi,\cdot)>0$ $\mu$-a.e. on $U$, then $\spt({\rm Dir}_{\varphi\ast}(\chi_G\bar\P_x))$ spans $\R^n$ for $\mu$-a.e. $x\in U$. Here ${\rm Dir}_{\varphi\ast}(\chi_G\bar\P_x)$ denotes the push-forward of the measure $\chi_G\bar\P_x$ by the map ${\rm Dir}_{\varphi}$.
\begin{proof}
	Let $x\in U$ be such that $\bar\P'_x$ exists and $\Phi'(\cdot,x)$ is a seminorm. Since $\varphi \circ \gamma$ is differentiable for a.e. $t\in {\rm dom}(\gamma)$, by disintegration, this occurs for a.e. $x\in U$. Let $V_x\subset \R^n$ be the subspace spanned by $\spt({\rm Dir}_{\varphi\ast}\chi_G\bar\P'_x)$. We have that $(\varphi\circ\gamma)_t'\in V_x$ for $\bar\P'_x$-a.e. $(\gamma,t)\in\overline\Fr(X)\cap G$, since $\spt({\rm Dir}_{\varphi\ast}\chi_G\bar\P'_x)$ equals the $\bar\P'_x$-essential image of ${\rm Dir}_\varphi(G)$. 
If $V_x\ne \R^n$, there exists $\lambda\in (\R^n)^*$ such that $|\lambda|^*=1$ and $\lambda|_{V_x}=0$ implying that either 
$| \lambda(\varphi\circ\gamma)'(t)|=0$ or $(\gamma,t)\not\in G$ for $\P'_x$-a.e. $(\gamma,t)\in \overline\Fr(X)$. In either case, we have $|\lambda((\varphi\circ\gamma)'(t))|\leq |(\varphi \circ \gamma)'(t)|\leq \frac{I'(\varphi,\gamma_t)}{2}|\gamma_t'|$ (or the derivative is not defined). Thus
\[
\left\|\frac{\lambda((\varphi\circ\gamma)'(t))}{|\gamma_t'|}\right\|_{L^\infty(\bar\P'_x)}\leq I'(\varphi,x)/2,
\]
which implies $I'(\varphi,x)\le I'(\varphi,x)/2$, i.e. $I'(\varphi,x)=0$.

Conversely, if $V_x=\R^n$ then for each $\lambda\in (\R^n)^*\setminus\{0\}$, $|\lambda((\varphi\circ\gamma)'(t))|>0$ for a set of $(\gamma,t)\in G$'s of positive $\bar\P'_x$-measure, and thus $\Phi'(\lambda,x)>0$. Since $\lambda\mapsto \Phi'(\lambda,x)$ is continuous it follows that $I'(\varphi,x)>0$.
\end{proof}

For the statement of the next proposition, denote $\displaystyle I(\varphi,x):=\inf_{\|\lambda\|_\ast=1}\Phi(\lambda,x)$, cf. Proposition \ref{prop:useful_properties}(ii).

\begin{proposition}\label{prop:indep-vs-Alberti-rep}
We have that $I(\varphi,\cdot)>0$ $\mu$-a.e. on $U$ if and only if there exists a Borel function $\bm\delta:U\to (0,\infty)$ and $\varphi$-independent Alberti representations $\mathcal A_j=\{\mu_\gamma^j,\P^j\}$ of $\mu|_U$ such that
\begin{align} \label{eq:speed-lower-bd}
	|(\varphi\circ\gamma)'(t)|\ge \bm\delta_{\gamma_t}|\gamma_t'|\quad  a.e.\ t\in \dom(\gamma)
\end{align}
$\P^j$-a.e. $\gamma$, for each $j=1,\ldots,n$. 
\end{proposition}

\begin{proof}  
Let $\P$ be the 1-plan in Proposition \ref{prop:canonical rep}. We assume that  $I(\varphi,\cdot)>0$ $\mu$-a.e. on $U$, and proceed to prove the existence of the Alberti representations. By Lemma \ref{lem:ess-span}, when 
\[
G=\Big\{(\gamma,t):\ \gamma_t\in U,\ \frac{| (\varphi\circ\gamma)'(t)|}{|\gamma_t'|}\ge \frac{I(\varphi,\gamma_t)}{2} \Big\}\subset\overline\Fr(X).
\]
we have that $\spt({\rm Dir}_{\varphi\ast}(\chi_G\bar\P_x))=\R^n$ for $\mu$-a.e. $x\in U$.  We also have
\begin{align}\label{eq:speed-lower-bd-plan}
|(\varphi\circ\gamma)'(t)|\ge \bm\delta_{\gamma_t}|\gamma_t'|\ a.e.\ t\in\dom(\gamma)\ \textrm{ for every } (\gamma,t)\in G.
\end{align}
where $\bm\delta=I(\varphi,\cdot)/2$.  

Now, let $S$ be a countable dense set in $S^{n-1}$ and define
\begin{align}
\mathcal{C}=\{((\xi_1,\dots, \xi_n), \vartheta)\ :\  & \xi_i \in S, \vartheta\in \Q \cap (0,\pi/2)  \nonumber \\
& \text{ and the cones } C(\xi_i, \vartheta) \text{ are independent for } i=1, \dots, n\}.
\end{align}
The set $\mathcal{C}$ is countable and we can enumerate its elements as $\mathcal{C} = \{((\xi_1^j,\dots, \xi_n^j), \vartheta^j) : j\in \N\}$. For each $j\in \N$, let $C_j$ be the set of $x\in X$ so that $\bar\P_x(\{(\gamma,t)\in G: (\phi \circ \gamma)'(t) \in C(\xi_i^j,\vartheta^j))\}>0$ for every $i=1, \dots, n$. The sets $C_j$ are Borel, by the measurability condition in the definition of the disintegration; recall the definition at the beginning of Section \ref{sec:repgradient}. 

Next, we will perform a measurable selection argument. This, and the similar arguments below, could alternatively be accomplished by quoting an appropriate measurable selection theorem, see \cite[Theorem 18.1]{kechris}. Instead, we will be explicit with the rather simple choices that we need. 
Since $\spt({\rm Dir}_{\varphi\ast}(\chi_G\bar\P_x))$ spans $\R^n$ for $\mu$-a.e. $x\in X$, we have that for $\mu$-a.e. $x\in U$, there exists a $j\in \N$ so that $x\in C_j$. Now, let $\tilde{C}_j = C_j \setminus \bigcup_{i<j} C_i$. Then, $\mu(A\setminus \bigcup_{j\in \N} \tilde{C}_j)=0$. Define the following Borel maps $\bm v_1,\ldots,\bm v_n:U\to S^{n-1}$ and $\bm\vartheta:U\to (0,\pi/2)$: $\bm v_i(x)=\xi_i^j$  and $\bm\vartheta(x)=\vartheta^j$ for $x\in \tilde{C}_j$. For $x\in U\setminus \bigcup_{j\in \N} \tilde{C}_j$ define $\bm v_i$ and $\bm\theta$ as an arbitrary constant. Now, by construction the cone fields $\bm C_j=C(\bm v_j,\bm\vartheta)$ are $\varphi$-independent
and $\bar\P_x({\rm Dir}_\varphi\inv\bm C_j(x)\cap G)>0$ $\mu$-a.e. $x\in U$ for each $j=1,\ldots,n$. 

We will conclude the proof by using \cite[Corollary 5.4]{bate15}. By this result, for every $j=1,\dots, n$ we can decompose $U$ as $U=A_j\cup S_j$, with $\mu|_{A_j}$ admitting an Alberti representation $\mathcal A_j$ in the cone direction $\bm C_j$ which satisfies \eqref{eq:speed-lower-bd}, and $S_j$ satisfies the following. For every $\gamma\in \Fr(X)$ which is in $\varphi$-direction $\bm C_j$ and which satisfies \eqref{eq:speed-lower-bd} we have $\Ha^1(\gamma \cap S_j)=0$. 

 Denote
\[
G_j={\rm Dir}_\varphi\inv\bm C_j \cap G=\{(\gamma,t)\in G:\ {\rm Dir}_\varphi(\gamma,t)\in \bm C_j(\gamma_t)\}\subset \overline\Fr(X).
\]
We have $\bar\P_x(G_j)>0$ for a.e. $x\in S_j$, since $S_j\subset U$. Note that $G_j$ is a Borel set and thus $x\mapsto \bar \P_x(G_j)$ is Borel by the properties of disintegration, cf. (2) in the beginning of Section \ref{sec:repgradient}. If $\mu(S_j)>0$, then from  Proposition \ref{prop:canonical rep}  and the disintegration identity we get
\begin{equation}\label{eq:Ginteg}
0<\int_{S_j}  \bar\P_x(G_j) \ud \P^\# =\int \int_{\gamma} \chi_{G_j}(\gamma,t) \chi_{S_j}(\gamma_t) ds \ud\P.
\end{equation}
Now, if $\gamma$ is any curve fragment, we claim that for a.e. $t\in \dom(\gamma)$ we have $(\gamma,t)\not\in G_j$ or $\gamma_t\not\in S_j$. First, let $K\subset \dom(\gamma)$ be a compact set, where $(\gamma,t)\in G_j$ and $\gamma_t\in S_j$. Since $(\gamma,t)\in G_j$ for all $t\in K$, we have that $\gamma|_K$ is in $\varphi$-direction of $\bm C_j$ and it satisfies \eqref{eq:speed-lower-bd}. Thus,  $\Ha^1(\gamma \cap S_j)=0$ by the previous paragraph, and since $\gamma$ is biLipschtiz we get $\leb^n(K)=0$. Thus, for a.e. $t\in \dom(\gamma)$ we have $(\gamma,t)\not\in G_j$ or $\gamma_t\not\in S_j$. This contradicts \eqref{eq:Ginteg}, since the previous implies that the integrand in the right hand side vanishes and thus
\[
\int \int_{\gamma} \chi_{G_j}(\gamma,t) \chi_{S_j}(\gamma_t) ds \ud\P=0.
\]
This concludes the first part of the proof.

Next assume that $\mathcal A_j=\{\mu_\gamma^j,\P^j\}$, $j=1,\ldots,n$ are $\varphi$-independent and satisfy \eqref{eq:speed-lower-bd} for some Borel function $\bm\delta:U\to (0,\infty)$. Let $N=\{x\in U: I(\varphi,x)=0\}$. We will show that $\mu(N)=0$. 

By definition, there exists independent cone fields $\bm C_j = C(\bm v_j, \bm \vartheta)$, so that $(\varphi \circ \gamma)_t' \in \bm C_{j,\gamma_t}$ for a.e. $t\in \dom(\gamma)$ and $\P^j$-a.e. $\gamma \in \Fr(X)$ for $j=1,\dots, n$. Since the cone fields $\bm C_j$ are independent, there exists a  positive function $\bm \varepsilon$ so that for a.e. $x\in U$ we have
\begin{equation}\label{eq:lambda}
\bm \varepsilon_x\|\lambda\|_{*} \leq \sum_{j=1}^n \inf_{w \in \bm C_{j,x}\cap S^{n-1}}|\lambda(w)| 
\end{equation}
for every $\lambda \in (\R^n)^*$. Arguing as in the construction of $\bm v_i$ and $\bm \theta$, we may find a Borel function $\bm\lambda:N\to (\R^n)^*$ so that $\|\bm\lambda\|_*=1$, and $\Phi(\bm\lambda,x)<\bm\delta \bm\varepsilon/(2n)$ for a.e. $x\in N$.  

We obtain from \eqref{eq:lambda} that for a.e. $x$ there exists a $j(x)$ so that 
\begin{equation}\label{eq:lambda''}
\bm \varepsilon/n|w| \leq |\bm\lambda(w)|
\end{equation}
for all unit vectors $w\in \bm C_j$. Then, by homogeneity this holds for all $w\in \bm C_j$. Now let $N_j = \{x\in N: j(x)=j\}$, and observe that $\mu(N)\leq \sum_{j=1}^n \mu(N_j)$. We claim that $\mu(N_j)=0$ for all $j=1,\dots, n$. First, by (iii) of Lemma \ref{prop:canonical rep}  there exists a null set $E$ so that $\Phi(\lambda, x)+\LIP(\varphi)\chi_E$ is an $\ast$-upper gradient for $\lambda \circ \varphi$. Thus for every $\gamma\in \Fr(X)$ we have for a.e. $t\in \gamma\inv(N_j\setminus E)$ 
        \begin{align}
        |\bm\lambda(\varphi \circ \gamma)'(t)| \leq & (\Phi(\bm\lambda, \gamma_t) +\chi_E(\gamma_t))|\gamma_t'| \label{eq:upperboundNj}\leq \bm\delta^{-1} \bm \delta \bm\varepsilon/(2n) |(\varphi \circ \gamma)'(t)|\\
        =&\bm\varepsilon/(2n) |(\varphi \circ \gamma)'(t)|.\nonumber
        \end{align}
        Here the second inequality comes from the fact that the curves have positive $\varphi$-speed. Since $\mu(E)=0$, for $\P^j$-a.e. $\gamma$ we have $\mu_{\gamma}^j(E)=0$. However, for $\P^j$ a.e. $\gamma$ and a.e. $t\in \dom(\gamma)$, we have $(\varphi \circ \gamma)_t' \in \bm C_j$ and thus by \eqref{eq:lambda''} we get for a.e. $t\in \gamma\inv(N_j)$ that
\[
|\bm\lambda(\varphi \circ \gamma)'(t)| \geq \bm \varepsilon/ n |(\varphi \circ \gamma)'(t)|.
\]
This estimate can not hold simultaneously with \eqref{eq:upperboundNj}. Thus, for a.e. $t\in \gamma\inv(N_j)$ we have $\gamma(t)\in E$ for $\P^j$-a.e. $\gamma$. Therefore, $\mu_{\gamma}(N_j) = \mu_{\gamma}(E)$ for a.e. $\gamma\in \P^j$.  since $\mathcal A_j$ is an Alberti representation, we conclude $\mu(N_j)=\mu(E)=0$. Since this holds for all $j=1, \dots, n$, we conclude that $\mu(N)=0$.

\end{proof}

We now establish the existence of differentials and pointwise norms associated to fragment-wise charts.

\begin{proof}[Proof of Theorem \ref{thm:diff-and-ptwise-norm}]
	Suppose $(U,\varphi)$ is a fragment-wise chart, and  define $|\cdot|_x=\Phi(\cdot,x)$ for $\mu$-a.e. $x\in U$. In the case $n=0$ thse are intrepreted vacuously. 
    
    First, consider $n\geq 1$. By Proposition \ref{prop:indep-vs-Alberti-rep} we have that $I(\varphi,\cdot)>0$ $\mu$-a.e. on $U$. By Proposition \ref{prop:useful_properties}(i), $|\cdot|_x$ is thus a norm on $(\R^n)^*$ for $\mu$-a.e. $x\in U$. Moreover $x\mapsto |\lambda|_x$ is a Borel representative of $|D(\lambda\circ\varphi)|_\ast$ for each $\lambda\in (\R^n)^*$.
	
	Given $f\in \LIP(X)$, let $\Psi:(\R^{n+1})^*\times X\to [0,\infty)$ be the map from Proposition \ref{prop:canonical rep} associated to $\psi:=(\varphi,f)\in \LIP(X,\R^{n+1})$. In what follows, given $\lambda\in (\R^n)^*$ and $a\in \R$, we use the notation $(\lambda,a)\in (\R^{n+1})^*$ for the linear map $\R^{n+1}\ni (x,x_{n+1})\mapsto \lambda(x)+ax_{n+1}$. It is straightforward to check, by applying Propositions \ref{prop:canonical rep} and \ref{prop:useful_properties} that
	\begin{align*}
		\Psi((\lambda,0),x)=\Phi(\lambda,x),\quad \lambda\in (\R^n)^*
	\end{align*}
	for $\mu$-a.e. $x\in U$. In particular, if $I(\varphi,x)>0$, the seminorm $\Psi^x$ restricts to a norm on the hyperplane $(\R^n)^*\times\{0\}\subset (\R^{n+1})^*$. On the other hand, since $\mu|_V$ does not admit $n+1$ $\psi$-independent Alberti representations for any positive measure $V\subset U$, by Proposition \ref{prop:indep-vs-Alberti-rep} we have that $I(\psi,x)=0$ $\mu$-a.e. $x\in U$. Thus $\Psi^x$ is degenerate for $\mu$-a.e. $x\in U$. Since it is a norm on the hyperplane $(\R^n)^*\times\{0\}$, the kernel $\ker(\Psi^x)$ is 1-dimensional. It follows that for $\mu$-a.e. $x\in U$, there exists a unique $\ud_xf\in (\R^n)^*$ such that $\Psi^x((\ud_xf,-1))=0$, and the map $x\mapsto \ud_xf$ is Borel. To see measurability, observe that whenever $A$ is a closed set, we have $\{x\in U : \ud_x f \in A\}=\{x\in U: \inf_{a\in A} \Psi(a,-1),x)=0\}$. From this together with the continuity of $a\mapsto \Psi(a,-1)$, measurability follows.  By Proposition \ref{prop:useful_properties}(i) this implies that $\ud_xf((\varphi\circ\gamma)_t')-(f\circ\gamma)_t'=0$ a.e. $t\in \gamma\inv(U)$ for $\Mod_\infty$-a.e. $\gamma$, i.e. $\ud f$ is a fragment-wise differential of $f$ with respect to $(U,\varphi)$.  The uniqueness follows, since if $\bm\xi_x$ was another differential, then $(\ud_xf-\bm\xi_x)(\varphi\circ\gamma)_t'=0$ for a.e. $t\in \gamma\inv(U)$ for $\Mod_\infty$-a.e. $\gamma$. Then, $\Phi(\ud_xf-\bm\xi_x,x)=0$ for $\mu$-a.e. $x\in U$ by (i) of Proposition \ref{prop:useful_properties}. This, together with $I(\varphi,x)>0$ implies $\ud_x f = \bm\xi_x$ for a.e. $x\in U$.  Moreover we have 
	\begin{align*}
		||\ud_xf|_x-|Df|(x)||=|\Psi((\ud_xf,0),x)-\Psi((\bar 0_{(\R^n)^*},1),x)|\le \Psi((\ud_xf,-1),x)=0
	\end{align*}
	$\mu$-a.e. $x\in U$, proving \eqref{eq:ptwise-norm}. This completes the proof of the theorem.

   In the $n=0$ case, the second paragraph and maximality give that for every $f\in \LIP(X)$ $\mu$-a.e. $x\in U$ satisfies $|Df|_*(x)=0$. This means that $(f\circ \gamma)'_t\leq |Df|_*(x)|\gamma'_t|$ for $\Mod_\infty$-a.e. $\gamma$ and a.e. $t\in \gamma\inv(U)$. In particular, $\bm\xi=0$ gives a differential and we have $|Df|_*(x)=\|\bm\xi\|_x$ for $\mu$-a.e. $x\in U$. with the trivial norm $\|\cdot\|_x=0$.
\end{proof}

\subsection{Basic properties and existence}\label{sec:basic-prop+exist}
In the following two lemmas we collect some basic properties of fragment-wise charts and differentials, such as the usual calculus rules for the differential.

\begin{lemma}\label{lem:basic-chart}
	Suppose $(U,\varphi)$ is a fragment-wise chart. If $V\subset U$ is Borel and $\mu(V)>0$, then $(V,\varphi)$ is a fragment-wise chart. Likewise if $N\subset X$ $\mu$-null, then $(U\cup N,\varphi)$ is a fragment-wise chart.
\end{lemma}
In the proof, given an Alberti representation $\mathcal A=\{\mu_\gamma,\P\}$, its restriction to a (Borel) set $E\subset X$ is defined as $\mathcal A|_E\defl \{\mu_\gamma|_E,\P|_{\{\gamma:\ \mu_\gamma(E)>0\}}\}$.
\begin{proof}
Any Alberti representation of $\mu|_U$ is an Alberti representation of $\mu|_{U\cup N}$ since $\mu|_U=\mu|_{U\cup N}$, implying the second claim. To see the first, note that the maximality in Definition \ref{def:frag-chart} holds for $(V,\varphi)$ and, if $\mathcal A_1,\ldots,\mathcal A_n$ are $\varphi$-independent Alberti representations of $\mu|_U$, then the restrictions $\mathcal A_j|_V$ form $n$ $\varphi$-independent Alberti representations of $\mu|_V$, proving independence. 
\end{proof}

\begin{lemma}\label{lem:label-diff}
Suppose $(U,\varphi)$ is a fragment-wise chart. Then 
\begin{itemize}
	\item[(1)] $\ud(af+bg)=a\ud f+b\ud g$ for every $f,g\in \LIP(X)$ and $a,b\in \R$ (linearity);
	\item[(2)] $\ud (fg)=f\ud g+g\ud f$ $\mu$-a.e. on $U$ for every $f,g\in \LIP(X)$ (product rule); 
	\item[(3)] $\ud (h\circ f)=(h'\circ f)\ud f$ for every $f\in \LIP(X)$, $h\in C^1(\R)$ (chain rule); and finally 
	\item[(4)] $\ud f=\ud g$ $\mu$-a.e. on $\{f=g\}\cap U$ for every $f,g\in \LIP(X)$ (locality).
\end{itemize}
\end{lemma}
\begin{proof}
Since the differential of a function is uniquely determined by \eqref{eq:diff}, we obtain the following fragment-wise identities (which are obtained by differentiation of single variable Lipschitz functions in the parameter $t$) 
\begin{align*}
	((af+bg)\circ\gamma)'(t)&=a(f\circ\gamma)'(t)+b(g\circ\gamma)'(t)\ a.e.\ t\in \dom(\gamma)\\
	(fg\circ\gamma)'(t)&=f(\gamma_t)(g\circ\gamma)'(t)+g(\gamma_t)(f\circ\gamma)'(t)\ a.e.\ t\in \dom(\gamma)\\
	(h\circ f\circ\gamma)'(t)&=h'(f(\gamma_t))(f\circ\gamma)'(t)\ a.e.\ t\in \dom(\gamma).
\end{align*}
These combined with uniqueness of the differential from Theorem \ref{thm:diff-and-ptwise-norm} imply (1)--(3). Similarly, since $((f-g)\circ\gamma)'(t)=0$ for a.e. $t\in \dom(\gamma)\cap \gamma\inv(\{f=g\})$, we have that $\ud (f-g)((\varphi\circ\gamma)'(t))=0$ for a.e. $t\in \gamma\inv(U\cap \{f=g\})$. Propositition \ref{prop:useful_properties}(i) with $\bm\lambda \defl \chi_{\{f=g\}}\ud(f-g)$ now implies (4).

\end{proof}

\begin{theorem}\label{thm:cover-by-frag-charts}
	Let $X=(X,d,\mu)$ be a metric measure space and suppose $X$ can be covered up to a $\mu$-null set by countably many Borel sets $X_i$ with finite Hausdorff dimension. Then there is a countable collection of fragment-wise charts $(U_j,\varphi_j)$ of dimension $n_j$ such that $X=\bigcup_jU_j$, and moreover $n_j\le \dim_HX_i$ if $\mu(U_j\cap X_i)>0$.
\end{theorem}

\begin{proof}
	Let $V\subset X$ be a Borel set with $\mu(V)>0$. There exists $i\in\N$ such that $\mu(V\cap X_i)>0$. Consider the family $\mathcal C$ of pairs $(U,\varphi)$ where $U\subset V\cap X_i$, $\mu(U)>0$, $\varphi\in \LIP(X,\R^{n_{\varphi}})$ for some $n_\varphi\in\N$ and $\mu|_U$ admits $n_\varphi$ $\varphi$-independent Alberti representations, or $n_\varphi=0$. By \cite[Theorem 5.3]{bkt19} we have $n_\varphi\le \dim_HX_i<\infty$ for every $(U,\varphi)\in \mathcal C$, and thus there exists a pair $(U,\varphi)\in\mathcal C$ attaining the maximum $N:=\sup_{(U,\varphi)\in \mathcal C} n_\varphi\le \dim_HX_i$. Thus $\mu|_U$ admits $N$ $\varphi$-independent Alberti representations and, since $N$ is maximal, $\mu|_W$ cannot admit more than $N$ independent Alberti representations for any $W\subset U$. That is, $(U,\varphi)$ is a fragment-wise chart. 
	
	We have shown that, given any $V\subset X$ with $\mu(V)>0$, there exists a Borel set $U\subset V$ with $\mu(U)>0$, and a map $\varphi\in\LIP(X,\R^n)$ such that $(U,\varphi)$ is a fragment-wise chart with dimension $n\le \dim_HX_i$ whenever $\mu(U\cap X_i)>0$. By  \cite[Proposition 3.1.1]{kei02} there exists a countable decomposition $X=N\cup\bigcup_jU_j$ and maps $\varphi_j\in\LIP(X,\R^{n_j})$ such that $\mu(N)=0$ and $(U_j,\varphi_j)$ are fragment-wise charts. This readily implies the claim.
\end{proof}

\subsection{Weak$^*$ continuity of differentials}
Assume that $(U,\varphi)$ is a fragment-wise chart. Given $p\in [1,\infty]$, we denote by $L^p(T^*(U,\varphi))$ the space of Borel functions $\bm\xi:U\to (\R^n)^*$ up to $\mu$-a.e. equivalence for which $\|\bm\xi\|_p<\infty$, where
\begin{align}\label{eq:Lp-norm}
	\|\bm\xi\|_{p}\defl \|x\mapsto |\bm\xi_x|_x\|_{L^p(\mu|_U)}.
\end{align}
The space $L^p(T^*(U,\varphi))$ equipped with the norm \eqref{eq:Lp-norm} is a Banach space. Denoting by $|\cdot|_x^*$ the dual norm of $|\cdot|_x$ on $\R^n$, we similarly define the space $L^p(T(U,\varphi))$ of ($\mu$-a.e. equivalence classes of) Borel functions $\bm X:U\to \R^n$ with finite norm
\begin{align}
	\|\bm X\|_{p}\defl \|x\mapsto |\bm X(x)|^*_x\|_{L^p(\mu|_U)}.
\end{align}
Note that if $\{W_i\}$ is a countable Borel partition of $U$ then we have an isometric identification 
\begin{align*}
L^1(T(U,\varphi))= \bigoplus_{l^1}L^1(T(W_i,\varphi))\textrm{ and thus }L^1(T(U,\varphi))^*=\bigoplus_{l^\infty}L^1(T(W_i,\varphi))^*,
\end{align*}
where $\bigoplus_{l^\infty}B_i$ consists of sequences $\bar v= (v_i)_i$ such that $v_i\in B_i$ and $\|\bar v\|=\sup_i\|v_i\|_{B_i}<\infty$. Using a partition $W_i$ of $U$ such that $\delta_i|\cdot|_{Eucl}\le |\cdot|_x^*\le \delta_i\inv|\cdot|_{Eucl}$ $\mu$-a.e. on $W_i$ for some $\delta_i>0$, we obtain the isomorphism $L^1(T(W_i,\varphi))\simeq L^1(W_i,\R^n)$, which by standard calculations implies the isometric identification $L^\infty(T^*(W_i,\varphi))=L^1(T(W_i,\varphi))^*$. Using basic functional analytic arguments we deduce from this the following lemma, whose detailed proof we omit.

\begin{lemma}\label{lem:L1-Linfty-duality}
We have the isometric identification $L^\infty(T^*(U,\varphi))=L^1(T(U,\varphi))^*$.
\end{lemma}

Likewise a simple change of variables reduces the following lemma to the well-known weak$^*$-continuity of derivatives of functions of one variable, see \cite[Chapter 7, Example 7.1.5]{wea99} for the statement. For the sake of completeness, we sketch a proof.
\begin{lemma}\label{lem:1-dim-weak-star-conv}
	Suppose $(f_j)\subset \LIP(X)$ converges pointwise to $f\in \LIP(X)$ with $\sup_j\LIP(f_j)<\infty$. Let $\gamma\in \Fr(X)$. Then 
	\begin{align*}
		\lim_{j\to\infty}\int_{\dom(\gamma)} h(t)(f_j\circ\gamma)'(t)\ud t=\int_{\dom(\gamma)} h(t)(f\circ\gamma)'(t)\ud t
	\end{align*}
	for all $h\in L^1(\dom(\gamma))$.
\end{lemma}
\begin{proof}
Let $\overline{\gamma}$ be the extension of $\gamma$ to a curve in $\sh(X)$. Extend $f_j$ and $f$ to Lipschitz functions on $\sh(X)$. When $h=\chi_{[a,b]}$ is a characteristic function of an interval, with $a,b\in \dom(\gamma)$, we have 
    \begin{align*}
    \lim_{j\to\infty}\int_{\dom(\overline{\gamma})} h(t)(f_j\circ\overline{\gamma})'(t)\ud t &= \lim_{j\to\infty}f_j(\gamma(b))-f_j(\gamma(a))=f(\gamma(b))-f(\gamma(a))\\
    &=\int_{\dom(\overline{\gamma})} h(t)(f\circ\gamma)'(t)\ud t.
    \end{align*}
    If $h\in L^1(\dom(\gamma))$, then we can approximate $h$ by a linear combination of functions of the previous type. This yields the claim.
\end{proof}
With these auxiliary results we are ready to prove the weak$^*$-continuity of differentials with respect to $(U,\varphi)$. This statement, and its proof, are quite similar to Lemma \ref{lem:lim-star-ug-weak}.

\begin{theorem}\label{thm:diff-w-star-cont}
	Let $(U,\varphi)$ be a fragment-wise chart. If $(f_j)\subset \LIP(X)$ converges pointwise to $f\in \LIP(X)$ and $\sup_j\LIP(f_j)<\infty$, then $\ud f_j \to \ud f$ weak-$\ast$ in $L^\infty(T^*(U,\varphi))=L^1(T(U,\varphi))^*$. In other words:
	\begin{align*}
		\lim_{j\to\infty}\int_U\ud f_j(\bm X)\ud\mu=\int_U\ud f(\bm X)\ud\mu,\quad \bm X\in L^1(T(U,\varphi)).
	\end{align*}
\end{theorem}

\begin{proof}
    Since we have $\sup_j\LIP(f_j)<\infty$, then $\ud f_j$ belong to some bounded subset of $L^\infty(T^*(U,\varphi))$. By Banach-Alaoglu, since $L^1(T(U,\varphi))$ is separable, a subsequence $\ud f_j$ converges weakly to some $\bm\xi$. We will show that this subsequential weak-$\ast$ limit is the differential of $f$, that is $\bm\xi = \ud f$. Therefore, the weak limit is unique, and thus the full sequence $\ud f_j$ converges weakly to $\bm\xi$. For the rest of the proof, suppose that we have passed to the subsequence so that $\ud f_j$ converges weak-$\ast$ to $\bm\xi$. 

    We first do a standard argument to convert the weak convergence to pointwise almost everywhere convergence. This employs Mazur's lemma applied for a reflexive $L^2$-space. Let $U_n\subset U$ is an increasing sequence of sets, $n\in \N$, with  $\mu(U_n)<\infty$, with $\mu(U\setminus \bigcup_{n\in \N} U_n)=0$. we have $\ud f_j|_{U_n} \in L^2(T^*(U_n,\varphi))$, and $\ud f_j|_{U_n} \to \bm\xi \chi_{U_n}$ weakly in $L^2(T^*(U_n,\varphi))$. Therefore, by Mazur's lemma, there exist values $\alpha_{k,n,m}\geq 0$ and $\sum_{k=m}^{N_{n,m}} \alpha_{k,n,m}=1$ and convex combinations $\ud g_{n,m} = \sum_{k=m}^{N_{n,m}} \alpha_{k,n,m} \ud f_k$, which converge strongly to $\bm\xi|_{U_n}$ in $L^2(T^*(U_n,\varphi))$. Moreover, by Lemma \ref{lem:label-diff} we have that $\ud g_{n,m}$ is the differential of  $g_{n,m} = \sum_{k=m}^{N_{n,m}} \alpha_{k,n} f_k$. Choose now $M_n$ so large that
    \[
    \mu(\{ x \in U_n : |\bm\xi - \ud g_{n,M_n}|_x\geq 2^{-n}\})<2^{-n}.
    \]
    We let $g_n = g_{n,M_n}$. We have $\ud g_n \to \bm \xi$ a.e. (by a standard Borel-Cantelli argument), that $g_n \to f$ pointwise and $\sup_n\LIP(g_n)\leq \sup_j\LIP(f_j)<\infty$.

    Since $\ud g_n \to \bm \xi$ a.e., there exists a null-set $E$ so that for all $x\not\in E$ we have $\ud g_{n,x}\to \bm\xi_x$.  For $\Mod_\infty$-a.e. $\gamma\in \Fr(X)$ we have $(g_n\circ\gamma)'(t)=\ud g_n(\varphi \circ \gamma)'_t$ and $(f\circ\gamma)'(t)\ud t=\ud f(\varphi \circ \gamma)'_t$ for a.e. $t$. Let $\Gamma_B$ be the set of curve fragments where this fails.
    
    Now, if $\gamma \in \Fr(X)\setminus (\Gamma_E^+\cup \Gamma_B)$, we have by dominated convergence and Lemma \ref{lem:1-dim-weak-star-conv} that for all $h\in L^1({\rm dom(\gamma)})$
    \begin{align*}
    \int_{\dom(\gamma)} h(t) \chi_U(\gamma_t)\bm \xi(\varphi \circ \gamma)'_t \ud t &=
    \lim_{n\to\infty}\int_{\dom(\gamma)} h(t) \chi_U(\gamma_t) \ud g_n(\varphi \circ \gamma)'_t \ud t \\
    &=\lim_{n\to\infty}\int_{\dom(\gamma)} h(t) \chi_U(\gamma_t) (g_n\circ\gamma)'(t)\ud t \\
    &=\int_{\dom(\gamma)} h(t) \chi_U(\gamma_t)(f\circ\gamma)'(t)\ud t \\
    &=\int_{\dom(\gamma)} h(t) \chi_U(\gamma_t) \ud f(\varphi \circ \gamma)'_t \ud t \\
    \end{align*}

    Since this holds for all $h\in L^1(\dom(\gamma))$, we conclude $\bm \xi(\varphi \circ \gamma)'_t=\ud f(\varphi \circ \gamma)'_t$ for a.e. $t\in \gamma\inv(U)$. Thus, by $\bm\xi$ is a differential for $f$ and by uniqueness from Theorem \ref{thm:diff-and-ptwise-norm} we have $\bm\xi=\ud f$. This concludes the proof.
\end{proof}

\subsection{Differentiable structure}\label{sec:diff-str} A \emph{measurable $L^\infty$-bundle} $\mathcal{T}$ over $X=(X,d,\mu)$ consists of collections $(\{U_i,V_{i,x}\})_{i\in I}$ and $(\{\phi_{i,j,x}\} )$, where $I$ is countable and
\begin{enumerate}
	\item $U_i\subset X$ are Borel for each $i\in I$, and cover $X$ up to a $\mu$-null set;
	\item for any $i\in I$ and $\mu$-a.e. $x\in U_i$, $V_{i,x}=(V_i,|\cdot|_{i,x})$ is a finite dimensional normed space so that $x\mapsto |v|_{i,x}$ is Borel for any $v\in V_i$;
	\item for any $i,j\in I$ and $\mu$-a.e. $x\in U_i\cap U_j$, $\phi_{i,j,x}\colon V_{i,x} \to V_{j,x}$ is an isometric bijective linear map satisfying the \emph{cocycle condition}: for any $i,j,k \in I$ and $\mu$-a.e. $x \in U_i \cap U_j \cap U_k$, we have $\phi_{j,k,x} \circ \phi_{i,j,x}=\phi_{i,k,x}$.
\end{enumerate}

\begin{remark}
The definition and terminology above are adopted from \cite{che99}, cf. \cite[Section 4]{che99}.
\end{remark}

A \emph{section} of $\mathcal T$ is a collection $\omega=\{\omega_i:U_i\to V_i \}$ of Borel measurable maps with $(\omega_i)_x\in V_{i,x}$ and $\omega_j=\phi_{i,j,x}\circ\omega_i$ $\mu$-a.e. $x$ for all $i,j\in I$. 
We remark that the map $x\mapsto |\omega(x)|_{x}$, given by 
\begin{align}\label{eq:norm-of-section}
|\omega(x)|_x\defl |\omega_i(x)|_{i,x}\ \mu-a.e. \ x\in U_i,
\end{align}
is well-defined up to $\mu$-negligible sets as a consequence of the cocycle condition and the isometry of $\phi_{i,j,x}$. We denote by $\Gamma(\mathcal T)$ the vector space of sections of $\mathcal T$, and by $\Gamma_p(\mathcal T)$ the elements of $\Gamma(\mathcal T)$  such that $x\mapsto |\omega(x)|_x\in L^p(\mu)$. 

\bigskip\noindent A countable collection $(U_i,\varphi_i)$ of fragment-wise charts (of dimension $n_i$) covering $X$ up to a null set gives rise to \emph{measurable tangent and cotangent bundles}, denoted $TX=(\{U_i,(\R^{n_i},|\cdot|_x^*)\})$ and $T^*X=(\{U_i,((\R^{n_i})^*,|\cdot|_x)\})$, respectively. Indeed, (1) holds by assumption, while (2) follows by Theorem \ref{thm:diff-and-ptwise-norm}, and (3) follows from the following lemma, which can be proven exactly as \cite[Proposition 5.1 and Lemma 5.2]{teriseb} . 

\begin{lemma}\label{lem:compatible}
Let $(U_i,\varphi_i)$, $i=1,2$, be fragment-wise charts of dimensions $n_i$. Denote by $|\cdot|_{i,x}$ and $\ud_x^if$ the pointwise norms and the differential of $f\in \LIP(X)$ with respect to $(U_i,\varphi_i)$, respectively, given by Theorem \ref{thm:diff-and-ptwise-norm}, for $\mu$-a.e. $x\in U_i$. If $\mu(U_1\cap U_2)>0$, then $n_1=n_2=:n$ and there exists a linear isometric bijection $\Phi_{1,2,x}:((\R^n)^*, |\cdot|_{1,x})\to ((\R^n)^*, |\cdot|_{2,x})$ such that $\displaystyle \Phi_{1,2,x}(\ud_x^1f)=\ud_x^2f$ $\mu$-a.e. $x\in U_1\cap U_2$ for any $f\in \LIP(X)$. Moreover the collection $(\Phi_{i,j,x})$ satisfies the cocycle condition.
\end{lemma}

It also follows from Lemma \ref{lem:compatible} that, given $f\in \LIP(X)$, the collection $\{ d_{U_i,\varphi_i}f\}$ is a section of $T^*X$, which we denote $\ud f\in \Gamma_\infty(T^*X)$ and call the \emph{fragment-wise differential} of $f$.

\begin{remark}\label{rmk:normed module}
It is possible to equip $\Gamma_p(T^*X)$ and $\Gamma_p(TX)$ with the structure of an $L^p$-normed module, cf. \cite[Section 5.2]{teriseb} and \cite{che99}. It follows from Lemma \ref{lem:L1-Linfty-duality} that $\Gamma_\infty(T^*X)=\Gamma_1(TX)^*$ and from Theorem \ref{thm:diff-w-star-cont} that the fragment-wise differential is w$^*$-continuous in the following sense. Further, from the proof of Theorem \ref{thm:diff-w-star-cont} we get that if $(f_j)\subset \LIP(X)$ converges pointwise to $f\in \LIP(X)$ and $\sup_j\LIP(f_j)<\infty$, then $\ud f_j\stackrel{w^*}{\rightharpoonup} \ud f$ in $\Gamma_\infty(T^*X)$. We will not need this fact in this paper and leave the details to the interested reader.
\end{remark}

\section{Density of Directions}\label{sec:dense-dir}
Recall the definition of the cone $C(w,\vartheta)$ from Section \ref{sec:1-plan+alb-rep}, and notice that $v\notin C(w,\vartheta)$ if and only if $|w\cdot v|\le \cot\vartheta|\pi^\perp_w(v)|$. We will use this observation repeatedly.

\begin{lemma}\label{lem:perp-approx}
Let $w\in S^{n-1}$, $\vartheta\in (0,\pi/2)$, $\varphi\in\LIP(X;\R^n)$ and let $U\subset X$ be a Borel set of positive $\mu$-measure. Suppose there is no 1-plan $\P$ such that $(\varphi\circ\gamma)_t'\in C(w,\vartheta)$ for $\bar\P$-a.e. $(\gamma,t)\in e\inv(U)$. Then there exists a compact set $K\subset U$ of positive $\mu$-measure and, for any $\delta>0$, a sequence $(f_j)\subset \LIP(X)$ with $\LIP(f_j)\le 1+\cot\vartheta\LIP(\varphi)$ satisfying the following properties.
\begin{itemize}
    \item[(i)] $f_j\to w\cdot\varphi$ uniformly in $K$;
    \item[(ii)] for every $j$ there exists $r_j>0$ so that 
  \begin{align}\label{eq:perp_approx}
    |f_j(x)-f_j(y)|\le \delta d(x,y)+\cot\vartheta|\pi_{w}^\perp(\varphi(x)-\varphi(y))|
\end{align}
    for all $x,y\in K$ with $d(x,y)<r_j$.
\end{itemize}
\end{lemma}

\begin{remark}
Functions satisfying an estimate of the type \eqref{eq:perp_approx} appear in \cite{bate15} and \cite{sch16b} where they were used to construct non-differentiable functions.
\end{remark}

\begin{proof}[Proof of Lemma \ref{lem:perp-approx}]
Let $\P$ be a 1-plan and $D\subset X$ a Borel set as in \ref{prop:canonical rep}, i.e. 
\begin{equation}\label{eq:difflam}
|D(\lambda\cdot\varphi)|(x)=\chi_D(x)\left\| \frac{\lambda\cdot(\varphi\circ\gamma)_t'}{|\gamma_t'|}\right\|_{L^\infty(\bar\P_x)},\quad \lambda\in \R^n,
\end{equation}
for $\mu$-a.e. $x\in U$. By assumption and the inner regularity of $\mu$ there exists a compact set $K\subset U$ of positive measure such that 
\begin{align}\label{eq:counterass}
|w\cdot(\varphi\circ\gamma)_t'|\le \cot\vartheta|\pi^\perp_w((\varphi\circ\gamma)_t')| \quad\bar\P_x-\textrm{a.e. }(\gamma,t)
\end{align}
for every $x\in K$. Consider the metric $\tilde d(x,y)=\delta d(x,y)+\cot\vartheta|\pi^\perp_w(\varphi(x)-\varphi(y))|$ on $X$ for $\varepsilon>0$. This metric is bi-Lipschitz to $d$ and thus the classes of  1-plans and curve fragments remain unchanged and in particular the representation above remains valid  for possibly a different plan. Denote the minimal $\ast$-upper gradients and metric differentials with respect to the metric $\tilde d$ by adding a subscript. Suppose $\widetilde \P$ is a 1-plan satisfying \eqref{eq:difflam} for the metric $\tilde d$. Then, by averaging $\P$ and $\widetilde \P$ we have that \eqref{eq:difflam} holds for both metrics $d$ and $\tilde d$, and \eqref{eq:counterass} holds for this new plan as well. We adjust notation and denote the average by $\P$. Note that
\begin{align*}
|\gamma_t'|_{\tilde d}=\delta|\gamma_t'|+\cot\vartheta|\pi^\perp_w((\varphi\circ\gamma)_t')|\quad a.e.\ t\in \dom(\gamma)
\end{align*}
for any $\gamma$. For $\mu$-a.e. $x\in K$, \eqref{eq:counterass} and the choice of plan yields 
\begin{align*}
|D(w\cdot\varphi)|_{\tilde d}(x)&=\left\| \frac{w\cdot(\varphi\circ\gamma)_t'}{|\gamma_t'|_{\tilde d}}\right\|_{L^\infty(\bar\P_x)}=\left\| \frac{w\cdot(\varphi\circ\gamma)_t'}{\delta|\gamma_t'|+\cot\vartheta|\pi^\perp_w((\varphi\circ\gamma)_t')|}\right\|_{L^\infty(\bar\P_x)} \\
&\le \left\| \frac{|w\cdot(\varphi\circ\gamma)_t'|}{\delta|\gamma_t'|+|w\cdot(\varphi\circ\gamma)_t'|}\right\|_{L^\infty(\bar\P_x)}<1.
\end{align*}

By Theorem \ref{thm:approximation} there is a uniformly Lipschitz sequence $(f_j)$ converging pointwise to $w\cdot \varphi$ such that $\Lip_{a,\tilde d}f_j\to |D(w\cdot \varphi)|_{\tilde d}$. Up to passing to a smaller set, and a subsequence, (with positive measure) we have that there exist $r_j\downarrow 0$ so that
\[
|f_j(y)-f_j(z)|\le \delta d(y,z)+\cot\vartheta|\pi^\perp_w(\varphi(y)-\varphi(z))|,\quad y,z\in  K,\ d(y,z)<r_j,
\]
as claimed.
\end{proof}

To derive a contradiction from \eqref{eq:perp_approx} we use the following lemma analogous to \cite[Lemma 9.1]{bate15}.
\begin{lemma}\label{lem:zigzagfrag}
Let $\varphi:X\to \R^n$ be Lipschitz and $U\subset X$ Borel so that there exists a function $\bm\delta:X\to (0,\infty)$ with $\bm\delta>0$ $\mu$-a.e. in $U$ and such that $\mu|_U$ admits $n$ $\varphi$-independent $\xi$-separated Alberti representations $\mathcal{A}_j=(\mu_j, \P_j)$ with $|(\varphi\circ\gamma)_t'|\ge \bm\delta_{\gamma_t}|\gamma_t'|$ for $\P_j$-a.e. $\gamma$ and a.e. $t\in \dom(\gamma)$, with $j=1,\dots,n$. Suppose moreover that $f\in \LIP(X)$ admits a fragment-wise differential $\ud f$ with respect to $(U,\varphi)$. Then, given $w\in S^{n-1}$ and any set $S\subset U$, for $\mu$-a.e. $x\in S$ there exists a sequence $(x_k)\subset S$ converging to $x$ such that $|\varphi(x_k)-\varphi(x)|\ge \bm\delta_x d(x_k,x)$,
\begin{align*}
    &\left|\frac{\varphi(x_k)-\varphi(x)}{d(x_k,x)}-\frac{|\varphi(x_k)-\varphi(x)|}{d(x_k,x)}w\right|\to 0\textrm{ and} \\
    &\left|\frac{f(x_k)-f(x)}{d(x_k,x)}-\ud_xf\bigg(\frac{\varphi(x_k)-\varphi(x)}{d(x_k,x)}\bigg)\right|\to 0.
\end{align*}
\end{lemma}

\begin{proof}
    The proof follows the argument in the proof of \cite[Lemma 9.1]{bate15} and we explain only the modifications needed in our case. 

For each $i\in \{1,\ldots,n\}$, $\varepsilon>0$ and $R>0$, define the set $G^i_R(S)$ to consist of those $x_0\in S$ for which there exists $\gamma\in \Fr(X)$ and $t_0\in \dom(\gamma)$ with $\gamma_{t_0}=x_0$ satisfying (1)-(5) in \cite[Lemma 9.1]{bate15}  and, in addition
\begin{align}
    &(f\circ\gamma)_{t_0}'=\ud_{x_0}f((\varphi\circ\gamma)_{t_0}')\textrm{ and }\label{eq:6}\\
    &|f(\gamma_t)-f(\gamma_{t_0})-(t-t_0)(f\circ\gamma)_{t_0}'|\le \varepsilon|(\varphi\circ\gamma)_{t_0}'|\label{eq:7}
\end{align}
for all $t\in \dom(\gamma)$ with $|t-t_0|\le \frac{5R}{|(\varphi\circ\gamma)_{t_0}'|}$. Since \eqref{eq:6} is satisfied a.e. for any Alberti representation, we conclude as in \cite{bate15} that $G^i_R(S)$ monotonically increases to a full measure set in $S$ as $R\to 0$.

Given $x_0\in G^i_R(S)$, there are $\gamma\in \Fr(X)$ and $t_0\in \dom(\gamma)$ satisfying (1)--(5) in the proof of \cite[Lemma 9.1]{bate15} as well as \eqref{eq:6} and \eqref{eq:7}, and we find $t\in \dom(\gamma)$ such that $x_0=\gamma_{t_0}$ and $x:= \gamma_t$ satisfy \cite[(9.4) and (9.5), p. 472]{bate15}. The additional properties \eqref{eq:6} and \eqref{eq:7} together with \cite[(2), p. 472]{bate15} further imply that
\begin{align}\label{eq:ast}
|f(x)-f(x_0)-\ud_{x_0}f(\varphi(x)-\varphi(x_0))|\le 2(\varepsilon+\LIP(\varphi)\varepsilon)|\varphi(x)-\varphi(x_0)|.
\end{align}
Now proceeding to define the sets $S_R^0=S$, $S_R^{i+1}=G_R^{i+1}(S_R^i)$. We can then for every $x_n\in S^n_{R}$ find a we find as in \cite{bate15} points $x_i\in S^i_R$ satisfying \cite[(9.6) and (9.8)]{bate15} and moreover
\begin{align}
|f(x_n)-f(x_0)- \ud_{x_0}f(\varphi(x_n)-\varphi(x_0))|\le C\varepsilon\sum_{0\le i<n}|\varphi(x_{i+1})-\varphi(x_i)|.
\end{align}
The proof is then concluded by the argument in \cite{bate15}.
\end{proof}

\begin{theorem}\label{thm:density_of_directions}
Let $(U,\varphi)$ be an $n$-dimensional weak chart. Given any $\Gamma_0\subset\Fr(X)$ with $\Mod_\infty(\Gamma_0)=0$ we have that, for $\mu$-a.e. $x\in U$, the collection
\begin{align*}
S_x(\Gamma_0):=\left\{\frac{(\varphi\circ\gamma)_t'}{|(\varphi\circ\gamma)_t'|}\in S^{n-1}: \gamma\notin\Gamma_0,\ \gamma_t=x\textrm{ and $(\varphi\circ\gamma)_t'\neq 0$ exists} \right\}
\end{align*}
is dense in $S^{n-1}$.
\end{theorem}
\begin{proof}
By Lemma \ref{lem:Nexistence} we have $\Gamma_0 \subset \Gamma^+_E$ for some Borel set $E$ with $\mu(E)=0$. This way $\Gamma^+_E$ is a Borel set of curve fragments (since $\gamma \mapsto \int_\gamma \chi_E \ud s$ is Borel for every Borel set $E$ -- this can be seen  by employing an extension to $\sh(X)$ and using results e.g. from \cite[Appendix]{teriseb}). Now, $S_x(\Gamma^+_E)\subset S_x(\Gamma_0)$. We claim that $S_x(\Gamma^+_E)$ is dense for a.e. $x\in U$. Suppose not, then let $S$ be a dense set in $S^{n-1}$. Let $\mathcal{C}=\{C(v,m^{-1}): v\in S, m\in \N\}$. Enumerate the cones in $\mathcal{C}=\{C_j: j\in \N\}$. The set 
\[
N=\{x \in U: S_x(\Gamma^+_E) \text{ is not dense in } \R^n\}\]
is covered by the sets $N_j = \{x\in U : S_x(\Gamma_E^+)\cap C_j=\emptyset\}$.
The sets $N_j$ are again measurable, and it suffices to prove $\mu(N_j)=0$.

Let $w\in S^{n-1},\varepsilon>0$ be such that $C_j=C(w,\varepsilon)$. From Proposition \ref{prop:indep-vs-Alberti-rep} we obtain independent Alberti representations for $N_j$ and so that \eqref{eq:speed-lower-bd} holds for some positive Borel function $\bm\delta:U\to (0,\infty)$. 
Fix $\varepsilon>0$. By Lemma \ref{lem:perp-approx} there is a sequence $(f_j)$ of uniformly Lipschitz functions converging to $w\cdot f$ in $K\subset N_j$ with $\mu(K)>0$  and satisfying \eqref{eq:perp_approx}. Lemma \ref{lem:zigzagfrag} implies that for  $\mu$-a.e. $x\in K$ there exists an $a>0$ and a sequence $(x_k)\subset K$ converging to $x$ such that 
\begin{align*}
&\frac{|\varphi(x_k)-\varphi(x)|}{d(x_k,x)}\to a\in [\bm\delta,\LIP(\varphi)],\ \frac{\varphi(x_k)-\varphi(x)}{d(x_k,x)}\to aw\ \textrm{ and }\\ 
& \limsup_{k\to \infty}\frac{|f_j(x_k)-f_j(x)-\ud_xf_j(\varphi(x_k)-\varphi(x))|}{d(x_k,x)}\to 0
\end{align*}
as $k\to\infty$ for all $j$. It follows from this and \eqref{eq:perp_approx} that $|\ud_x f_j(aw)|\le \varepsilon + \cot\vartheta|\pi_w^\perp(w)|$ and thus $|\ud f_j(w)|\le \varepsilon/\bm\delta$ for $\mu$-a.e. $x\in K$. Since $f_j\to w\cdot \varphi$ pointwise with uniformly bounded Lipschitz constants, Theorem \ref{thm:diff-w-star-cont} implies that $\ud_xf_j(w)\rightharpoonup \ud_x(w\cdot\varphi)(w)=1$ in the weak$^\ast$ sense in $L^\infty(\mu|_U)$, as $j\to \infty$. This yields a contradiction by choosing $\varepsilon$ small enough.
\end{proof}

\section{Comparisons}\label{sec:comparisons}
\subsection{Weaver derivations}

Recall the notion of Weaver derivations from Section \ref{sec:weaver-der}.

\begin{lemma}
	Let $\P$ be a fragment 1-plan and $G:\overline\Fr(X)\to \R$ a bounded Borel function. Given $f\in \LIP(X)$, set
	\begin{align}\label{eq:derivation-from-plan}
		\bm b(f)_x=\chi_{\{\rho_{\P^\#}>0\}}(x)\int G(\gamma,t)\frac{(f\circ\gamma)'(t)}{|\gamma_t'|}\ud\bar\P_x,\quad \mu-a.e.\ x.
	\end{align}
Then \eqref{eq:derivation-from-plan} defines a derivation $\bm b:\LIP(X)\to L^\infty(\mu)$ with $\|\bm b\|_{\mathcal X(\mu)}\le \|G\|_\infty$. 
\end{lemma}
\begin{proof} Recall that $\P^\#=e_\ast\bar\P$ and $\rho_\mathbb{\P^\#}=\frac{d\P^\#}{d\mu}$ is its density. 
Since $|\bm b(f)|\le \|G\|_\infty\LIP(f)$ $\mu$-a.e. we have that $\bm b(f)\in L^\infty(\mu)$ for every $f\in \LIP(X)$ as well as the estimate $\|\bm b\|_{\mathcal X(\mu)}\le \|G\|_\infty$. For any $g\in L^1(\mu)$ and $f\in\LIP(X)$ we have the identity
\begin{align}\label{eq:b-P}
\int g\bm b(f)\ud\mu=&\int_{\{\rho_{\P^\#}>0\}}\int \frac{g(\gamma_t)G(\gamma,t)}{\rho_{\P^\#}(\gamma_t)}\frac{(f\circ\gamma)'(t)}{|\gamma_t'|}\ud\bar\P_x\ud\P^\#\ud\mu\\
=&\int\int_{\dom(\gamma)}\frac{g(\gamma_t)G(\gamma,t)}{\rho_{\P^\#}(\gamma_t)}(f\circ\gamma)'(t)\ud t \ud\P(\gamma)\nonumber.
\end{align}
Suppose $(f_j)\subset \LIP_\infty(X)$ converges pointwise to $f\in \LIP(X)$ and $\sup_j\LIP(f_j)<\infty$.  The estimate 
\begin{align*}
	\int \int_{\dom(\gamma)}\frac{|g(\gamma_t)G(\gamma,t)|}{\rho_{\P^\#}(\gamma_t)}|\gamma_t'|\ud t\ud\P(\gamma)\le \int\int |g G| \ud\bar\P_x\ud\mu(x)\le \|G\|_\infty\|g\|_{L^1(\mu)}
\end{align*}
implies that
\begin{align}
	&\gamma\mapsto \int_{\dom(\gamma)}\frac{g(\gamma_t)G(\gamma,t)}{\rho_{\P^\#}(\gamma_t)}|\gamma_t'|\ud t\in L^1(\P) \label{eq:L1-P},\textrm{ and} \\
	&t\mapsto \frac{g(\gamma_t)G(\gamma,t)}{\rho_{\P^\#}(\gamma_t)}\in L^1(\dom(\gamma)) \label{eq:L1-gamma}\quad \P-a.e.\ \gamma.
\end{align}
By \eqref{eq:L1-gamma} and Lemma \ref{lem:1-dim-weak-star-conv} we have that
\begin{align}\label{eq:conv-P}
	\lim_{j\to\infty}\int_{\dom(\gamma)}\frac{g(\gamma_t)G(\gamma,t)}{\rho_{\P^\#}(\gamma_t)}(f_j\circ\gamma)'(t)\ud t=\int_{\dom(\gamma)}\frac{g(\gamma_t)G(\gamma,t)}{\rho_{\P^\#}(\gamma_t)}(f\circ\gamma)'(t)\ud t
\end{align}
for $\P$-a.e. $\gamma$. The identities \eqref{eq:b-P} and \eqref{eq:conv-P}, together with \eqref{eq:L1-P} and the dominated convergence theorem, yield
\begin{align*}
\lim_{j\to\infty}\int g\bm b(f_j)\ud\mu=\int g\bm b(f)\ud \mu,
\end{align*} 
proving that $\bm b$ is weak$^\ast$-continuous, and thus a Weaver derivation.
\end{proof}

\begin{proposition}\label{prop:normchar}
	Let $f\in \LIP(f)$ and $\bm b\in \mathcal X(\mu)$. Then \[|\bm b(f)|\le |\bm b|_{\mathcal X,\loc}|Df|_\ast\quad\textrm{and}\quad |Df|_\ast=\esssup\{ |\bm b(f)|:\ \|\bm b\|_{\mathcal X(\mu)}\le 1 \} \]
	$\mu$-a.e. on $X$.
\end{proposition}
\begin{proof}
Denote $h=\esssup\{ |\bm b(f)|:\ \|\bm b\|_{\mathcal X(\mu)}\le 1 \}$. The inequality $|Df|_\ast\ge h$ follows from the claimed estimate. To see that $|Df|_\ast\le h$, let $\P$ be the 1-plan given by Theorem \ref{thm:min_star_ug}. Let $R:\Fr(X)\to \Fr(X)$ be the map which reverses each curve fragment. We may assume $R_*\P=\P$ by averaging with the reversion $R_*\P$ of $\P$.Now
	\begin{align*}
		|Df|_\ast(x)=\bigg\|\frac{(f\circ\gamma)'(t)}{|\gamma_t'|}\bigg\|_{L^\infty(\bar\P_x)}\quad \mu-a.e.\ x.
	\end{align*}
Since we averaged with the reversion, we have that the set $G_\varepsilon\defl \{ (\gamma,t)\in \overline\Fr(X):\ (f\circ\gamma)'(t)\ge (1-\varepsilon)|Df|_\ast(\gamma_t)|\gamma_t'| \}$ has positive $\bar\P_x$-measure for every $\varepsilon>0$, for $\mu$-a.e. $x\in \{|Df|_\ast>0\}$. The Weaver derivation $\bm b_\varepsilon$ \eqref{eq:derivation-from-plan} associated to $\P$ and $G=\chi_{G_\varepsilon}$ satisfies $\|\bm b_\varepsilon\|_{\mathcal X(\mu)}\le 1$ and $\bm b_\varepsilon(f)\ge (1-\varepsilon)|Df|_\ast$ $\mu$-a.e. on $\{|Df|_\ast>0\}$. Thus $h\ge \bm b_\varepsilon(f)\ge (1-\varepsilon)|Df|_\ast$. Since $\varepsilon>0$, the equality $h=|Df|_\ast$ $\mu$-a.e. follows. 

It remains to prove the estimate
\begin{align}\label{eq:derivation-minfrag-bound}
|\bm b(f)|\le |\bm b|_{\mathcal X,\loc}|Df|_\ast.
\end{align}
If we prove this inequality for all Lipschitz functions with bounded support, then the claim follows from locality of Weaver derivatives (see \cite[Lemma 2.120]{sch16}) and minimal weak-$\ast$ upper upper gradients (see e.g. Lemma \ref{lem:label-diff}). Thus, assume that $f\in \LIP_{bs}(X)$. If $(f_j)\subset \LIP_{bs}(X)$ is given by Theorem \ref{thm:approximation} then by the weak$^\ast$ continuity of $\bm b$ have
	\begin{align*}
		\int g\bm b(f)\ud\mu\le \liminf_{j\to\infty}\int g|\bm b|_{\mathcal X,\loc}\Lip_af_j\ud\mu=\int g|\bm b|_{\mathcal X,\loc}|Df|_\ast\ud \mu
	\end{align*}
for every $g\in L^1(\mu)$. From this \eqref{eq:derivation-minfrag-bound} easily follows.
\end{proof}

\begin{proof}[Proof of Theorem \ref{thm:weaver-der}]
Firstly, $\iota(\bm X)$ is a Weaver derivation for any $\bm X\in L^\infty(TX)$. Indeed, the product rule follows from the product rule of differentials $\ud (fg)=f\ud g+g\ud f$, see Lemma \ref{lem:label-diff}, and weak$^*$-continuity follows from the weak$^*$-continuity of differentials, cf. Theorem \ref{thm:diff-w-star-cont}. More precisely, suppose $(f_j)\subset \LIP_\infty(X)$ converges pointwise to $f\in\LIP(X)$ and $\sup_j\max\{\LIP(f_j),\|f_j\|_\infty\}<\infty$. For every $h\in L^1(\mu)$ we have 
\begin{align*}
	\lim_{j\to \infty}\int h\iota(X)(f_j)\ud \mu=\sum_m\lim_{j\to\infty}\int_{U_m}\ud f_j(h\bm X)\ud\mu=\int h\iota(\bm X)(f)\ud \mu
\end{align*}
by Theorem \ref{thm:diff-w-star-cont}, where $(U_j,\varphi_j)$ is a countable partition up to a $\mu$-null set of $X$ into fragment-wise charts such that $U_j$ is compact for each $j$.
Since 
\begin{align*}
	\iota(h\bm X+g\bm Y)(f)=\ud f(h\bm X+g\bm Y)=h\ud f(\bm X)+g\ud f(\bm Y)\quad \mu-a.e.
\end{align*}
for all $f\in \LIP(X)$ and $h,g\in L^\infty(\mu)$, it follows that $\iota$ is a homomorphism of $L^\infty(\mu)$-modules. 

We show that $\iota$ is 1-to-1. Suppose $\bm X\in \ker\iota$. If $\bm X\ne 0$, there exists $\delta_0>0$ and a fragment-wise chart $(U,\varphi)$ of dimension $n$ such that $|\bm X(x)|_x> \delta_0$ on $U$, and consequently a piecewise constant Borel map $\bm \xi:U\to (\R^n)^*$ with $|\bm \xi_x|_x^*=1$ such that $\bm \xi_x(\bm X(x))\ge \delta_0$ for $\mu$-a.e. $x\in U$ (cf. the proof of Proposition \ref{prop:indep-vs-Alberti-rep}). Thus, we find a positive measure subset $U'\subset U$ and $\xi\in(\R^n)^*$ such that $\bm\xi = \xi$ and $\xi(\bm X)\ge \delta_0$ on $U'$. However choosing $f=\xi\circ\varphi$ this implies that $\ud f(\bm X)=\xi(\bm X)\ne 0$ on $\mu$-a.e. $U'$, contradicting the fact that $\iota(\bm X)=0$. This shows that $\iota$ is injective.

Next we prove that $\iota$ is surjective. Let $\bm b$ be a Weaver derivation, and let $(U,\varphi)$ be a fragment-wise chart. Define the vector field $\bm X:X\to \R^n$, $\bm X=(\bm b(\varphi_1),\ldots,\bm b(\varphi_n))$. We claim that $\chi_U\iota(\bm X)=\chi_U\bm b$. By the locality of derivations (see \cite[Lemma 2.120]{sch16}), and since $X$ can be covered by countably many fragment-wise charts, this implies $\iota(\bm X)=\bm b$. Let $f\in\LIP(X)$, and let $\displaystyle \bm \xi_j=\sum_m\xi_j^m\chi_{V_j^m}:U\to (\R^n)^*$ be piecewise constant vector fields for each $j\in\N$ such that $\bm\xi_j\to \ud f$ converges $\mu$-a.e on $U$. By the linearity and locality of $\bm b$ we have that
\begin{align*}
\bm \xi_j(\bm X)=\sum_m\chi_{V_j^m}\xi_j^m((\bm b(\varphi_1),\ldots,\bm b(\varphi_n))=\sum_m\chi_{V_j^m}\bm b(\xi_j^m\circ\varphi)\quad \mu-a.e.\textrm{ on }U.
\end{align*}
Consequently we obtain 
\begin{align*}
	\chi_U|\bm\xi_j(f)-\bm b(f)|&\le \sum_m\chi_{V_j^m}|\bm b(\xi_j^m\circ\varphi-f)|\le \sum_m\chi_{V_j^m}|\bm b|_{\mathcal X, \loc}|D(\xi_j^m\circ\varphi-f)|_\ast\\
	&\le |\bm b|_{\mathcal X, \loc}\sum_m\chi_{V_j^m}|\xi_j^m-\ud f|^*\stackrel{j\to\infty}{\longrightarrow} 0\quad \mu-a.e.
\end{align*}
Since also $\bm\xi_j(\bm X)\to \ud f(\bm X)$ $\mu$-a.e. on $U$, we conclude that $\chi_U\bm b(f)=\chi_U\ud f(\bm X)$. 

It remains to prove \eqref{eq:isom}. Since $|\iota(\bm X)f|=|\ud f(\bm X)|\le |D f|_*|\bm X|$ for all $f\in \LIP(X)$, it follows that
\begin{align*}\label{eq:iota-bd-1}
	|\iota(\bm X)|_{\mathcal X, \loc}\le |\bm X|.
\end{align*}
On the other hand, by the linearity of $\iota(\bm X)\in \mathcal X(\mu)$, for every $\xi\in (\R^n)^*$ we have
\begin{align*}
	\xi(\bm X)=\xi(\iota(\bm X)(\varphi_1),\ldots,\iota(\bm X)(\varphi_n))=\iota(\bm X)(\xi\circ\varphi)\quad \mu-a.e.,
\end{align*}
and thus $\chi_U|\xi(\bm X)|\le \chi_U||\iota(\bm X)|_{\mathcal X,\loc}|D(\xi\circ\varphi)|_\ast=\chi_U||\iota(\bm X)|_{\mathcal X,\loc}|\xi|^*$ whenever $(U,\varphi)$ is a fragment-wise chart. It follows that
\begin{align*}
	|\bm X|\le |\iota(\bm X)|_{\mathcal X,\loc}\quad\mu-a.e.
\end{align*}
These two inequalities imply \eqref{eq:isom} and complete the proof.
\end{proof}

\subsection{Lipschitz differentiability spaces} We say that $X=(X,d,\mu)$ is a Lipschitz differentiability space if there is a countable collection $\{(U,\varphi)\}$, where $U\subset X$ is Borel with $\mu(U)>0$ and $\varphi\in \LIP(X,\R^n)$ ($n$ may depend on $U$), such that 
\begin{itemize}
    \item[(a)] the sets $U$ cover $X$ up to a $\mu$-null set, and 
    \item[(b)] every $f\in \LIP(X)$ admits a unique differential $\ud_xf\in (\R^n)^*$ with respect to $(U,\varphi)$ at a.e. $x\in U$ satisfying 
\[
\limsup_{y\to x}\frac{|f(y)-f(x)-\ud_xf(\varphi(y)-\varphi(x))|}{d(x,y)}=0.
\]
\end{itemize}
A pair $(U,\varphi)$ satisfying (b) is called a Lipschitz differentiability chart, or \emph{Cheeger} chart, and $\ud f$ is known as the Cheeger differential of $f$. See \cite{che99,kei02,bate15} for further background. The following lemma is a straightforward consequence of \cite[Theorem 1.15]{sch16b}. 
\begin{lemma}\label{lem:Lip-lsc} Let $X$ be a Lipschitz differentiability space. 
Suppose $(f_j)\subset \LIP(X)$ converges pointwise to $f\in \LIP(X)$ and $\sup_j\LIP(f_j)<\infty$. Then
\begin{align*}
	\int g\Lip f\ud\mu\le \liminf_{j\to\infty}\int g\Lip f_j\ud\mu
\end{align*}	
for every $g\in L^1(\mu)$.
\end{lemma}
\begin{proof}
By \cite[Theorem 1.15]{sch16b} we have that $\Lip f=\esssup\{\bm b(f):\|\bm b\|_{\mathcal X(\mu)}\le 1\}$ $\mu$-a.e. for every $f\in \LIP(X)$. Theorem \ref{thm:weaver-der} and the properties of the measure theoretic supremum further imply that there is a sequence $(\bm X_m)\subset L^\infty(TX)$ with $\|\bm X_m\|_{L^\infty(TX)}\le 1$ such that
\begin{align*}
	\Lip f(x)=\sup_{m\in \N}\ud_x f(\bm X_m)\quad \text{ for } \mu-\text{a.e. } x\in X.
\end{align*}
It follows that, for any $\varepsilon>0$, the sets $E_m=\{\ud f(\bm X_m)\ge (1-\varepsilon)\Lip f \}$ cover $X$ up to a $\mu$-null set. For each fixed $m$, the lower semicontinuity of differentials (Theorem \ref{thm:diff-w-star-cont}) implies that
\begin{align*}
	(1-\varepsilon)\int_{E_m}g\Lip f\ud\mu\le \int_{E_m}  g\ud f(\bm X_m)\ud\mu\le \liminf_{j\to\infty}\int_{E_m} g\Lip f_j\ud\mu
\end{align*}
for every non-negative $g\in L^1(\mu)$. Since the sets $E_m$ cover $X$ and $\varepsilon$ is arbitrary, the claim follows.

\end{proof}

\begin{proof}[Proof of Theorem \ref{thm:LDS-char}]
	Assume that (i) holds, i.e. that $X$ is a Lipschitz differentiability space. By \cite[Theorem 6.6]{bate15}, $X$ can be covered up to a $\mu$-null set by compact Cheeger charts $(U,\varphi)$ such that $\mu|_U$ admits $n_U$ $\varphi$-independent Alberti representations, and the dimension $n_U$ of $(U,\varphi)$ is the maximal number of independent Alberti representations of $\mu|_V$ for any $V\subset U$ with $\mu(V)>0$. Thus $(U,\varphi)$ is a fragment-wise chart. Moreover, for any $f\in \LIP(X)$, the fragment-wise and Cheeger differentials of $f$ with respect to $(U,\varphi)$ agree, cf. \cite[Lemma 3.8 or Corollary 6.7]{bate15}. To prove the implication (i)$\implies$ (ii) it suffices to show that $\Lip f=|Df|_\ast$ $\mu$-a.e. on $U$ for each such Cheeger chart $(U,\varphi)$.

	Since the claim is local, we may assume without loss of generality that $\spt f$ is bounded. Let $(f_j)\subset \LIP(X)$ be a sequence given by Theorem \ref{thm:approximation} with
	\begin{align*}
		\sup_j \LIP(f_j)<\infty
	\end{align*}
	converging pointwise to $f$ such that $\Lip f_j\to |Df|_\ast$. By Lemma \ref{lem:Lip-lsc} we have 
	\begin{align*}
		\int_U\Lip f\ud\mu\le \liminf_{j\to\infty}\int_U\Lip f_j\ud\mu=\int_U|Df|_\ast\ud \mu.
	\end{align*}
	Since $\Lip f$ is am $\ast$-upper gradient, we have by Proposition \ref{prop:min-star-ug-exists} the inequality $|Df|_\ast\le \Lip f$ $\mu$-a.e., this implies the claimed equality.
	
	The implication (ii)$\implies$(iii) is trivial, and thus it remains to prove (iii)$\implies$(i). Assume that (iii) holds. Applying (iii) to the function $f=\dist(S,\cdot)$ we obtain that $\Lip f(x)=0$ $\mu$-a.e. $x\in S$, since $|Df|_\ast \chi_S=0$ by Lemma \ref{lem:basic-prop-min-star-ug}. However, if $S$ is $\eta$-porous for some $\eta>0$ then $\Lip f(x)\ge \eta$ for every $x\in S$. Thus (iii) implies that porous sets are $\mu$-null. By \cite[Lemma 8.3]{bate15}  it follows that $X$ can be expressed as $X=N\cup \bigcup_iX_i$ where $\mu(N)=0$ and $\dim_HX_i<\infty$ for each $i$. Thus $X$ can be covered by a countable collection of fragment-wise charts, cf. Theorem \ref{thm:cover-by-frag-charts}. To complete the proof it suffices to show that every fragment-wise chart is a Cheeger chart.
	
	Let $(U,\varphi)$ be a fragment-wise chart (of dimension $n$) in $X$. In particular 
	\[
	0<I(\varphi,x)\le \inf_{\|\lambda\|_\ast=1}\Lip(\lambda\circ\varphi)(x)\quad\mu-a.e.\ x\in U,
	\]
	so that uniqueness of Cheeger differentials holds on $U$ (cf. \cite[Lemma 3.3]{bate15}). Let $f\in \LIP(X)$, and let $\ud f$ be the fragment-wise differential of $f$ with respect to $(U,\varphi)$. If $\Psi$ is given by Proposition \ref{prop:canonical rep} for $\psi=(\varphi,f)\in\LIP(X,\R^{n+1})$, then $\ud f:U\to (\R^n)^*$ satisfies $\Psi((\ud_xf,-1),x)=0$ $\mu$-a.e. $x\in U$. By our assumption (iii), we have 
    \begin{align*}
		\Lip(f- \bm\xi \circ\varphi)(x)\le \omega_x\left(\Psi((\bm\xi ,-1),x)\right)
	\end{align*}
 a.e. for every piecewise constant $\bm\xi$. By approximating $\ud_x f$ by such $\bm\xi$ we get
	\begin{align*}
		\Lip(f-\ud_x f\circ\varphi)(x)\le \omega_x\left(\Psi((\ud_xf,-1),x)\right)=0
	\end{align*}
	$\mu$-a.e. $x\in U$. Consequently $\ud f$ is the Cheeger differential of $f$ by definition.
\end{proof}

\bibliographystyle{plain}
\bibliography{abib}

@Preamble{"\def\cprime{$'$}"}

@misc{schulli,
      title={Characterizing rectifiability via biLipschitz pieces of Lipschitz mappings on the space}, 
      author={Sean Li and Raanan Schul},
      year={2025},
      eprint={2510.13525},
      archivePrefix={arXiv},
      primaryClass={math.MG},
      url={https://arxiv.org/abs/2510.13525}, 
}

@article{batepure,
 author = {Bate, David},
 title = {Purely unrectifiable metric spaces and perturbations of {Lipschitz} functions},
 fjournal = {Acta Mathematica},
 journal = {Acta Math.},
 issn = {0001-5962},
 volume = {224},
 number = {1},
 pages = {1--65},
 year = {2020},
 language = {English},
 doi = {10.4310/ACTA.2020.v224.n1.a1},
 zbMATH = {7194389},
 Zbl = {1440.28004}
}

@book {kechris,
    AUTHOR = {Kechris, A. S.},
     TITLE = {Classical descriptive set theory},
    SERIES = {Graduate Texts in Mathematics},
    VOLUME = {156},
 PUBLISHER = {Springer-Verlag, New York},
      YEAR = {1995},
     PAGES = {xviii+402},
      ISBN = {0-387-94374-9},
   MRCLASS = {03E15 (03-01 03-02 04A15 28A05 54H05 90D44)},
  MRNUMBER = {1321597},
MRREVIEWER = {Jakub Jasi\'{n}ski},
       DOI = {10.1007/978-1-4612-4190-4},
       URL = {https://doi.org/10.1007/978-1-4612-4190-4},
}

@article {sion58,
    AUTHOR = {Sion, M.},
     TITLE = {On general minimax theorems},
   JOURNAL = {Pacific J. Math.},
  FJOURNAL = {Pacific Journal of Mathematics},
    VOLUME = {8},
      YEAR = {1958},
     PAGES = {171--176},
      ISSN = {0030-8730,1945-5844},
   MRCLASS = {52.00},
  MRNUMBER = {97026},
MRREVIEWER = {Wendell\ H.\ Fleming},
       URL = {http://projecteuclid.org/euclid.pjm/1103040253},
}

@article {car-jar10,
    AUTHOR = {Durand-Cartagena, E. and Jaramillo, J. A.},
     TITLE = {Pointwise {L}ipschitz functions on metric spaces},
   JOURNAL = {J. Math. Anal. Appl.},
  FJOURNAL = {Journal of Mathematical Analysis and Applications},
    VOLUME = {363},
      YEAR = {2010},
    NUMBER = {2},
     PAGES = {525--548},
      ISSN = {0022-247X,1096-0813},
   MRCLASS = {46E35 (30L99)},
  MRNUMBER = {2564873},
MRREVIEWER = {Jeremy\ T.\ Tyson},
       DOI = {10.1016/j.jmaa.2009.09.039},
       URL = {https://doi.org/10.1016/j.jmaa.2009.09.039},
}

@article {bate15,
    AUTHOR = {Bate, D.},
     TITLE = {Structure of measures in {L}ipschitz differentiability spaces},
   JOURNAL = {J. Amer. Math. Soc.},
  FJOURNAL = {Journal of the American Mathematical Society},
    VOLUME = {28},
      YEAR = {2015},
    NUMBER = {2},
     PAGES = {421--482},
      ISSN = {0894-0347,1088-6834},
   MRCLASS = {46G05 (30L99 49J52 53C23)},
  MRNUMBER = {3300699},
MRREVIEWER = {Riikka\ Korte},
       DOI = {10.1090/S0894-0347-2014-00810-9},
       URL = {https://doi.org/10.1090/S0894-0347-2014-00810-9},
}

@book{dellacherie2011probabilities,
    AUTHOR = {Dellacherie, C. and Meyer, P.-A.},
     TITLE = {Probabilities and potential. {C}},
    SERIES = {North-Holland Mathematics Studies},
    VOLUME = {151},
      NOTE = {Potential theory for discrete and continuous semigroups,
              Translated from the French by J. Norris},
 PUBLISHER = {North-Holland Publishing Co., Amsterdam},
      YEAR = {1988},
     PAGES = {xiv+416},
      ISBN = {0-444-70386-1},
   MRCLASS = {60-02 (31-02 31C20 60J45)},
  MRNUMBER = {939365},
}

@article {che-kle-sch16,
    AUTHOR = {Cheeger, J. and Kleiner, B. and Schioppa, A.},
     TITLE = {Infinitesimal structure of differentiability spaces, and
              metric differentiation},
   JOURNAL = {Anal. Geom. Metr. Spaces},
  FJOURNAL = {Analysis and Geometry in Metric Spaces},
    VOLUME = {4},
      YEAR = {2016},
    NUMBER = {1},
     PAGES = {104--159},
      ISSN = {2299-3274},
   MRCLASS = {58C20 (28A15 30L05)},
  MRNUMBER = {3543676},
MRREVIEWER = {Riikka\ Korte},
       DOI = {10.1515/agms-2016-0005},
       URL = {https://doi.org/10.1515/agms-2016-0005},
}

@article{teriseb,
	author={Eriksson-Bique, S. and Soultanis, E.},
	title={Curvewise characterizations of minimal upper gradients and the construction of a {S}obolev differential},
journal = {To appear in {A}nalysis and PDE},
	note = "arXiv:2102.08097",
	year={2021}
}

@article{bate18,
title = {Differentiability and Poincaré-type inequalities in metric measure spaces},
journal = {Advances in Mathematics},
volume = {333},
pages = {868-930},
year = {2018},
issn = {0001-8708},
doi = {https://doi.org/10.1016/j.aim.2018.06.002},
url = {https://www.sciencedirect.com/science/article/pii/S0001870818302263},
author = {Bate, D. and Li, S.}
}

@book {rud,
    AUTHOR = {Rudin, W.},
     TITLE = {Function theory in the unit ball of {$\Bbb C^n$}},
    SERIES = {Classics in Mathematics},
      NOTE = {Reprint of the 1980 edition},
 PUBLISHER = {Springer-Verlag, Berlin},
      YEAR = {2008},
     PAGES = {xiv+436},
      ISBN = {978-3-540-68272-1},
   MRCLASS = {32-02 (32Axx 46E15 46E20)},
  MRNUMBER = {2446682},
}

@article{kell,
 author = {Kell, Martin},
 title = {On {Cheeger} and {Sobolev} differentials in metric measure spaces},
 fjournal = {Revista Matem{\'a}tica Iberoamericana},
 journal = {Rev. Mat. Iberoam.},
 issn = {0213-2230},
 volume = {35},
 number = {7},
 pages = {2119--2150},
 year = {2019},
 language = {English},
 doi = {10.4171/rmi/1114},
 zbMATH = {7144996},
 Zbl = {1433.53069}
}

@article{pasqualetto20,
 author = {Di Marino, Simone and Lu{\v{c}}i{\'c}, Danka and Pasqualetto, Enrico},
 title = {A short proof of the infinitesimal hilbertianity of the weighted {Euclidean} space},
 fjournal = {Comptes Rendus. Math{\'e}matique. Acad{\'e}mie des Sciences, Paris},
 journal = {C. R., Math., Acad. Sci. Paris},
 issn = {1631-073X},
 volume = {358},
 number = {7},
 pages = {817--825},
 year = {2020},
 language = {English},
 doi = {10.5802/crmath.88},
 zbMATH = {7278835},
 Zbl = {1466.46025}
}

@article{keith04,
title = "A differentiable structure for metric measure spaces",
journal = "Advances in Mathematics",
volume = "183",
number = "2",
pages = "271 - 315",
year = "2004",
issn = "0001-8708",
doi = "https://doi.org/10.1016/S0001-8708(03)00089-6",
url = "http://www.sciencedirect.com/science/article/pii/S0001870803000896",
author = "Keith, S.",
}

@article {seb2020,
    AUTHOR = {Eriksson-Bique, S.},
     TITLE = {Density of {L}ipschitz functions in energy},
   JOURNAL = {Calc. Var. Partial Differential Equations},
  FJOURNAL = {Calculus of Variations and Partial Differential Equations},
    VOLUME = {62},
      YEAR = {2023},
    NUMBER = {2},
     PAGES = {Paper No. 60, 23},
      ISSN = {0944-2669,1432-0835},
   MRCLASS = {46E36 (26B30 30L99)},
  MRNUMBER = {4525741},
       DOI = {10.1007/s00526-022-02395-1},
       URL = {https://doi.org/10.1007/s00526-022-02395-1},
}

@article{amb13,
  title={On the duality between p-modulus and probability measures},
  author={Ambrosio, L. and Di Marino, S. and Savar{\'e}, G.},
  journal={J. Eur. Math. Soc.},
  volume={17}, 
  pages={1817--1853},
  year={2015}
}

@article {cheegerkleiner,
    AUTHOR = {Cheeger, J. and Kleiner, B.},
     TITLE = {Differentiability of {L}ipschitz maps from metric measure
              spaces to {B}anach spaces with the {R}adon-{N}ikod\'{y}m property},
   JOURNAL = {Geom. Funct. Anal.},
  FJOURNAL = {Geometric and Functional Analysis},
    VOLUME = {19},
      YEAR = {2009},
    NUMBER = {4},
     PAGES = {1017--1028},
      ISSN = {1016-443X},
   MRCLASS = {30L05 (46G05 58C20)},
  MRNUMBER = {2570313},
MRREVIEWER = {Jeremy T. Tyson},
       DOI = {10.1007/s00039-009-0030-6},
       URL = {https://doi.org/10.1007/s00039-009-0030-6},
}

@misc{almarflat,
 author = {Alberti, Giovanni and Marchese, Andrea},
 title = {On the structure of flat chains with finite mass},
 year = {2023},
 howpublished = {Preprint, {arXiv}:2311.06099 [math.{CA}] (2023)},
 url = {https://arxiv.org/abs/2311.06099},
 arXiv = {arXiv:2311.06099}
}

@article {almar16,
    AUTHOR = {Alberti, G. and Marchese, A.},
     TITLE = {On the differentiability of {L}ipschitz functions with respect
              to measures in the {E}uclidean space},
   JOURNAL = {Geom. Funct. Anal.},
  FJOURNAL = {Geometric and Functional Analysis},
    VOLUME = {26},
      YEAR = {2016},
    NUMBER = {1},
     PAGES = {1--66},
      ISSN = {1016-443X},
   MRCLASS = {49Q15 (26A27 26B05 28A75)},
  MRNUMBER = {3494485},
MRREVIEWER = {Gareth Speight},
       DOI = {10.1007/s00039-016-0354-y},
       URL = {https://doi.org/10.1007/s00039-016-0354-y},
}

@article{shanmunsemmes,
  title={Equivalence of two BV classes of functions in metric spaces, and existence of a Semmes family of curves under a 1-{P}oincar{\'e} inequality},
  author={Durand-Cartagena, E. and Eriksson-Bique, S. and Korte, R. and Shanmugalingam, N.},
  journal={Advances in Calculus of Variations},
  volume={1},
  number={ahead-of-print},
  year={2019},
  publisher={De Gruyter}
}

@article{exnerova2019plans,
title = {Plans on measures and AM-modulus},
journal = {Journal of Functional Analysis},
volume = {281},
number = {10},
pages = {109205},
year = {2021},
issn = {0022-1236},
doi = {https://doi.org/10.1016/j.jfa.2021.109205},
url = {https://www.sciencedirect.com/science/article/pii/S0022123621002871},
author = {Vendula {Honzlová Exnerová} and Ondřej F.K. Kalenda and Jan Malý and Olli Martio}
}

@Unpublished{bkt19,
author = {Bate, D. and Kangasniemi, I. and Orponen, T.},
title = {Cheeger's differentiation theorem via the multilinear {K}akeya inequality},
journal = {preprint (arXiv:1904.00808)},
note = {To appear, Pure Appl. Funct. Anal.},
year = {2019},
}

@book{AGS08,
	Address = {Basel},
	Author = {Ambrosio, L. and Gigli, N. and Savar{\'e}, G.},
	Edition = {Second},
	Isbn = {978-3-7643-8721-1},
	Mrclass = {49-02 (28A33 35K55 35K90 49Q20 60B05)},
	Mrnumber = {2401600 (2009h:49002)},
	Mrreviewer = {Pietro Celada},
	Pages = {x+334},
	Publisher = {Birkh\"auser Verlag},
	Series = {Lectures in Mathematics ETH Z\"urich},
	Title = {Gradient flows in metric spaces and in the space of probability measures},
	Year = {2008}}

@article{keith03,
author={Keith, S.},
year={2003},
title={Modulus and the {P}oincar\'e inequality on metric measure spaces}, 
journal={Mathematische Zeitschrift},
volume={245},
number={2}, 
pages={255--292}
}

@article{davideb20,
author={David, G. C.   and Eriksson-Bique,  S.},
title={Infinitesimal splitting for spaces with thick curve families and Euclidean embeddings},
journal={To appear in Annales de l’institut Fourier (arXiv:2006.10668)},
year={2020},
}

@article {gig15,
	AUTHOR = {Gigli, N.},
	TITLE = {On the differential structure of metric measure spaces and
	applications},
	JOURNAL = {Mem. Amer. Math. Soc.},
	FJOURNAL = {Memoirs of the American Mathematical Society},
	VOLUME = {236},
	YEAR = {2015},
	NUMBER = {1113},
	PAGES = {vi+91},
	ISSN = {0065-9266},
	ISBN = {978-1-4704-1420-7},
	MRCLASS = {53C23 (30L05 30L10 49Q15 58C20)},
	MRNUMBER = {3381131},
	MRREVIEWER = {Davide Vittone},
	DOI = {10.1090/memo/1113},
	URL = {https://doi.org/10.1090/memo/1113},
}

@article {gig18,
	AUTHOR = {Gigli, N.},
	TITLE = {Nonsmooth differential geometry---an approach tailored for
	spaces with {R}icci curvature bounded from below},
	JOURNAL = {Mem. Amer. Math. Soc.},
	FJOURNAL = {Memoirs of the American Mathematical Society},
	VOLUME = {251},
	YEAR = {2018},
	NUMBER = {1196},
	PAGES = {v+161},
	ISSN = {0065-9266},
	ISBN = {978-1-4704-2765-8; 978-1-4704-4266-8},
	MRCLASS = {49J52 (46E25 53C23)},
	MRNUMBER = {3756920},
	DOI = {10.1090/memo/1196},
	URL = {https://doi.org/10.1090/memo/1196},
}

@article{hino,
 author = {Hino, Masanori},
 title = {Energy measures and indices of {Dirichlet} forms, with applications to derivatives on some fractals},
 fjournal = {Proceedings of the London Mathematical Society. Third Series},
 journal = {Proc. Lond. Math. Soc. (3)},
 issn = {0024-6115},
 volume = {100},
 number = {1},
 pages = {269--302},
 year = {2010},
 language = {English},
 doi = {10.1112/plms/pdp032},
 zbMATH = {5658780},
 Zbl = {1185.60089}
}

@misc{AIDP2,
 author = {Luigi Ambrosio and Toni Ikonen and Danka Lu{\v{c}}i{\'c} and Enrico Pasqualetto},
 title = {Metric {Sobolev} spaces {II}: dual energies and divergence measures},
 year = {2025},
 howpublished = {Preprint, {arXiv}:2510.12424 [math.{FA}] (2025)},
 url = {https://arxiv.org/abs/2510.12424},
 arXiv = {arXiv:2510.12424}
}

@article{AIDP1,
 author = {Ambrosio, Luigi and Ikonen, Toni and Lu{\v{c}}i{\'c}, Danka and Pasqualetto, Enrico},
 title = {Metric {Sobolev} spaces. {I}: {Equivalence} of definitions},
 fjournal = {Milan Journal of Mathematics},
 journal = {Milan J. Math.},
 issn = {1424-9286},
 volume = {92},
 number = {2},
 pages = {255--347},
 year = {2024},
 language = {English},
 doi = {10.1007/s00032-024-00407-7},
 zbMATH = {7960925},
 Zbl = {1565.46025}
}

@article{che99,
	author = "Cheeger, J.",
	coden = "GFANFB",
	doi = "10.1007/s000390050094",
	fjournal = "Geometric and Functional Analysis",
	issn = "1016-443X",
	journal = "Geom. Funct. Anal.",
	mrclass = "53C23 (49J52)",
	mrnumber = "1708448 (2000g:53043)",
	mrreviewer = "William P. Minicozzi, II",
	number = "3",
	pages = "428--517",
	title = "{Differentiability of {L}ipschitz functions on metric measure spaces}",
	url = "http://dx.doi.org/10.1007/s000390050094",
	volume = "9",
	year = "1999"
}

@Book{ambrosiotilli,
 Author = {Ambrosio, L. and T., Paolo},
 Title = {Topics on analysis in metric spaces},
 FSeries = {Oxford Lecture Series in Mathematics and its Applications},
 Series = {Oxf. Lect. Ser. Math. Appl.},
 Volume = {25},
 ISBN = {0-19-852938-4},
 Year = {2004},
 Publisher = {Oxford: Oxford University Press},
 Language = {English},
 zbMATH = {2024202},
 Zbl = {1080.28001}
}

@article{heinonen2007nonsmooth,
  title={Nonsmooth calculus},
  author={Heinonen, J.},
  journal={Bulletin of the American mathematical society},
  volume={44},
  number={2},
  pages={163--232},
  year={2007}
}

@book{gigenr,
 Author = {Gigli, N. and Pasqualetto, E.},
 Title = {Lectures on nonsmooth differential geometry},
 FSeries = {SISSA Springer Series},
 Series = {SISSA Springer Ser.},
 ISSN = {2524-857X},
 Volume = {2},
 ISBN = {978-3-030-38612-2; 978-3-030-38615-3; 978-3-030-38613-9},
 Year = {2020},
 Publisher = {Cham: Springer},
 Language = {English},
 DOI = {10.1007/978-3-030-38613-9},
 zbMATH = {7183685},
 Zbl = {1452.53002}
}

@book{hei01,
	address = "New York",
	author = "Heinonen, J.",
	isbn = "0-387-95104-0",
	mrclass = "30C65 (28A75 28A78 46E35)",
	mrnumber = "1800917 (2002c:30028)",
	mrreviewer = "Christopher Bishop",
	pages = "x+140",
	publisher = "Springer-Verlag",
	series = "{Universitext}",
	title = "{Lectures on analysis on metric spaces}",
	year = "2001"
}

@book{kei02,
	author = "Keith, S. J.",
	isbn = "978-0493-55682-6",
	mrclass = "Thesis",
	mrnumber = "2703118",
	note = "Thesis (Ph.D.)--University of Michigan",
	pages = "80",
	publisher = "ProQuest LLC, Ann Arbor, MI",
	title = "{A differentiable structure for metric measure spaces}",
	url = "http://gateway.proquest.com/openurl?url_ver=Z39.88-2004&rft_val_fmt=info:ofi/fmt:kev:mtx:dissertation&res_dat=xri:pqdiss&rft_dat=xri:pqdiss:3042096",
	year = "2002"
}

@article{sha00,
	author = "Shanmugalingam, N.",
	doi = "10.4171/RMI/275",
	fjournal = "Revista Matem{\'a}tica Iberoamericana",
	issn = "0213-2230",
	journal = "Rev. Mat. Iberoamericana",
	mrclass = "46E35",
	mrnumber = "1809341 (2002b:46059)",
	mrreviewer = "Daniele Morbidelli",
	number = "2",
	pages = "243--279",
	title = "{Newtonian spaces: an extension of {S}obolev spaces to metric measure spaces}",
	url = "http://dx.doi.org/10.4171/RMI/275",
	volume = "16",
	year = "2000"
}

@book{HKST07,
	address = "United Kingdom",
	author = "Heinonen, J. and Koskela, P. and Shanmugalingam, N. and Tyson, J.",
	edition = "First",
	isbn = "9781107092341",
	mrclass = "",
	mrnumber = "",
	pages = "448",
	publisher = "Cambridge University Press",
	series = "{New Mathematical Monographs}",
	title = "{Sobolev spaces on metric measure spaces: an approach based on upper gradients}",
	year = "2015"
}

@book{bjo11,
	author = "Bj{\"o}rn, A. and Bj{\"o}rn, J.",
	doi = "10.4171/099",
	isbn = "978-3-03719-099-9",
	mrclass = "31-02 (31C45 31E05 35B45 35D30 35J70 35J92)",
	mrnumber = "2867756",
	mrreviewer = "Nageswari Shanmugalingam",
	pages = "xii+403",
	publisher = "European Mathematical Society (EMS), Z{\"u}rich",
	series = "{EMS Tracts in Mathematics}",
	title = "{Nonlinear potential theory on metric spaces}",
	url = "http://dx.doi.org/10.4171/099",
	volume = "17",
	year = "2011"
}

@article{bateclosability,
 author = {Alberti, Giovanni and Bate, David and Marchese, Andrea},
 title = {On the closability of differential operators},
 fjournal = {Journal of Functional Analysis},
 journal = {J. Funct. Anal.},
 issn = {0022-1236},
 volume = {289},
 number = {7},
 pages = {22},
 note = {Id/No 111029},
 year = {2025},
 language = {English},
 doi = {10.1016/j.jfa.2025.111029},
 zbMATH = {8053122}
}

@book{wea99,
	address = "River Edge, NJ",
	author = "Weaver, N.",
	isbn = "981-02-3873-8",
	mrclass = "46-02 (46Bxx 46Exx 46H05 46J10)",
	mrnumber = "1832645 (2002g:46002)",
	mrreviewer = "S. J. Sidney",
	pages = "xiv+223",
	publisher = "World Scientific Publishing Co. Inc.",
	title = "{Lipschitz algebras}",
	year = "1999"
}

@article{hei98,
	author = "Heinonen, J. and Koskela, P.",
	coden = "ACMAA8",
	doi = "10.1007/BF02392747",
	fjournal = "Acta Mathematica",
	issn = "0001-5962",
	journal = "Acta Math.",
	mrclass = "30C65 (46E99)",
	mrnumber = "1654771 (99j:30025)",
	mrreviewer = "M. Yu. Vasil{\cprime}chik",
	number = "1",
	pages = "1--61",
	title = "{Quasiconformal maps in metric spaces with controlled geometry}",
	url = "http://dx.doi.org/10.1007/BF02392747",
	volume = "181",
	year = "1998"
}

@misc{eidconj,
 author = {Eriksson-Bique, Sylvester and Murugan, Mathav},
 title = {On the energy image density conjecture of {Bouleau} and {Hirsch}},
 year = {2025},
 howpublished = {Preprint, {arXiv}:2510.13659 [math.{PR}] (2025)},
 url = {https://arxiv.org/abs/2510.13659},
 arXiv = {arXiv:2510.13659}
}

@article{marchesemerlo,
 author = {Marchese, Andrea and Merlo, Andrea},
 title = {Characterization of rectifiability via {Lusin}-type approximation},
 fjournal = {Analysis \& PDE},
 journal = {Anal. PDE},
 issn = {2157-5045},
 volume = {17},
 number = {6},
 pages = {2109--2121},
 year = {2024},
 language = {English},
 doi = {10.2140/apde.2024.17.2109},
 zbMATH = {7904879}
}

@article{wil12,
	author = "Williams, Marshall",
	fjournal = "Annali della Scuola Normale Superiore di Pisa. Classe di Scienze. Serie V",
	issn = "0391-173X",
	journal = "Ann. Sc. Norm. Super. Pisa Cl. Sci. (5)",
	mrclass = "58A25 (30L10 43A05 49Q15)",
	mrnumber = "3011992",
	mrreviewer = "David Matthew Freeman",
	number = "2",
	pages = "259--302",
	title = "{Metric currents, differentiable structures, and {C}arnot groups}",
	volume = "11",
	year = "2012",
}

@article{sch16b,
  title={Derivations and Alberti representations},
  author={Schioppa, A.},
  journal={Advances in Mathematics},
  volume={293},
  pages={436--528},
  year={2016},
  publisher={Elsevier}
}

@article{sch16,
	author = "Schioppa, A.",
	doi = "10.1016/j.jfa.2016.08.022",
	fjournal = "Journal of Functional Analysis",
	issn = "0022-1236",
	journal = "J. Funct. Anal.",
	mrclass = "53C23 (49Q15)",
	mrnumber = "3554700",
	mrreviewer = "Enrico Le Donne",
	number = "11",
	pages = "3007--3081",
	title = "{Metric currents and {A}lberti representations}",
	url = "https://doi.org/10.1016/j.jfa.2016.08.022",
	volume = "271",
	year = "2016",
}
\end{document}